\documentclass[letterpaper,12pt,reqno]{amsart}
\usepackage[margin=1.5in]{geometry}

\usepackage[dvips]{color}
\usepackage{graphicx}
\usepackage{epsfig}
\usepackage{tikz}
\usepackage{tikz-cd}
\usepackage{enumerate}
\usepackage{enumitem}
\usepackage{color}
\usepackage{hyperref}
\hypersetup{
     colorlinks   = true,
     citecolor    = black,
     linkcolor    = blue,
     urlcolor     = blue
}

\usepackage{latexsym}
\usepackage{amssymb}

\usepackage{amsmath}
\usepackage{amsthm}
\usepackage{amscd}
\usepackage{latexsym}
\usepackage{amssymb}
\usepackage{amsmath}
\usepackage{amsthm}
\usepackage{cite}
\usepackage{amscd}

\theoremstyle{plain}
\newtheorem{thm}{Theorem}[section]
\newtheorem{lem}[thm]{Lemma}

\newtheorem{prop}[thm]{Proposition}
\theoremstyle{definition}
\newtheorem{cor}[thm]{Corollary}
\newtheorem{rem}[thm]{Remark}
\newtheorem{claim}[thm]{Claim}

\newtheorem{ques}[thm]{Question}

\newtheorem*{acknowledgement}{Acknowledgements}

\newcommand{\N}{\mathbb{N}}
\newcommand{\Z}{\mathbb{Z}}
\newcommand{\Q}{\mathbb{Q}}
\newcommand{\R}{\mathbb{R}}

\newcommand{\mc}{\mathcal}

\newcommand{\mf}{\mathfrak}
\newcommand{\mbf}{\mathbf}
\newcommand{\mrm}{\mathrm}

\renewcommand{\a}{\alpha}
\renewcommand{\b}{\beta}
\newcommand{\g}{\gamma}
\newcommand{\G}{\Gamma}
\renewcommand{\d}{\delta}
\newcommand{\e}{\varepsilon}
\renewcommand{\l}{\lambda}

\newcommand{\w}{\omega}
\newcommand{\s}{\sigma}
\newcommand{\vp}{\varphi}
\renewcommand{\t}{\tau}
\renewcommand{\th}{\theta}

\renewcommand{\k}{\kappa}

\newcommand{\set}[1]{\left\{#1\right\}}
\renewcommand{\r}{\rightarrow}

\newcommand{\norm}[1]{\left\lVert#1\right\rVert}

\newcommand{\Lcal}{\mc{L}}
\newcommand{\Kcal}{\mc{K}}

\newcommand{\Pcal}{\mc{P}}
\newcommand{\Wcal}{\mc{W}}

\newcommand{\Dcal}{\mc{D}}

\newcommand{\Ocal}{\mc{O}}
\newcommand{\Qcal}{\mc{Q}}

\newcommand{\Ucal}{\mc{U}}
\newcommand{\Vcal}{\mc{V}}

\newcommand{\GLd}{\mrm{GL}_2(\R)}

\newcommand{\RP}{\R\mathbb{P}^1}

\newcommand{\KZ}[2]{\mrm{KZ}(#1,#2)}
\newcommand{\inj}{\mrm{inj}}
\newcommand{\Wu}{W^{\mrm{u}}_{\mrm{loc}}}
\newcommand{\Wcs}{W^{\mrm{cs}}_{\mrm{loc}}}

\newcommand{\tO}{{\widetilde{\mathcal O}}}
\newcommand{\GLp}{\mrm{GL}_2^+(\R)}

\newcommand\hol{\mathrm{hol}}
\newcommand{\dist}{{\rm dist}}

\newcommand{\HH}{{\mathcal{H}}}

\newcommand{\HHH}{{\mathbb{H}}}

\newcommand{\TT}{{\mathcal{T}}}
\newcommand{\Dom}{{\mathrm{Dom}}}

\newcommand{\future}[1]{{}}

\newcommand{\bal}[1]{{\mathrm T}^{\mrm{bal}}_{#1}}
\newcommand{\can}[1]{{\mathrm T}^{\mrm{st}}_{#1}}
\newcommand{\cyl}{\mathfrak{C}}

\newcommand{\SL}{\mathrm{SL}_2(\R)}
\newcommand{\tq}{\tilde{q}}

\newcommand{\tK}{\tilde{\Kcal}}
\newcommand{\Mod}{\mrm{Mod}(S,\Sigma)}
\newcommand{\Hom}{\mrm{Hom}}

\newcommand{\tVV}{\tilde{\Vcal}}
\newcommand{\VV}{\Vcal}
\newcommand{\HHu}{\HH_{\mrm{u}}}
\newcommand{\HHm}{\HH_{\mrm{m}}}
\newcommand{\tp}{\tilde{p}}
\newcommand{\tw}{\tilde{\w}}
\newcommand{\distu}{\dist_{\mrm{u}}}
\newcommand{\distm}{\dist_{\mrm{m}}}

\numberwithin{equation}{section}

\title[Failure of Mozes-Shah Phenomenon on Moduli Spaces]{On the Space of Ergodic Measures for the Horocycle Flow on Strata of Abelian Differentials}

\author{Jon Chaika}
\address{Department of Mathematics, University of Utah, Salt Lake City, UT.}
\email{chaika@math.utah.edu}

\author{Osama Khalil}
\address{Department of Mathematics, Statistics, and Computer Science, University of Illinois Chicago, Chicago, IL.}
\email{okhalil@uic.edu}
 
\author{John Smillie}
\address{Mathematics Institute, University of Warwick, Coventry CV4 7AL, U.K.}
\email{J.Smillie@warwick.ac.uk}

\begin{document}

\maketitle

    \begin{abstract}
We study the  horocycle flow
on the stratum of abelian differentials 
$\HH(2)$. We show that there is a sequence of horocycle ergodic measures, each supported on a periodic horocycle orbit, which weakly converges to an invariant, but non-ergodic, measure by $\SL$.
As a consequence, we show that there are points in $\HH(2)$ whose horocycle flow orbits do not equidistribute towards any invariant measure.
\end{abstract}

\section{Introduction}
\subsection{Context and main results}
A central topic in the areas of ergodic theory and geometry of surfaces is that of the dynamics of the $\SL$-action on  moduli spaces of abelian differentials. 
Many of the results in this area have been inspired by the analogy with the theory of actions of Lie groups on homogeneous spaces. 
The works~\cite{EM,EMM} establish strong rigidity results for the action of $\SL$ and its full upper triangular parabolic subgroup on moduli spaces which mirror Ratner's fundamental results in homogeneous dynamics; cf.~\cite{ratner}. 
Ratner's orbit and measure classification theorems for unipotent flows on homogeneous spaces  show that these enjoy strong rigidity properties.
On the other hand, dynamics of horocycle flows on moduli spaces given by the action of the subgroup
$U:=\set{u_s=\left(\begin{smallmatrix} 1 & s\\0&1\end{smallmatrix}\right):s\in\R }$ of $ \SL$
remain largely mysterious. There are some positive results in relatively simple settings cf.~\cite{EskinMarklofMorris,EskinMasurSchmoll,BSW} for measures and orbit closures and striking non-rigidity results in slightly more complicated settings for orbit closures~\cite{CSW}.

Ratner's results in homogeneous dynamics have found numerous deep applications. Among the key tools that enable such applications is the work of Mozes-Shah~\cite{MozSha} regarding limits of unipotent invariant measures, which builds in an essential way on the results of Dani-Margulis~\cite{DaniMargulis}. The goal of this article is to explore the validity of the analogs of these fundamental results for the horocycle flow $u_s$ on moduli spaces. We begin by stating our main results and defer the detailed definitions to Section~\ref{sec:background}.

The moduli space of abelian differentials is a union of strata. We consider the horocycle flow on the stratum $\HH_1(2)$ consisting of unit area translation surfaces of genus $2$ with one cone point.
This $7$-dimensional space is an affine orbifold admitting an action by $\SL$ and is the support of a $U$-ergodic and invariant probability measure of Lebesgue class known as the Masur-Veech measure and denoted $\mu_{\mrm{MV}}$.

\begin{thm} \label{thm:intro}

There exists a sequence of $U$-invariant and ergodic probability measures $\nu_n$ on $\HH_1(2)$ such that
\begin{enumerate}
    \item For every $n\in\N$, $\nu_n$ is the uniform measure on a $U$-periodic orbit.
    \item The sequence $(\nu_n)$ weak-$\ast$ converges to a measure $\nu$ which is a non-trivial convex combination of $\mu_{\mrm{MV}}$ and the $\SL$-invariant measure on a Teichm\"uller curve in $\HH_1(2)$.
\end{enumerate}
In particular, $\nu$ is $\SL$-invariant but is not ergodic for the action of either $U$ or $\SL$.
\end{thm}

Let $V$ be a one-parameter unipotent flow on a homogeneous space. The result of Mozes-Shah asserts that every limit point of a sequence of $V$-ergodic measures is ergodic for a possibly larger group, which is also generated by unipotents \cite{MozSha}. One potential larger group in our setting would be $\SL$ but the theorem shows that the limit measure is not ergodic for the action of this group.
We will show in Remark~\ref{rem:DM} that no larger group of any sort  can appear thus Theorem~\ref{thm:intro} shows that the analog of Mozes-Shah does not hold for moduli spaces.
The key ingredient in our proof is a statement showing roughly that the analog of the fundamental results of Dani-Margulis~\cite{DaniMargulis} ceases to hold for the $U$-action on $\HH(2)$; cf.~Theorem~\ref{thm:sufficient} and the discussion in Section~\ref{sec:motivations1}.

In \cite{CSW} it is shown that horocycle {\em orbit closures} can  have dramatically differently behavior from unipotent orbit closures in homogeneous dynamics but Theorem~\ref{thm:intro} is the first result asserting that horocycle {\em measures} display different behavior from those in homogeneous dynamics.
 
 Theorem~\ref{thm:intro} has additional dynamical consequences.
Ratner showed that  if $V$ is any one-parameter unipotent flow on a (finite volume) homogeneous space then the orbit of every point equidistributes towards some $V$-ergodic measure \cite{ratner}. This property does not hold for $U$ acting on $\HH_1(2)$.

\begin{cor}\label{cor:nongen} 
There
exists a dense $\mrm{G}_\d$ subset $S\subset\HH_1(2)$, so that for every $x \in S$, there is a function $f\in C_c\big(\HH_1(2)\big)$, such that the limit of $\frac 1 T \int_0^T f(u_tx)\;dt$ does not exist as $T \to \infty$. In particular, every $x\in S$ fails to be generic with respect to any $U$-invariant measure.
\end{cor}

A non-genericity result analogous to Corollary~\ref{cor:nongen} for the 9-dimensional stratum $\HH_1(1,1)$ was established in~\cite{CSW} by fundamentally different methods. That work relied on the existence of surfaces with minimal non-uniquely ergodic horizontal foliations, while no such surfaces exist in $\HH(2)$.

\subsection{Comparison with homogeneous dynamics}\label{sec:motivations1}

As discussed above, Theorem~\ref{thm:intro} and Corollary~\ref{cor:nongen} are in contrast with the results of~\cite{MozSha} and~\cite{ratner}. It is worth noting that, in general, the set of $V$-ergodic measures need not be closed in the set of $V$-invariant probability measures for any one-parameter unipotent group $V$. However, it follows from powerful results of Dani and Margulis, that if a sequence of periodic orbits becomes dense within the support of a $V$-ergodic measure $\mu$, then the sequence of ergodic measures supported on these periodic orbits weak-$*$ converges to $\mu$ (see~\cite{DaniMargulis}).
This is in contrast to the behavior of the sequence of periodic $U$-orbits described in Theorem~\ref{thm:intro} where the orbits become dense in the support of $\mu_{\mrm{MV}}$ but the uniform measures $\nu_n$ on these orbits do not converge to $\mu_{\mrm{MV}}$.

\begin{rem} \label{rem:DM} We have observed that the limit measure $\nu$ in Theorem~\ref{thm:intro} is not ergodic for the $\SL$-action.
We remark here that this measure fails to  be ergodic for the action of \textit{any} group acting by locally bilipschitz maps on $\HH_1(2)$.
Recall that $\HH_1(2)$ has dimension $7$ while a Teichm\"uller curve has dimension $3$ so the set
$$\left\{ x\in \HH_1(2):\underset{r \to 0}{\limsup}\, \frac{\nu(B(x,r))}{r^7}<\infty\right\}$$ is invariant under locally bilipschitz diffeomorphisms that preserve $\nu$, but it does not have measure $0$ or $1$ with respect to $\nu$.  
\end{rem}

\subsection{Comparison with the $P$-action}\label{sec:motivations2} 
As mentioned above, Theorem~\ref{thm:intro} and Corollary~\ref{cor:nongen} are in contrast with the rigid behavior of the action of the upper triangular group $P$ on moduli spaces established in~\cite{EMM}.
Indeed,
in spaces of translation surfaces every point equidistributes under $P$ (which is an amenable group) to a $P$-ergodic probability measure and the set of $P$-ergodic probability measures is closed in the set of probability measures \cite{EMM}. These results leave open the question of whether the measures obtained by pushing the uniform measure on a bounded piece of a horocycle orbit through a point $q$ by the geodesic flow for a sequence of times $t_n\to \infty$ converge to the unique $\SL$-ergodic measure fully supported on $\overline{\SL  q}$. It is shown in~\cite{EMM} that this convergence holds if one additionally averages over the geodesic flow.
Theorem~\ref{thm:intro} does not directly address this question but it shows that if we are given the freedom to let the initial point $q$ depend on $n$, then ergodicity of the limit cannot be expected.

\subsection{Connections with counting}\label{sec:motivations3} 

Understanding the aforementioned dynamical properties of the action of $P$ and $u_t$ on strata has applications to the geometry and dynamics of translation surfaces and billiards in rational polygons. 
Among the primary examples of these applications is the counting problem for saddle connections~\cite{EskinMasur}.
This requires some explanation. Let $q$ be a translation surface, and let
$$r_\theta=\begin{pmatrix}\cos(\theta)&-\sin(\theta)\\ \sin(\theta)&\cos(\theta)\end{pmatrix}, \quad g_t=\begin{pmatrix}e^t&0\\0&e^{-t}\end{pmatrix}.
$$
Many of the geometric properties of the translation surface can be understood by studying the distribution of the circles $\{g_tr_\theta q:\theta \in [0,2\pi)\}$ (or the uniform measure on this set) for a fixed $q$ as $t$ goes to infinity.
For instance, it is shown in~\cite{EskinMasur} that if these measures converge to an $\SL$-invariant measure, then the number of saddle connections on $q$ of length at most $L$ grows quadratically in $L$, with a rate that can in theory be calculated. Because the pushforward by the geodesic flow of $r_{\theta}$-segments track long (reparametrized) $u_s$-orbits, understanding horocycle limits is a natural approach to this problem.\footnote{Note that for any $\SL$-ergodic measure $\nu$, for $\nu$-almost every $q$, the uniform measure on $\{g_tr_\theta q:\theta \in [0,2\pi)\}$ converges to $\nu$ as $|t|\to \infty$. However, one would like to understand the limit for \emph{every} $q$, since the sets of translation surfaces arising from billiards in rational polygons typically have zero $\nu$-measure.}

This strategy was carried out successfully in certain special settings where it is shown that, in fact, every $U$-orbit is generic for some $U$-ergodic measure~\cite{EskinMarklofMorris,EskinMasurSchmoll,BSW}.
Additionally, these works put strong restrictions on possible limits of $U$-ergodic measures (see for example~\cite[Proposition 11.6]{BSW}) and take advantage of these restrictions to prove their counting results.
In that spirit, it is natural to ask whether there are any restrictions on the closure of the set of horocycle ergodic measures in general. We make this more precise for $\HH_1(2)$ in the following question. 

\begin{ques}\label{q:horocycle limits general}
Is the weak-$*$ limit of horocycle ergodic measures in $\HH_1(2)$ either a horocycle ergodic measure that gives full measure to translation surfaces with a horizontal saddle connection or a convex combination of $\mathrm{SL}_2(\mathbb{R})$-ergodic measures?
\end{ques}
By \cite[Proposition 2.13]{EMM}, a positive answer to this question would show that any limit as $t \to \infty$ of uniform measures on $\{g_tr_\theta q:\theta \in [0,2\pi)\}$ is the unique $\SL$-ergodic measure whose support is $\overline{\mathrm{SL}_2(\mathbb{R})q}$. This would show that every surface in $\HH_1(2)$ satisfies quadratic growth of saddle connections where the constant is the appropriate Siegel-Veech constant; cf.~\cite{EskinMasur}.

The following questions concern the special cases of Question~\ref{q:horocycle limits general}, obtained by pushing a fixed horocycle ergodic measure by $g_t$.

\begin{ques} Given a horocycle ergodic measure $\mu$ on $\HH_1(2)$, does there exist an  $\mathrm{SL}_2(\mathbb{R})$-ergodic measure $\nu$ on the one-point compactification of $\HH_1(2)$ so that $(g_t)_\ast\mu$ converges to $\nu$?
\end{ques}
By results of Forni \cite{ForniDensity1}, for each $\mu$, there is a density $1$ subset $S_\mu\subseteq\mathbb{R}^+$
so that $(g_t)_\ast\mu$ 
converges to $\nu$ along every sequence of $t\in S_\mu$ tending to $\infty$, where $\nu$ is the unique $\mathrm{SL}_2(\mathbb{R})$-ergodic and invariant probability measure supported on the $\mathrm{SL}_2(\mathbb{R})$-orbit closure of $\mrm{supp}(\mu)$.
\begin{ques}Let $\mu$ be a horocycle invariant ergodic measure on a stratum. Assume that $\mu$ is not $\mathrm{SL}_2(\mathbb{R})$-invariant.
Do there exist $\mathrm{SL}_2(\mathbb{R})$-ergodic measures $\nu_1$ and $\nu_2$ on the one-point compactification of the stratum so that $\mrm{supp}(\nu_1)$ is a proper subset of $\mrm{supp}(\nu_2)$ and $(g_{-t})_\ast\mu$ converges to $\nu_1$ while $(g_t)_\ast\mu$ converges to $\nu_2$ as $t\to\infty$?
\end{ques} 
This question can also be asked for orbit closures. The answer is positive for the examples of horocycle-invariant measures and orbit closures constructed in~\cite{BSW}. We conjecture that the answer is also positive for the orbit closures constructed in \cite{CSW}.
It would also be interesting to know if there exists a horocycle ergodic measure $\mu$, that is not $\mathrm{SL}_2(\mathbb{R})$-invariant, but so that $(g_t)_\ast\mu$ and $(g_{-t})_\ast\mu$ converge to the same $\mathrm{SL}_2(\mathbb{R})$-ergodic measure
as $t \to +\infty$?

\subsection{Outline of the article}
The periodic horocycles in Theorem~\ref{thm:intro} are constructed using certain cylinder twists which are special examples of tremors introduced in~\cite{CSW}.
In Section~\ref{sec:background}, we recall several facts regarding tremors and their interaction with the action of the upper triangular subgroup on $\HH_1(2)$. We also give a more precise form of Theorem~\ref{thm:intro} in Theorem~\ref{thm:main} and deduce Corollary~\ref{cor:nongen} from Theorem~\ref{thm:main}. In Sections~\ref{sec:near teich} and \ref{sec:general}, we reduce Theorem~\ref{thm:main} to Theorem~\ref{thm:sufficient}.
A key step in this deduction, carried out in Section~\ref{sec:near teich}, is showing that tremor orbits do not concentrate near proper $\SL$-orbit closures.

Theorem~\ref{thm:sufficient} is the main technical part of our arguments and concerns the non-concentration of the norm of the Kontsevich-Zorich cocycle. Its proof occupies Sections~\ref{sec:oscillation}-\ref{sec:proof of oscillations}. As this outline suggests, the proof of Theorem~\ref{thm:main} splits into two main parts, deducing Theorem~\ref{thm:main} from Theorem~\ref{thm:sufficient} and proving Theorem~\ref{thm:sufficient}.
These arguments are outlined in Sections~\ref{sec:outline of main proof from oscillations} and \ref{sec:outline of oscillations} respectively.

\begin{acknowledgement}
The first author is supported in part by a Simons fellowship, a Warnock chair and NSF grant DMS-145762. The second author is partially supported by NSF grant number DMS-2247713. 
We thank Hamid Al-Saqban and Barak Weiss for helpful conversations that inspired this project.
We thank the anonymous referees for many corrections and comments that helped improve the exposition.
\end{acknowledgement}

    \section{Background}\label{sec:background}

\subsection{Strata, the mapping class group and the GL(2,R)-action}\label{sec:basics}

In this section we recall some basic definitions. The main reference is \cite{CSW} and we make an effort to follow the notation used there. Translation surfaces and their markings are defined in \cite[Section 2.1]{CSW}.
Let $\HHm$ and $\HH_{\mrm{m},1}$ denote the corresponding strata of marked translation surfaces and unit area marked translation surfaces of genus 2 and one cone point of angle $6 \pi$. 
In general a  marked translation surface is given by a map $\phi:(S,\Sigma)\to (M,\Sigma')$ where $S$ is a model surface, $\Sigma$ is a finite subset of $S$ and $\Sigma'$ is the set of cone points in $M$. 
Let $\mathrm{Mod}(S,\Sigma)$ denote the mapping class group of $(S,\Sigma)$, that is the group of isotopy classes of homeomorphisms of $S$ that fix $\Sigma$ . The quotients of $\HHm$ and $\HH_{\mrm{m},1}$ by the right action of $\mathrm{Mod}(S,\Sigma)$ are denoted by $\HHu$ and $\HH_1$ respectively and these are the corresponding strata of (unmarked) translation surfaces (of genus 2 with one cone point of angle $6\pi$) and (unmarked) unit area translation surfaces (of genus 2 with one cone point of angle $6\pi$).  See \cite[Section 2.2]{CSW} for a detailed description of these objects. As in \cite{CSW} we often denote a marked surface by $\tilde{q}$ and the corresponding unmarked surface by $q$.

 We denote by $\GLp$ the subgroup of $\GLd$ consisting of matrices with positive determinant.
 The group $\GLp$ acts on $\HHu$ and $\HHm$ and  the subgroup of elements of determinant one  acts on $\HH_1$ and $\HH_{\mrm{m},1}$.
 We write $A$ and $U$ for the subgroups generated by $g_t$ and $u_s$.
 We use the notation
 \begin{equation*}
     \hat{u}_t=\begin{pmatrix}1&0\\t&1 \end{pmatrix}
 \end{equation*}
  and we denote the corresponding subgroup by $U^-$.
 There is a unique $\mathrm{SL}_2(\mathbb{R})$-ergodic and invariant probability measure on $\HH_1$ of full support, known as the Masur-Veech measure, and we denote it by $\mu_{\mrm{MV}}$. 

\subsection{Invariant splittings of the tangent bundle}
\label{sec:splittings}

Let $\tq\in \HHm$ be a point representing a marked flat surface. We denote the marking map by $\phi:(S,\Sigma)\to (M_q,\Sigma')$.
Then, $\tq$ determines a holonomy homomorphism $\hol_{\tq}$ from $H_1(M_q,\Sigma';\Z)$ to $\R^2$ which we can also interpret as an element of $H^1(M_q,\Sigma';\R^2)$. 
We denote the $x$ and $y$ components of this map by $\hol_{\tilde q}^{(x)}$ and $\hol_{\tilde q}^{(y)}$ respectively.
The cohomology class $\hol_{\tilde q}^{(x)}$ is represented by the 1-form $dx_{\tq}$, viewed as the real part of the holomorphic $1$-form determined by $\tq$. As a map on homology, it is given by $\hol_{\tilde q}^{(x)}[\gamma]=\int_\gamma dx_{\tq}$. See \cite[Section 2.1]{CSW} for more information.

There is a developing map $\mrm{dev}$ from $\HHm$ to $H^1(S,\Sigma;\R^2)$ which sends the marked surface $\tilde q$, with marking given by $\phi:(S,\Sigma)\to (M,\Sigma')$, to $\phi^*(\hol_{\tq})$. The developing map is a local diffeomorphism and identifies the tangent space at $\tilde q$ with the tangent space to $H^1(S,\Sigma;\R^2)$ at $\mrm{dev}(q)$, which we can identify with $H^1(S,\Sigma;\R^2)$. Throughout the article, we identify $\hol_{\tq}$ with its image $\phi^*(\hol_{\tq})$ in $H^1(S,\Sigma;\R^2)$.

A locus $\Lcal$ in $\HHu$ is the orbit closure of a point under the $\GLp$-action. Strata are examples of loci. Eskin, Mirzakhani, and Mohammadi showed in~\cite{EMM} that a proper locus is essentially an affine suborbifold of $\HHm$ in particular it has a well defined dimension.  We will consider a single stratum in this paper which is the stratum of surfaces of genus two with one singular point with cone angle $6\pi$. Since we are only dealing with one stratum we will use the notation $\HHm$ to denote this stratum. The proper loci contained in this stratum are 3-dimensional loci corresponding to Teichm\"uler curves of which the most important for us is the (closed) orbit of the regular octagon surface.

We denote the tangent bundle of $\HHm$ by $T(\HHm)$. Using the developing map, we see that $T(\HHm)$ admits a trivialization as a product bundle (cf.~\cite[Section 2.2]{CSW}):
\begin{equation*}
    T(\HHm)\cong \HHm \times H^1(S,\Sigma;\R^2).
\end{equation*}

If $\Lcal$ is a locus in $\HHu$ then the inverse image of $\Lcal$ in $\HHm$ consists of countably many components. Let $\Lcal_m$ denote one of these components. The restriction of the developing map to $\Lcal_m$ is a local diffeomorphism into a vector space $V\subset H^1(S,\Sigma;\R^2)$. We can identify the tangent bundle of $\Lcal_m$, $T(\Lcal_m)\subset T(\HHm)$ with $\Lcal_m\times V$.

Using the Universal coefficient theorem, we can identify $H^1(S,\Sigma;\R^2)$ with $\Hom(H_1(S,\Sigma;\Z),\R^2)$ or with $H^1(S,\Sigma;\R)\otimes\R^2$. Each of these identifications will play a role in our analysis of the tangent bundle.
The identification of $H^1(S,\Sigma;\R^2)$ with $H^1(S,\Sigma;\R)\otimes\R^2$ is useful in analyzing $\GLp$-invariant subspaces of $H^1(S,\Sigma;\R^2)$ such as the vector space $V$ connected to a locus $\Lcal$.

If we choose a basis for $H_1(S,\Sigma;\Z)$ then we can identify an element of $\Hom(H_1(S,\Sigma;\Z),\R^2)$ with a matrix with 2 rows and $d$ columns. The left action of $\GLp$ is given by left multiplication of matrices and the right action by the monodromy group is given by the right multiplication by $d\times d$ integral matrices. This picture generalizes the classical picture for the torus where the moduli space is identified with a set of $2\times 2$ matrices.

In view of the above description of the $\GLp$-action in coordinates, the derivative $Dg$ of the action of an element $g\in \GLp$ on $\HHm$ takes the following form according to the above identification as
\begin{equation}\label{eq:G derivative}
    Dg = \mrm{Id}\otimes g,
\end{equation}
where $g$ acts on $\R^2$ via the standard left action of $\GLp$ and $\mrm{Id}$ is the identity mapping.
The group $\Mod$ acts on $\HHm \times H^1(S,\Sigma;\R^2)$ by changing the marking on the first factor and its induced action on cohomology on the second.
The quotient space is the (orbifold) tangent bundle of $\HHu$.

We denote the tangent space at $\tilde q \in \HHm$ by $\mrm{T}_{\tilde q}$. 
Let $\hol_{\tilde q}\in \mrm{T}_{\tq}\cong H^1(S,\Sigma;\R^2)$ denote the element corresponding to the holonomy homomorphism from $H_1(S,\Sigma;\Z)$ to $\R^2$.  
Hence, we have $\hol_{\tilde q}^{(x)},\hol_{\tilde q}^{(y)}\in H^1(S,\Sigma;\R)$. We refer to the subspace spanned by $\hol_{\tilde q}^{(x)}$ and $\hol_{\tilde q}^{(y)}$ as \textit{the $\mrm{GL}_2(\R)$-subspace} at $\tq$.

We denote by $\can{\tilde q}\subset\mrm{T}_{\tilde q}$ the tangent space to the $\GLp$-orbit of $\tq$ (where ``st" stands for ``standard"). It is given by 
\begin{equation*}
    \can{\tq} \cong  ( \R\cdot \hol_{\tilde q}^{(x)}\oplus \R\cdot\hol_{\tilde q}^{(y)})\otimes\R^2 ,
\end{equation*}
where, for $v\in H^1(S,\Sigma;\R)$, $\R\cdot v$ denotes the $\R$ span of $v$.
In particular, $\can{\tq}$ is $4$-dimensional.

Given $\tq$, we use the terminology \emph{balanced space} at $\tq$ to refer to the subspace of $H^1(S,\Sigma;\R)$ consisting of cohomology classes whose cup product with all elements of the $\mrm{GL}_2(\R)$-subspace of $\tq$ vanishes.
These two subspaces are complementary and span $H^1(S,\Sigma;\R)$ for all $\tq$.
We define $\bal{\tq}\subset \mrm{T}_{\tq}$ to be the tensor product of the balanced space at $\tq$ with $\R^2$.
By a mild abuse of terminology, we also refer to $\bal{\tq}$ as the balanced space at $\tq$.
We thus obtain a splitting of the tangent space at $\tq$:
\begin{equation}\label{eq:standard and balanced splitting}
    \mrm{T}_{\tq} = \can{\tq} \oplus \bal{\tq},
\end{equation}
This splitting is invariant under the actions of $\GLp$ (via the derivative maps~\eqref{eq:G derivative}) and $\mathrm{Mod}(S,\Sigma)$.  This splitting varies continuously but is not invariant under parallel translation in general. This splitting is discussed in \cite[Section 1.2]{MatYoc} and our notation is motivated by theirs. In one particular case this splitting is invariant under parallel translation: when we restrict to a closed $\GLp$-orbit. In this case the splitting represents the bundle with fiber $H^1(S,\Sigma;\R)$ as a sum of two flat subbundles.
Flatness of this splitting over closed $\GLp$-orbits will be used in Section~\ref{sec:pseudoanosov} to give an explicit description of these bundles in terms of their monodromy over the octagon locus.

Write $\R_x$ and $\R_y$ for the subspaces of $\R^2$ spanned by the standard basis vectors $e_1$ and $e_2$ respectively. By composing with the dual projections $\pi_1$ and $\pi_2$ from $\R^2$ to $\R$, we obtain a splitting of the tangent bundle via the following splitting of the fibers:
\begin{equation}\label{eq:splitting into horizontal and vertical}
    H^1(S,\Sigma;\mathbb{R}^2)=H^1(S,\Sigma;\mathbb{R}_x) \oplus H^1(S,\Sigma;\mathbb{R}_y).
\end{equation}
We refer to the summands in the above splitting as the \textit{horizontal} and \textit{vertical} subspaces respectively. 
Viewing cohomology as maps on the corresponding homology group, the induced action of $\Mod$ on homology induces its action on $H^1(S,\Sigma;\mathbb{R}^2)$ via precomposition. In particular, $\Mod$ leaves the above splitting invariant. We have thus shown:
\begin{lem}\label{lem: split}
The splitting~\eqref{eq:splitting into horizontal and vertical} is $\mathrm{Mod}(S,\Sigma)$-invariant and, hence, induces a splitting of the tangent bundle over $\HHu$.
\end{lem}
With respect to the isomorphism $H^1(S,\Sigma;\R^2)\cong  H^1(S,\Sigma;\R)\otimes\R^2$, we have
\begin{equation}\label{eq:horizontal space as tensor}
    H^1(S,\Sigma;\mathbb{R}_x) \cong   H^1(S,\Sigma;\R)\otimes\R_x .
\end{equation}
The vertical space has the analogous description with $e_1$ replaced with $e_2$.
Given $\tq\in\HHm$, we identify the holonomy components $\hol_{\tq}^{(x)}$ and $\hol_{\tq}^{(y)} \otimes e_1$ with the vectors $ \hol_{\tq}^{(x)}$ and $ \hol_{\tq}^{(y)}\otimes e_2$ respectively.

\subsection{The sup norm metric}
In \cite[Section 2.2.2]{AGY}, Avila-Gou\"ezel-Yoccoz define a family of Finsler norms on $\mrm{T}_{\tilde q}$ and is denoted $\|\cdot \|_{\tilde q}$. We call this family the {\em sup norm}. 
These norms vary continuously in $\tq$~\cite[Proposition 2.11]{AGY}.
They give a complete metric on $\HHm$ which we denote $\distm(\cdot,\cdot)$. The norms and thus the metric are $\mathrm{Mod}(S,\Sigma)$-invariant and so they descend to $\HHu$. We denote the corresponding norm on $\mrm{T}_q$ by $\|\cdot\|_q$ and the metric by $\distu(\cdot,\cdot)$.

\subsection{Tremors}
Tremors were introduced and studied in~\cite{CSW}. In this article, we will use a special case of the general construction in~\cite{CSW} arising from cylinder shears, which we now describe.
We first recall the relevant notation of tremors from \cite[Sections 4.1.2 and 4.1.4]{CSW}. 
Given $\tilde{q}\in \HHm$ the space of tremor vectors or tremor space, denoted $\TT_{\tq} $, is a subspace of the horizontal space $  H^1(S,\Sigma;\mathbb{R}_x)\subset \mrm{T}_{\tilde{q}}$. This space corresponds to tangent vectors constructed from signed measures transverse to the horizontal foliation of $M_{\tilde{q}}$.
Its precise definition is given in \cite[Section 4.1.1]{CSW}. The collection of tremor spaces is $\Mod$-equivariant.
At almost every surface $\tq$ in $\HHm$, $\TT_{\tq}$ consists entirely of scalar multiples of $dy$.

If $q\in\HHu$ is the image of $\tq$, we denote by $\TT_q$ the image of the space $\TT_{\tq}$.
We let
\begin{equation*}
    \TT^{\mrm{bal}}_{\tq}=\TT_{\tq} \cap \bal{\tilde{q}}
\end{equation*}
 and denote by $\TT^{\mrm{bal}}_{q}$ its image in $T(\HHu)$. 
Given $\tau \in \TT_{\tq}$ or $\tau \in \TT_{q}$, the tremor path is a $1$-parameter family\footnote{Technically, this is only true for $s$ in the domain of definition of the tremor, denoted $ \Dom(\tilde{q},\tau)$ or $s \in \Dom(q,\tau)$ in~\cite{CSW}, but the tremors we consider are non-atomic and thus by \cite[Proposition 4.8]{CSW} $\Dom(\tilde{q},\tau)$ and $\Dom(q,\tau)$ are $\mathbb{R}$.} of marked translation surfaces $\mathrm{Trem}_{\tilde q,\tau}(s)$ and (unmarked) translation surfaces $\mathrm{Trem}_{ q,\tau}(s)$ for all $s\in \mathbb{R}.$
By~\cite[Eq.~(4.7)]{CSW}, we have the following description of tremor paths in holonomy coordinates:
\begin{equation}\label{eq:tremor holonomy}
    \hol_{\mathrm{Trem}_{\tilde q,\tau}(s)}=\hol_{\tq}+s\t.
\end{equation}
We can interpret this formula as saying that tremor paths are straight lines (i.e.~affine geodesics) in period coordinates.

A typical tremor of the type that we consider in this paper is obtained from a surface with a horizontal cylinder decomposition. The tremor path is a  family of surfaces obtained by shearing in horizontal cylinders at different rates. A second example of a tremor path is a horocycle orbit. In this case the corresponding transverse measure is given by $dy$. If a surface has a uniquely ergodic horizontal foliation then every tremor vector is a multiple of $dy$ and every tremor path is a reparametrized horocycle orbit. These examples will suffice for the purposes of this paper; cf.~Section~\ref{sec:octagon-tremors}.

The following two basic properties of tremors will be useful.

\begin{lem} 
\label{lem: tremor}
Let $\tq\in\HHm$ and $\tau\in \TT_{\tq}^{\mrm{bal}}$. Then, 
$\tau\in\TT^{\mrm{bal}}_{\mathrm{Trem}_{\tq,\tau}(s)}$ for all $s\in\R$.
\end{lem}
\begin{proof}

Let $q(s)=\mathrm{Trem}_{\tq,\tau}(s)$ and let $(dx,dy)$ be the horizontal and vertical coordinates of $\hol_{q(0)}$ relative to the splitting~\eqref{eq:splitting into horizontal and vertical}. 
That $\t$ remains in the tremor space of $q(s)$ follows by the first bullet point of~\cite[Corollary 6.2]{CSW}.
We show that $\t\in \bal{q(s)}$.
Recalling that $\TT_{\tq}$ is contained in the horizontal subspace, we write $\t=(\t_x,0)$, with $\t_x$ in the horizontal space.
The holonomy of $q(s)$ has coordinates $(dx+s\tau_x,dy)$. The balanced condition at time $s$ is that
$\int \tau \wedge (dx+s\tau_x)=
\int \tau_x\wedge(dx+s\tau_x) =0$ and $\int \tau \wedge dy=\int \tau_x\wedge dy=0$.
Now note that since $\t\in \bal{q}$, by definition we have $\int\t \wedge dx = 0=\int \t\wedge dy$. This immediately implies the second condition.
This also implies $\int \tau_x\wedge(dx+s\tau_x)= \int \tau_x\wedge dx+s\int \tau_x\wedge\tau_x=\int \tau_x\wedge dx=0$ yielding the first condition.

\end{proof}

For the next lemma, recall the notation in equations~\eqref{eq:G derivative} and~\eqref{eq:horizontal space as tensor}. In this notation, we note that the derivative of $u_t$, denoted $Du_t= \mrm{Id} \otimes u_t$, preserves the horizontal space $H^1(S,\Sigma;\R_x)$. Moreover, since $u_t$ fixes $e_1$, $Du_t$ acts trivially on the horizontal space.
In particular, we have that for $\tq\in\HHm$,
\begin{equation}\label{eq:horocycle fixes tremor}
    Du_t(\t) =\t ,\qquad \forall \t\in \TT_{\tq}.
\end{equation}
With this set up, we have the following lemma.
\begin{lem}[Proposition 6.5,~\cite{CSW}] \label{lem:commute}
Tremor maps commute with the horocycle flow in the following sense. For all $s,a\in \mathbb{R}$ and $\tau \in \TT_{q}$, we have that $u_s\mathrm{Trem}_{q,\tau}(a)=\mathrm{Trem}_{u_sq,\tau}(a)$. In particular, the tremor of a periodic $u_s$-orbit is a periodic $u_s$-orbit.
 \end{lem}

\subsection{The octagon locus}\label{sec:octagon}
We can build a translation surface from the regular octagon by identifying pairs of opposite sides. 
We denote by $M_0$ the flat surface obtained in this manner.
We choose a marking $\vp: S\to M_0$ and we denote the corresponding point in $\HHm$ by $\tilde{\w}_0$.
Denote by $\tilde{\Ocal}\subset\HHm$ the $\GLp$-orbit through $\tilde{\w}_0$.
We let $\omega_0\in \HHu$ be the image of $\tilde{\w}_0$. Let $\mathcal{O}\subset \HHu$ denote the orbit
\begin{equation*}
    \Ocal = \GLp\cdot\w_0 \subset \HHu
\end{equation*}
which we call the \textit{octagon locus}.

We recall that $\Ocal$ is a closed subset of $\HHu$ and can be identified with $\GLp/\Gamma$ where $\Gamma$ is the stabilizer of the regular octagon. The group $\Gamma$ meets $\SL$ in a lattice and is in fact a triangle group (see \cite{VeechLat}). The octagon locus is a specific example of a \emph{Teichm\"uller curve}, that is a closed $\GLp$-orbit in a stratum of translation surfaces. Let $\mathcal{O}_1$ denote the closed $\mathrm{SL}_2(\mathbb{R})$-orbit consisting of the area $1$ surfaces in $\mathcal{O}$. On the locus of area $1$ surfaces of a Teichm\"uller curve, there is a unique $\mathrm{SL}_2(\mathbb{R})$-ergodic and invariant probability measure and we let $\mu_{\mathcal{O}}$ denote this measure on $\mathcal{O}_1$. Note that the support of $\mu_{\mathcal{O}}$ is $\mathcal{O}_1$ and not $\mathcal{O}$.

The group $\Gamma$ has an alternate definition in terms of affine automorphisms.
An affine map of a translation surface is a map which is smooth away from the singular points and has constant derivative. The collection of affine automorphisms forms a group and $\Gamma$ is isomorphic to the affine automorphism group of $M_0$.  
  
 Each affine automorphism determines a mapping class in $\Mod$ using the marking $\varphi$.
Hence, we can view the Veech group $\G$ as a subgroup of $\mrm{Mod}(S,\Sigma)$ leaving $\tO$ invariant.
 
\begin{figure}
\includegraphics[scale=0.68]{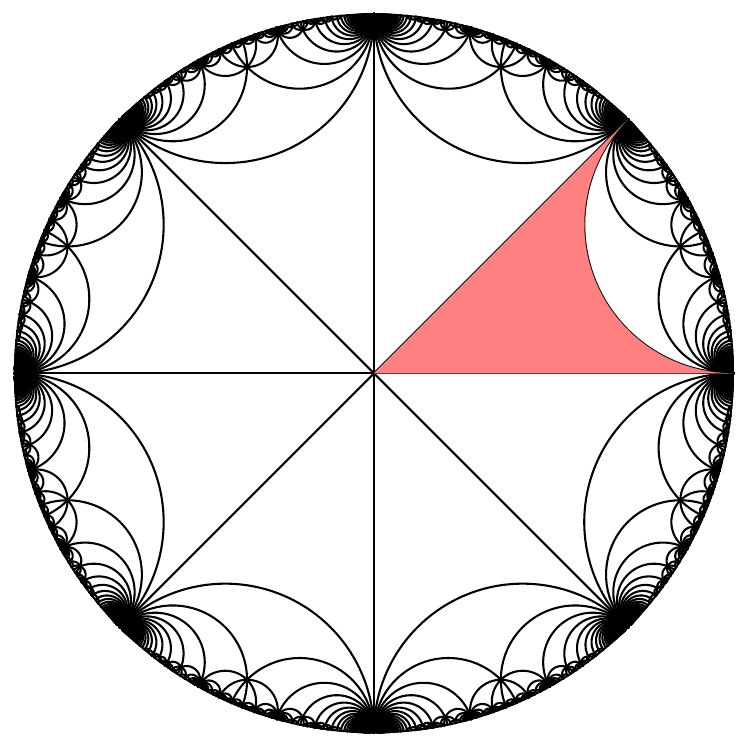}
\caption{A fundamental domain for the Veech group of the regular octagon. The action of the extended Veech group is generated by reflections in the sides of the shaded triangle.}
\label{fig: oct}
\end{figure} 

Note that $\tilde{\Ocal}$ can be identified with the tangent bundle of the hyperbolic plane.
We fix a fundamental domain $\Dcal_0\subset \tO$ for the action of $\G$ on $\tO$.
It will be convenient for arguments in Section~\ref{sec:pseudoanosov} to let $\Dcal_0$ be the fundamental domain constructed in~\cite{SU-Geodesic,SU-Coding}. 
Figure~\ref{fig: oct} shows the projection of $\Dcal_0$ to the hyperbolic plane. Note that the translates of $\Dcal_0$ by all of $\Mod$ are disjoint.

\subsection{The Kontsevich-Zorich (KZ) cocycle over the octagon}\label{sec:cocycle}
 
Recall that the tangent space $\mrm{T}_{\tq}$ at each point $\tq\in \HHm$ is identified with $H^1:= H^1(S,\Sigma;\R^2)\cong  H^1(S,\Sigma;\R)\otimes\R^2 $.
Hence, we identify the tangent bundle over the image of $\Dcal_0$ inside $\Ocal$ with $\Dcal_0\times H^1$.
Moreover, $\Mod$ (and hence $\G$) acts (on the right) by linear automorphisms on $H^1(S,\Sigma;\R)$.
Let $\phi:\mathrm{Mod}(S,\Sigma)\to  \mathrm{Aut}(H^1(S,\Sigma;\R))$ denote this right action.
We note that $\phi$ is a right action since cohomology is contravariant. More explicitly, for $\g_1,\g_2\in \Mod$, and $v\in H^1(S,\Sigma;\R)$,
\begin{equation}\label{eq:contravariant action}
    \phi(\g_1\g_2)(v)= \phi(\g_2)(\phi(\g_1)(v)).
\end{equation}
This action is induced from the action of representatives of the isotopy classes of its elements by homeomorphisms on $S$. This linear action agrees with the derivative of elements of $\Mod$ acting by diffeomorphisms of $\HHm$.

The above trivialization allows us to describe the derivative of the $\GLp$-action on $\Ocal$ on the tangent bundle to $\HHu$ as follows.
Let $\tq\in \Dcal_0\subset\HHm$ and $g\in \GLp$. Let $\g\in \G$ be the unique element satisfying $g\tq\cdot\g \in \Dcal_0$.
The KZ cocycle is defined as follows
\[ \mrm{KZ}(g,\tq):= \phi(\g).\]
In view of~\eqref{eq:contravariant action}, we obtain the following cocycle relation
\begin{equation*}
    \mrm{KZ}(gh,\tq) = \mrm{KZ}(g,h\tq)\circ\mrm{KZ}(h,\tq), \qquad \forall g,h\in \GLp.
\end{equation*}
In the remainder of the article, we drop the composition notation and simply write $\mrm{KZ}(gh,\tq) = \mrm{KZ}(g,h\tq)\mrm{KZ}(h,\tq)$.

A key component in our arguments is the action of the geodesic flow $g_t$ on tremors.
To describe this action precisely, let
$\tq_1\in \Dcal_0$ and suppose that $\tq_2\in\Dcal_0$ is a representative of $g_t\tq_1$.
Let $q_1\in \Ocal$ be the image of $\tq_1$ in $\HHu$ and set $q_2=g_tq_1$.
The restriction of the derivative $Dg_t:\mrm{T}_{q_1}\to \mrm{T}_{q_2}$ to the horizontal space admits the following description.
Given a horizontal vector $\t \in H^1(S,\Sigma;\R_x)$ at $\tq_1$, we write $\t= \t'\otimes e_1\in  H^1(S,\Sigma;\R)\otimes\R^2$ via the identification~\eqref{eq:horizontal space as tensor}.
Then, $Dg_t(\t) =  \KZ{g_t}{\tq_1}\t'\otimes e^t\cdot e_1$, where we identify $\t$ with its image in $\mrm{T}_{q_1}$ and both sides of the equation with their respective images in $\mrm{T}_{g_t q_1}$.
For simplicity, when dealing with the action of the upper triangular group $P$, and in particular, with the action of $g_t$, we will use the notation
\begin{equation}\label{eq:cocycle on horizontal}
    Dg_t(\t) = e^t \cdot \KZ{g_t}{\tq_1}\t \in H^1(S,\Sigma;\mathbb{R}_x)
\end{equation}
throughout the remainder of the article.
This is justified by the invariance of the horizontal space under $P$.

Eq.~\eqref{eq:cocycle on horizontal} describes the action of $g_t$ on $H^1(S,\Sigma;\mathbb{R}_x)\subset H^1(S,\Sigma;\mathbb{R}^2)$ via the KZ cocycle.
It is shown in~\cite{CSW} that tremor subspaces are equivariant under this action. 
The next lemma summarizes these results and is our key tool to apply renormalization dynamics in our setting.

\begin{lem} \label{lem:renorm}
Let $\tilde{q} \in \Dcal_0$ and $\tau  \in \TT_{\tq}$. 
Then, $Dg_t(\t)\in\TT_{g_t\tq}$ and $g_t\mathrm{Trem}_{\tilde{q},\tau}(s)=\mathrm{Trem}_{g_t\tilde{q},\tau}(e^ts)$.
If $q\in \Ocal$ is the image of $\tq$, then
\begin{equation*}
   g_t\mathrm{Trem}_{q,\tau}(s)=
    \mathrm{Trem}_{g_tq,\tau_2}(e^ts),
\end{equation*}
 where $\tau_2=\mrm{KZ}(g_t,\tq)\tau$.
 In particular, $\tau_2$ is an element of $\TT_{g_t q}$.

\end{lem}

\begin{proof}
That $Dg_t(\t)\in \TT_{g_t\tq}$ follows by the proof of~\cite[Proposition 6.5]{CSW}.
The second claim follows by~\eqref{eq:G derivative} and~\eqref{eq:horizontal space as tensor} along with the fact that tremors belong to the horizontal subspace.
In particular, this is the marked stratum version of the first equation in~\cite[Eq.~(6.4)]{CSW}. 
The last claim follows from this fact and~\eqref{eq:cocycle on horizontal} since tremors are equivariant under $\Mod$.

\end{proof}

\subsection{Tremors of the octagon locus} 
\label{sec:octagon-tremors}

The unit area octagon locus $\Ocal_1\cong \mathrm{SL}_2(\mathbb{R})/\Gamma$ has two cusps corresponding to the two ideal vertices of the fundamental domain. Each cusp corresponds to a family of closed horocycles. In each family the closed horocycles are parametrized by the periods of the horocycles.
Let $\w_1\in\Ocal_1$ denote a surface that is on a horocycle of period exactly 1; i.e.
\begin{equation}\label{eq:omega1}
    u_1 \w_1= \w_1
\end{equation}
and denote by $\mc{C}$ the $u_s$-orbit of $\w_1$.
Let $\tilde{\omega}_1\in \Dcal_0$ and $\tilde{\mathcal{C}}=\set{u_s\tilde{\w}_1:0\leq s\leq 1}$
denote lifts of $\w_1$ and $\mc{C}$ respectively.

The surface $\tilde \omega_1$ has two cylinders, with different circumferences. Let $\cyl_a$ be the cylinder with shorter circumference and denote its area by $a$ and let $\cyl_b$ be the cylinder with longer circumference and denote its area by $b$.
Let $\chi_{\cyl_a}$ and $\chi_{\cyl_b}$ denote the indicator functions of $\cyl_a$ and $\cyl_b$ respectively.

Let $\sigma\in H^1(S,\Sigma;\R^2)$ be the cohomology class
\begin{equation}\label{eq:sigma}
   \s= [b\chi_{\cyl_a} dy-a\chi_{\cyl_b}dy],
\end{equation}
where $b\chi_{\cyl_a} dx-a\chi_{\cyl_b}dx$ is a (signed) $1$-form and $dx$ and $dy$ are the $1$-forms representing the cohomology classes $\hol_{\tilde{\w}_1}^{(x)}$ and $\hol_{\tilde{\w}_1}^{(y)}$ respectively; cf.~Section~\ref{sec:splittings}. Notice that $\sigma\in \TT^{\mrm{bal}}_{\tilde{\omega}_1}$. 
 Indeed, it is clear that $$(b\chi_{\cyl_a} dy-a\chi_{\cyl_b}dy)\wedge dy=0.$$
 Moreover, by our choice of $a,b$ we have that
 $$\int_S(b\chi_{\cyl_a} dy-a\chi_{\cyl_b}dy)\wedge dx=ba-ab=0.$$
 We extend this definition to $\omega_1$ and the surfaces in $\tilde{\mathcal{C}}$ and $\mathcal{C}$. Note we can extend this to $\tilde{\mathcal{C}}$ and $\mathcal{C}$ because $u_s$ preserves the  direction, area and circumference of horizontal cylinders.

\subsection{Restatement of Theorem \ref{thm:intro}}

We establish Theorem \ref{thm:intro} by proving a more specific theorem, which requires some notation. Denote by $\nu_t$ the $u_s$-ergodic and invariant probability measure supported on the periodic horocycle $g_t\mathcal{C}$ and set $\nu=\nu_0$. 
For $a, t\in \R$, we let $\nu_{t,a}$ denote the $u_s$-ergodic probability measure on the periodic horocycle through $g_t\mathrm{Trem}_{\w_1,\sigma}(a) = \mathrm{Trem}_{g_t\w_1,\sigma_t}(e^ta)$, where $\s$ is given by~\eqref{eq:sigma} and
\begin{align*}
    \sigma_t = \KZ{g_t}{\w_1}\s \in \TT_{g_t\w_1}.
\end{align*}
Here, we applied Lemma~Lemma \ref{lem:renorm} applied with $\tau=\sigma$.
We also used Lemma \ref{lem:commute} to ensure that $\mathrm{Trem}_{g_t\w_1,\sigma_t}(e^ta)$ is $u_s$-periodic.
In particular, $(\nu_{t,a})$ is a $2$-parameter family of $u_s$-ergodic probability measure and we show that Theorem~\ref{thm:intro} holds for a sequence of elements of this family. More precisely, we prove the following.

\begin{thm}\label{thm:main}
There exist two sequences of real numbers $t_i\to \infty$ and $(a_i)_i$ so that the weak-$\ast$ limit of $\nu_{t_i,a_i}$ is a non-trivial convex combination of $\mu_{\mathcal{O}}$ and $\mu_{\mrm{MV}}$.
\end{thm} 

Corollary~\ref{cor:nongen} can be deduced from this theorem as follows.

\begin{proof}[Proof of Corollary \ref{cor:nongen}]

We prove that there are two dense $\mrm{G}_\d$-subsets, $V_1$ and $V_2$, where the $u_s$-orbits of points in $V_1$ equidistribute along some subsequence to a measure other than $\mu_{\mrm{MV}}$ and the $u_s$-orbits of points in $V_2$ equidistribute along some subsequence to  $\mu_{\mrm{MV}}$. By the Baire category theorem the intersection, $V_1\cap V_2$, is a dense $\mrm{G}_\d$ subset where the $u_s$-orbits of points in it equidistribute along (different subsequences) to different measures, giving Corollary \ref{cor:nongen}. 

The construction of $V_1$ is more involved and we do it first. Let $V_1'$ be the union of the periodic horocycles supporting the measures $\nu_{t_i,a_i}$ provided by Theorem~\ref{thm:main}.
Since the measures $\nu_{t_i,a_i}$ converge to a measure of full support, the set $V_1'$ is dense.
Denote by $d_\ast$ a complete metric on the space of probability measures that induces the same topology as the weak-$\ast$ topology\footnote{Concretely, given a countable dense set $\set{\vp_n}$ of elements of $C_c(\HH_1)$, we may define $d_\ast(\mu,\nu)$ to be $ \sum_n 2^{-n}\left|\int \vp_n \;d\mu - \int \vp_n\;d\nu\right|$ for any probability measures $\mu$ and $\nu$.}.
Since the limit of the measures $(\nu_{t_i,a_i})_i$ is different from $\mu_{\mrm{MV}}$, we can find a constant $c>0$ so that $d_\ast(\nu_{t_i,a_i},\mu_{\mrm{MV}})> c$ for all $i$.
Note that for every $i$, the periodic horocycle supporting $\nu_{t_i,a_i}$ has period $e^{2t_i}$ by~\eqref{eq:omega1}.
Hence, one checks that the set
\begin{align*}
    V_1 
    := \bigcap_{m\in\N} \bigcup_{T\geq m, T\in\R}
    \left\lbrace 
    x\in \HHu:
    d_\ast\left(\frac 1 {T}\int_{0}^{T}\d_{u_sx}ds,\mu_{\mrm{MV}} \right)>c
    \right\rbrace
\end{align*}
is a $\mrm{G}_\d$-set containing $V_1'$.
In particular, $V_1$ is a dense $\mrm{G}_\d$-set.

Define $V_2$ to be the set of points whose orbits equidistribute along a subsequence to $\mu_{\mrm{MV}}$. This is a $\mrm{G}_\d$ set. Because $\mu_{\mrm{MV}}$ has full support and is ergodic, $V_2$ contains a dense set of points whose $u_s$-orbits in fact equidistribute to $\mu_{\mrm{MV}}$. In particular, $V_2$ is a dense $\mrm{G}_\d$ set as desired, thus completing the proof of the corollary.
\end{proof}


\section{Avoidance of Teichm\"uller curves}
\label{sec:near teich}

    The goal of this section is to show that, asymptotically, balanced tremor orbits do not give mass to any Teichm\"uller curve, Proposition~\ref{prop:near teich}.
    This allows us to show that the Masur-Veech measure occurs as one of the ergodic components of the measure in the conclusion of Theorem~\ref{thm:main}. The key mechanism of the proof is the transversality between the tangent space $\can{}$ to $\GLp$-orbits and the balanced space $\bal{}$; cf.~discussion at the beginning of Section~\ref{sec:avoidance proof} for a more detailed outline of the argument.

\begin{prop}\label{prop:near teich}
For every $\e>0$, Teichm\"uller curve $\VV$ and open precompact set 
$\Kcal\subset \HHu$ there exists $\d>0$ so that for all $T\geq 1$, $q\in \HHu$ and $\beta \in \TT_{q}^{\mrm{bal}}$ with $\norm{\b}_q=1$, we have 
\begin{equation*}
\left| \{ |\ell| \leq T:\distu(\mathrm{Trem}_{q,\b}(\ell),\VV)<\d  
 \text{ and }\mathrm{Trem}_{q,\b}(\ell) \in \Kcal \} \right|<2\e T.
\end{equation*}
\end{prop} 

    We deduce Proposition~\ref{prop:near teich} from the following local estimate.
    
    \begin{prop}\label{prop:near teich conn comp}
     Given a Teichm\"uller curve $\VV$ and open precompact set $\Kcal\subset \HHu$, there exist $\d_0>0$ and $C\geq 1$ such that the following holds.
     Let $T\geq 1$, $q\in\HHu$ and $0\neq \b\in\TT_q^{\mrm{bal}}$ be arbitrary.
     Given $\d>0$, define
     \begin{equation}\label{eq:local exceptional sets}
         E(\d,T) = \set{|\ell|\leq T: \distu(\mathrm{Trem}_{q,\b}(\ell),\VV) < \d \text{ and }\mathrm{Trem}_{q,\b}(\ell) \in \Kcal}.
     \end{equation}
     Let $I$ be any connected component of $E(\d_0,T)$.
     Then, for every $0<\d<\d_0$, we have that
         \begin{equation*}
        |I\cap E(\d, T)| \leq  \frac{C\d }{\sup_{\ell\in I} \distu(\mathrm{Trem}_{q,\b}(\ell),\VV)} |I|.
    \end{equation*}
    Moreover, if $I=[-T,T]$,  then
    \begin{equation*}
        |I\cap E(\d, T)|=|E(\d,T)| \leq  \frac{C\d }{\norm{\b}_{q}}\cdot 2T.
    \end{equation*}

    \end{prop}

\begin{proof}[Proof of Proposition \ref{prop:near teich} from Proposition~\ref{prop:near teich conn comp}] 
    
    Let $\d_0>0$ and $C\geq 1$ be the constants provided by Proposition~\ref{prop:near teich conn comp}.
    For $0<\d<\d_0$ and $T\geq 1$, denote by $\mc{I}=\set{I_n}$ the collection of connected components of $E(\d_0,T)$. Let $q_\ell=\mathrm{Trem}_{q,\b}(\ell)$. 
    
    If $\mc{I}=\emptyset$, there is nothing to prove.
    If $\mc{I}$ consists of a single element $I$ such that $I=[-T,T]$, then the statement follows by the second assertion in Proposition~\ref{prop:near teich conn comp} by taking $\d<\e/C$ since $\norm{\b}_q=1$.

    Finally, suppose $I\in \mc{I}$ satisfies $I\neq [-T,T]$.
    Then, there is a boundary point $\ell$ of $I$ such that $\distu(q_\ell,\VV)=\d_0$.
    In this case, we see that $|I\cap E(\d,T)|\leq C\d |I|/\d_0$. 
    Hence, since $E(\d,T)\subseteq E(\d_0,T)$, we obtain
    \begin{equation*}
        |E(\d,T)| = \sum_{I\in\mc{I}} |E(\d,T)\cap I|
        \leq  C\d/\d_0 \sum_{I\in\mc{I}} |I|
        \leq  C\d/\d_0 \cdot 2T.
    \end{equation*}
    This concludes the proof by taking $\d=\e \d_0/C$. \qedhere
\end{proof}

Before the proof of Proposition~\ref{prop:near teich conn comp}, we need the following lemma.
It relates the sup-norm distance $\distu$ between nearby points in $\HHu$ to the norm of a suitable vector in the tangent space. The reader is referred to~\cite[Propositions 5.3 and  5.5]{AG} for related results.

\begin{lem}
\label{lem:neighborhood control}
For every $\tq\in \HHm$, there exist $r_{\tq}>0$ so that for all $\tq_1,\tq_2\in B(\tq,r_{\tq})$, we have
$$\frac 1 4 \distm(\tq_1,\tq_2)\leq \|\hol_{\tq_1}-\hol_{\tq_2}\|_{\tq}\leq 4 \distm(\tq_1,\tq_2). $$ 
Moreover, given a compact set $\tK\subset \HHm$, there exists $r_{\tK}>0$ such that $r_{\tq}\geq r_{\tK}$ for all $\tq\in \tK$.
\end{lem}
\begin{proof} 
We choose a small neighborhood $U$ of $\tq\in \HHm$ where $\|\cdot \|_{\tq}$ changes by at most a factor of 2 and so that $\hol:U \to H^1(S,\Sigma;\mathbb{R}^2)$ is injective. In particular, for every $\tq'\in U$, $\beta \in H^1(S,\Sigma;\mathbb{R}^2)$ we have 
$\frac 1 2 \|\beta\|_{\tq}\leq\|\beta\|_{\tq'}\leq 2\|\beta\|_{\tq}.$ Let $V$ be  a small neighborhood of $\tq$ such that any sequence of  $\distm$ minimizing paths 
between two points in $V$ has that its elements eventually stay in $U$.
Now, in the normed vector space $(H^1(S,\Sigma;\mathbb{R}^2),\|\cdot\|_{\tq})$, the straight line segment between $\hol_{\tq_1}$ and $\hol_{\tq_2}$ (which is a translate of the line from $0$ to $\hol_{\tq_2}-\hol_{\tq_1}$) is a geodesic.
If $\tq_1, \, \tq_2\in V$, the corresponding path in $\HHm$ has at most twice this length. That is, if $\gamma$ is the line segment in $H^1(S,\Sigma;\mathbb{R}^2)$ joining $\hol_{\tq_1}$ and $\hol_{\tq_2}$, then 
$$ 
\text{Length}(\hol^{-1}\gamma)\leq 2\|\hol_{\tq_2}-\hol_{\tq_1}\|_{\tq}, $$ 
    where for a $C^1$ map $\k: [a,b]\to \HHm$, $\mrm{Length}(\k)=\int_a^b \norm{\k'(t)}_{\k(t)} \;dt$.

 Similarly, any $C^1$ path $\g$ in $U$ has the property that the length of the curve $\hol(\g)$ in $H^1(S,\Sigma;\mathbb{R}^2)$ with respect to the metric coming from $\|\cdot\|_{\tq}$ is at most $2\mrm{Length}(\g)$. This gives the lower bound.
 
 The uniformity over compact sets follows by choosing a finite cover of $\tK$ by open sets $U$ as above. The lemma follows.
\end{proof}

\subsection{Proof of Proposition \ref{prop:near teich conn comp}}
\label{sec:avoidance proof}

    The argument has two main steps.
    First, we take advantage of the local nature of the problem to reduce the statement to one regarding estimating the proportion of time a connected tremor path spends near a lift of the Teichm\"uller curve $\VV$ in the marked stratum $\HHm$ (cf.~Claim~\ref{claim:lifts of T are locally finite} below).
    The second step is to linearize the latter estimate using the fact that tremor paths are straight lines in period coordinates (cf.~\eqref{eq:tremor holonomy}).
    With this setup in hand, the desired estimate will follow from the observation that balanced tremor straight lines are transverse to the image of $\VV$ in period coordinates (cf.~Claim~\ref{claim:coeff of tremor poly } below). 
    
    To simplify notation, let
    \begin{align*}
        H=H^1(S,\Sigma;\mathbb{R}^2).
    \end{align*}
    The marking maps provide an identification of fibers of the tangent bundle of $\HHm$ with $H$ so that we may regard all the spaces $\can{\bullet}$ and $\bal{\bullet}$ as subspaces of $H$.
    Moreover, for all $g\in \GLp$ and $\tq\in \HHm$, $\can{g\tq}=\can{\tq}$ and, hence, $\bal{g\tq}=\bal{\tq}$.
    In particular, the spaces $\can{\bullet}$ and $\bal{\bullet}$ are constant over $\GLp$-orbits.
    If $\mathfrak{T}\subset\HHm$
    is one such $\GLp$-orbit, we denote by $V^{\mrm{st}}_{\mathfrak{T}}$ and $V^{\mrm{bal}}_{\mathfrak{T}}$ the common value over $\mathfrak{T}$ of $\can{}$ and $\bal{}$ respectively.
    Finally, by Lemma~\ref{lem: tremor}, for $\tq\in\HHm$ and $\beta\in \TT^{\mrm{bal}}_{\tq}$, we have
    \begin{equation*}
       \b\in \TT^{\mrm{bal}}_{\mathrm{Trem}_{\tq,\beta}(\ell)}, \qquad \forall \ell\in\R.
    \end{equation*}
    
    Fix an open precompact set $\Kcal$ and a Teichm\"uller curve $\VV$ in $\HHu$.
    Let $\tilde{\Kcal}\subset \HHm$ denote a compact set inside the closure of some fixed fundamental domain for $\mrm{Mod}(S,\Sigma)$ and which projects to the closure of $\Kcal$.
    Let $\tK_1$ denote the closure of the $1$-neighborhood of $\tK$ in $\HHm$.

    Let $\d_0>0$ be a small parameter whose value is to be determined. Over the course of the proof, we will assume $\d_0$ to be small enough depending only on $\Kcal$ and $\VV$.
    
    Fix $q\in\HHu$, $0\neq \beta\in \TT^{\mrm{bal}}_q$, and $T\geq 1$. Recall the sets $E(\cdot,T)$ defined in~\eqref{eq:local exceptional sets} and fix a connected component $I$ of $E(\d_0,T)$.
    If $\ell_0\in [-T,T]$ is the center of the interval $I$, by replacing $q$ with $\mrm{Trem}_{q,\b}(\ell_0)$, we may assume in the sequel that $I$ is centered at $0$ and that $q\in \Kcal$.

    Denote by $\tq\in\tK\subset\HHm$ a lift of $q$. 
    For $\ell\in\R$, let $\tq_\ell=\mathrm{Trem}_{\tq,\b}(\ell)$ and $q_\ell=\mathrm{Trem}_{q,\b}(\ell)$.
    Given a $\GLp$-orbit $\tVV \subset \HHm$ that projects to $\VV$ and $0<\d<1$, let
    \begin{align}\label{eq:lift of E(delta,T)}
        E(\d,I,\tVV) = 
        \set{\ell\in I: \distm(\tq_\ell,\tVV) < \d }.
    \end{align}
    It follows that
    \begin{align}\label{eq:union over lifts}
        I\cap E(\d,T) \subseteq \bigcup_{\tVV} E(\d,I,\tVV),
    \end{align}
    where the union runs over lifts $\tVV$ of $\VV$.

    \begin{claim}\label{claim:lifts of T are locally finite}
     
    If $\d_0$ is small enough, depending on $\Kcal$ and $\VV$, then there exists a unique lift\footnote{This unique lift $\tVV$ depends on our earlier choices of lifts $\tq$ and $\tK$.} $\tVV$ of $\VV$ such that $I=E(\d_0,I,\tVV)$.
    In particular, for all $\d\leq\d_0$, $I\cap E(\d,T) =E(\d,I,\tVV) $.

    \end{claim}

    \begin{proof}
        
        Note that the second assertion follows from the first in light of~\eqref{eq:union over lifts} and the fact that $E(\d,I,\tVV)\subseteq E(\d_0,I,\tVV)$ whenever $\d\leq \d_0$.
        
        To prove the first assertion, observe that the countable collection $\set{\tVV}$ of lifts of $\VV$ is locally finite, i.e. only finitely many of these lifts meet any given compact subset of $\HHm$.
        This follows by proper discontinuity of the action of $\mrm{Mod}(S,\Sigma)$ on $\HHm$ and the fact that $\VV$ is closed in $\HHu$.
        
        Recall that $\tK_1$ denotes the closure of the $1$-neighborhood of $\tK$ in $\HHm$ and note that if $E(\d_0,I,\tVV)$ is non-empty for some $\tVV$, then $\tVV$ meets $\tK_1$.
        Let $\tVV_1,\dots, \tVV_n$ denote the lifts of $\VV$ which meet $\tK_1$.
        Since the sets $\tVV_i\cap\tK_1$ are closed and disjoint, they are uniformly separated by some $0<\d_1<1$ in the metric $\distm$.

        Let $\d_0=\d_1/2$.
        Recall that $I$ is a connected component of $E(\d_0,T)$.
        Hence, applying~\eqref{eq:union over lifts} with $\d=\d_0$, we see that the sets $E(\d_0,I,\tVV)$ form a cover of $I$ by disjoint open subsets of $I$.
        Thus, at most one of them can be non-empty by connectedness of $I$ which concludes the proof of the claim. 
        \qedhere
        
    \end{proof}

    For the remainder of the proof, we assume that $\d_0$ is small enough so that Claim~\ref{claim:lifts of T are locally finite} holds and denote by $\tVV$ the unique lift provided by the Claim.
    In particular, to prove the first assertion of the proposition, it suffices to show that for all $0<\d <\d_0$
    \begin{align*}
        |E(\d, I,\tVV)| \leq \frac{C \d |I|}{\sup_{\ell\in I} \distu(q_\ell,\VV) },
    \end{align*}
    for a constant $C\geq 1$ depending only on $\tK$.
    Having lifted our problem from $\HHu$ to $\HHm$, our next step is to transfer the estimates into the linear space $H$ using period coordinates.
    This requires some preparation.

     It will be convenient to fix some Euclidean inner product $\langle\cdot,\cdot\rangle$ on $H$ in which the spaces $V^{\mrm{st}}_{\tVV}$ and $V^{\mrm{bal}}_{\tVV}$ are orthogonal. Denote the induced norm by $\norm{\cdot}_0$.
    Let $d$ denote the common dimension of $\bal{\tw}$ for $\tw\in \HHm$ and let $\mrm{Gr}(d,H)$ be the Grassmannian of $d$-dimensional planes in $H$.
    Our Euclidean structure induces a metric on $\mrm{Gr}(d,H)$, which we denote $\dist_0$, given by
    \begin{align*}
        \dist_0(W_1,W_2) := \sup_{v\in W_1} \inf_{w\in W_2} \measuredangle(v,w), \qquad
        \forall W_1, W_2\in \mrm{Gr}(d,H),
    \end{align*}
    where $\measuredangle (v,w) := \frac{\langle v,w\rangle}{\norm{v}_0 \norm{w}_0}$ denotes the Euclidean angle between $v$ and $w$.
    Continuity and compactness then imply that the map $\tw\mapsto \can{\tw}$ is uniformly continuous as a map from $\tK_1$ to the  $\mrm{Gr}(d,H)$.
    In particular, for every $\eta_1>0$, we can find $\eta_2 = \eta_2(\eta_1,\tK,\tVV)>0$ so that  for all $\tw_1,\tw_2\in \tK_1$, 
      \begin{equation}\label{eq:tautological continuity}
       \distm(\tw_1,\tw_2) < \eta_2 
        \Longrightarrow
        \dist_0( \bal{\tw_1}, \bal{\tw_2}) <\eta_1.
  \end{equation}

    Since the AGY norms vary continuously over $\HHm$, we can also find a constant $C_{\tK}\geq1$ such that
    \begin{equation}\label{eq:AGY equiv to ad hoc}
        C_{\tK}^{-1} \norm{\cdot}_{\tw} \leq
        \norm{\cdot}_0 \leq C_{\tK} \norm{\cdot}_{\tw},
        \qquad \forall \tw \in \tK.
    \end{equation}

    Denote by $\pi_{\tVV}^{\mrm{st}}$ and $\pi_{\tVV}^{\mrm{bal}}$ the orthogonal projections (relative to our fixed inner product) from $H$ onto $V_{\tVV}^{\mrm{st}}$ and $V_{\tVV}^{\mrm{bal}}$ respectively.
    We observe that for every $\tilde{\w}\in\HHm$ and for every $\beta\in \TT_{\tilde{\w}}^{\mrm{bal}}$, the map
    \begin{equation}\label{eq:polynomial distance}
       \ell\mapsto P(\tilde{\w},\tVV;\ell):=  \norm{\hol_{\mathrm{Trem}_{\tilde{\w},\beta}(\ell)}-\pi^{\mrm{st}}_{\tVV} (\hol_{\mathrm{Trem}_{\tilde{\w},\beta}(\ell)}) }_0^2 
    \end{equation}
    is a polynomial in $\ell$ of degree $\leq2$.
    This follows from the linear relation $\hol_{\mathrm{Trem}_{\tilde{\w},\beta}(\ell)} = \hol_{\tilde{\w}} + \ell \beta $; cf.~\eqref{eq:tremor holonomy}.
    This polynomial is nothing but the squared distance of $\hol_{\mathrm{Trem}_{\tilde{\w},\beta}(\ell)}$ to $V^{\mrm{st}}_{\tVV}$ inside $H$.
    Note that we suppress the dependence on $\b$ in the notation $P(\tilde{\w},\tVV;\ell)$.

    \begin{claim}\label{claim:coeff of tremor poly }
    For every $\tilde{\w}\in\HHm$ and $\b\in\TT_{\tilde{\w}}^{\mrm{bal}}$,
     the coefficient of the quadratic term in $P(\tilde{\w},\tVV;\cdot)$ is $\norm{\pi^{\mrm{bal}}_{\tVV}(\beta)}_0^2$.
     Moreover, there exists $\eta = \eta(\tK,\tVV)$ such that
     \begin{equation}\label{eq:nontrivial proj}
        \norm{\pi^{\mrm{bal}}_{\tVV}(\b)}_0 \geq  \norm{\b}_0/2,
    \end{equation}
     whenever $\tilde{\w}\in \tK$ and $\distm(\tilde{\w},\tVV)<\eta$. 
    \end{claim}
    
    \begin{proof}
    Note that $V^{\mrm{st}}_{\tVV}$ and $V^{\mrm{bal}}_{\tVV}$ are orthogonal and span $H$ (cf.~\eqref{eq:standard and balanced splitting}).
    Hence, the equation $\hol_{\mathrm{Trem}_{\tilde{\w},\beta}(\ell)} = \hol_{\tilde{\w}} + \ell \beta $ implies that
    \begin{equation*}
        \hol_{\mathrm{Trem}_{\tilde{\w},\beta}(\ell)}-\pi^{\mrm{st}}_{\tVV} (\hol_{\mathrm{Trem}_{\tilde{\w},\beta}(\ell)})
        = \pi^{\mrm{bal}}_{\tVV}(\hol_{\tilde{\w}}) 
        +\ell\cdot \pi^{\mrm{bal}}_{\tVV}(\beta).
    \end{equation*}
    It follows that the coefficient of $\ell^2$ is $\norm{\pi^{\mrm{bal}}_{\tVV}(\beta)}_0^2$.
   The second assertion follows by
    ~\eqref{eq:tautological continuity}.
    \end{proof}

    Let $r_{\tK}>0$ and $\eta>0$ be the constants provided by Lemma~\ref{lem:neighborhood control} and Claim~\ref{claim:coeff of tremor poly } respectively.
    By making the constant $\d_0$ chosen in Claim~\ref{claim:lifts of T are locally finite} smaller, we may assume that
    \begin{equation}\label{eq:delta0}
        \d_0 < \min\{\eta,r_{\tK} \}.
    \end{equation}
    Let $0<\d<\d_0$ be given and let
    $\ell\in E(\d,I,\tVV)$.
    Recall the notation $q_\ell$ and $\tq_\ell$ introduced above~\eqref{eq:lift of E(delta,T)}.
    Let $p_\ell \in \VV$ be such that $\distu(q_\ell,\VV)=\distu(q_\ell,p_\ell)$.
    Let $\tp_\ell\in\tVV\cap \tK_1$ be a lift of $p_\ell$ such that
    \begin{align}\label{eq:dist to VV and tVV are equal}
        \distu(q_\ell,\VV) = \distm(\tq_\ell,\tp_\ell) =
        \distm(\tq_\ell,\tVV).
    \end{align}
    By Lemma~\ref{lem:neighborhood control} and~\eqref{eq:AGY equiv to ad hoc}, we have
    \begin{align*}
        \distu(q_\ell,\VV) =
        \distm(\tq_\ell,\tp_\ell) 
        \geq\frac{1}{4} 
        \norm{\hol_{\tq_\ell} -\hol_{\tp_\ell}}_{\tq_\ell} 
        \geq \frac{1}{4 C_{\tK}}
        \norm{\hol_{\tq_\ell} -\hol_{\tp_\ell} }_0.
    \end{align*}
    Note that $\hol_{\tp_\ell}\in V^{\mrm{st}}_{\tVV}$.
    Denoting by $\pi^{\mrm{st}}_{\tVV} $ the orthogonal projection onto $V^{\mrm{st}}_{\tVV}$, it follows that
    \begin{align}\label{eq:linearizing distance}
       \distm(\tq_\ell,\tVV)
        \geq \frac{1}{4 C_{\tK}}
        \norm{\hol_{\tq_\ell}  - \pi^{\mrm{st}}_{\tVV}(\hol_{\tq_\ell}) }_0 
        = \frac{1}{4 C_{\tK}}\sqrt{P(\tq,\tVV;\ell)}.
    \end{align}
    Hence, we obtain 
   \begin{equation*}
        E(\d,I,\tVV) \subseteq 
        \left\{ \ell\in I: |P(\tq,\tVV;\ell)| < 16 C_{\tK}^2 \d^2
        \right\}.
    \end{equation*}
     
    In particular, we can apply what is commonly called
the \emph{$(C,\alpha)$-good} property of polynomials, cf.~\cite[Proposition 3.2]{Kleinbock-Clay} and~\cite[Lemma 4.1]{DaniMargulis}, to get
    \begin{equation}\label{eq:linearization}
       | E(\d,I,\tVV)| \leq |I| \cdot \frac{16C_{\tK} \d}{\left(\sup_{\ell\in I} |P(\tq,\tVV;\ell)|\right)^{1/2}}.  
    \end{equation}

    It remains to estimate the above supremum from below using the quantity $\sup_{\ell\in I}\distu(q_\ell,\VV)$.
    By~\eqref{eq:dist to VV and tVV are equal} and~\eqref{eq:linearizing distance}, for every $\ell\in I$, we have 
    \begin{equation*}
        \norm{\hol_{\tq_\ell} -\pi^{\mrm{st}}_{\tVV}(\hol_{\tq_\ell})}_0
        \leq 4C_{\tK} \distm(\tq_\ell,\tVV)
        = 4C_{\tK} \distu(q_\ell,\VV)
        \leq 4C_{\tK}\d_0,
    \end{equation*}
    where the last inequality follows since $I$ is a subset of the $E(\d_0,T)$ defined in~\eqref{eq:local exceptional sets}.
    We may assume that $\d_0$ is small enough, depending on $\tK$, so that for every $\tw\in\tK$, the map $\tilde{x}\mapsto \hol_{\tilde{x}}$ is invertible on any ball of radius $4C_{\tK}\d_0$ centered around $\hol_{\tw}$ in the norm $\norm{\cdot}_0$.
    Hence, for each $\ell\in I$, there is $\tilde{x}_\ell\in \HHm$ such that
    \begin{equation*}
        \hol_{\tilde{x}_\ell} = \pi^{\mrm{st}}_{\tVV}(\hol_{\tq_\ell}).
    \end{equation*}
    Since $\pi^{\mrm{st}}_{\tVV}(\hol_{\tq_\ell})$ belongs to $V^{\mrm{st}}_{\tVV}$ (which is the image of $\tVV$ under the holonomy map $\tw\mapsto \hol_{\tw}$), it follows that $\tilde{x}_\ell$ can be chosen to belong to $\tVV$.
    Thus, by Lemma~\ref{lem:neighborhood control} and~\eqref{eq:AGY equiv to ad hoc}, for every $\ell\in I$, we have that
    \begin{align*}
        \distu(q_\ell,\VV)
        \leq \distm(\tq_\ell,\tilde{x}_\ell) 
        &\leq 4 
        \norm{\hol_{\tq_\ell} -\hol_{\tilde{x}_\ell}}_{\tq_\ell}
         \nonumber
         \leq 4 C_{\tK}
        \norm{\hol_{\tq_\ell} -\pi^{\mrm{st}}_{\tVV}(\hol_{\tq_\ell})}_0.
    \end{align*}
    This implies that
    \begin{equation}\label{eq:sup lower bound}
        P(\tq,\tVV;\ell) \geq (4C_{\tK})^{-2}  \dist^2_u(q_\ell,\VV), \qquad \forall \ell\in I.
    \end{equation}
    Combined with~\eqref{eq:linearization}, we obtain
    \begin{equation*}
        | E(\d,I,\tVV)| \leq \d |I| \cdot\frac{64 C_{\tK}^2}{\sup_{\ell\in I} \distu(q_\ell,\VV)}.
    \end{equation*}
    This completes the verification of the first assertion of the proposition in light of Claim~\ref{claim:lifts of T are locally finite}.

    To prove the second assertion, suppose that $I=[-T,T]$ so that $|I|\geq 2$.
    Hence, since $\d_0 < \eta$ by~\eqref{eq:delta0}, we have that $\distm(\tq,\tVV)<\eta$, and, thus, the supremum in~\eqref{eq:linearization} is non-zero in view of Claim~\ref{claim:coeff of tremor poly }.
    Since polynomials of degree at most $2$ form a finite dimensional vector space, all norms on this space are equivalent. In particular, there is $\d_2 = \d_2(\tK,\tVV)>0$ such that the supremum of the absolute value of any such polynomial over the interval $[-1,1]$ is at least $\d_2$ multiplied by the maximal magnitude of its coefficients.

    Recall that the polynomial $P(\tq,\tVV;\cdot)$ has leading coefficient $\norm{\pi^{\mrm{bal}}_{\tVV}(\b)}_0^2$ by Claim~\ref{claim:coeff of tremor poly }.
    By~\eqref{eq:nontrivial proj}, this coefficient is at least $\norm{\b}_0^2/4$.
    Moreover, Eq.~\eqref{eq:AGY equiv to ad hoc} implies that $\norm{\b}_0\geq \norm{\b}_q/C_{\tK}$.
    Combined with our choice of $\d_2$, we get
    \begin{equation*}
        \sup_{\ell\in I} |P(\tq,\tVV;\ell)| \geq \frac{\d_2 \norm{\b}_q^2}{4C_{\tK}^2},
    \end{equation*}
   Together with~\eqref{eq:linearization}, this estimate completes the proof of the proposition.


\section{A sufficient condition}\label{sec:general}

The goal of this section is to reduce the proof of Theorem~\ref{thm:main} to that of Theorem~\ref{thm:sufficient} below, which establishes the non-concentration of the norms of our cocycle.
This result is the main technical step in the proof of Theorem~\ref{thm:intro}. The proof of Theorem~\ref{thm:sufficient} is outlined in Section~\ref{sec:outline of oscillations} and occupies Sections~\ref{sec:oscillation}-\ref{sec:proof of oscillations}.

Recall that the balanced space is constant along $\GLp$-orbits. We denote the common balanced space over $\tilde{\Ocal}$ by $\bal{\Ocal}\subset H^1(S,\Sigma;\R)$.
Recall our fixed surface $\w_1\in\Ocal$ in~\eqref{eq:omega1} with periodic horocycle orbit of period $1$.

\begin{thm}\label{thm:sufficient}
    For all $0<\varrho<1$, there exists $t_\varrho >0$ so that for all $t\geq t_\varrho$, and $v\in \bal{\Ocal}$ with $\|v\|_{\w_1}=1$, there exist $\frac 1 {100}\leq \gamma_{t,v}\leq \frac {99} {100}$ and $C_t>0$ such that 
    \begin{equation*}
    \left| \left\{ s\in [0,1]: \|\mathrm{KZ}(g_t,u_s\w_1)v\|_{g_t u_s\omega_1}
    \geq C_t \varrho^{-1} \right\} \right|
     \geq \gamma_{t,v}-\varrho, 
    \end{equation*}    
        and 
    \begin{equation*}
        |\{s\in [0,1]: \|\mathrm{KZ}(g_t,u_s\w_1)v\|_{g_t u_s\omega_1}< {C_t}\varrho\}| \geq (1-\gamma_{t,v})-\varrho.
    \end{equation*}
    
\end{thm}

We refer the reader to~\cite{Hamid} for related results.

\subsection{Outline of the proof of Theorem \ref{thm:main} from Theorem \ref{thm:sufficient}} \label{sec:outline of main proof from oscillations}
Recall the notation preceding Theorem~\ref{thm:main}.
In Section \ref{sec:tight} we establish tightness, that is 
we show that any weak-$*$ limit of $\{\nu_{t,a}:t\geq 0\}$ is a probability measure.
For $t\geq 0$ and $C_t$ as in Theorem \ref{thm:sufficient}, let
\begin{align}\label{eq:B and S}
    B_t(\varrho) &=\left\{s\in[0,1]: \lVert \mrm{KZ}(g_t,u_s\w_1) \s \rVert_{g_t u_s\omega_1} \geq C_t \varrho^{-1}\right\},\nonumber \\
    S_t(\varrho) &= \left\{s\in[0,1]: \lVert \mrm{KZ}(g_t,u_s\w_1) \s \rVert_{g_t u_s\omega_1} <C_t\varrho\right\}.
    \end{align}
Roughly, the strategy is as follows. 
We find sequences $\varrho_k\to 0$, $L_k>0$, and $t_k\to \infty$ so that
\begin{enumerate}
    \item For all $\ell$ with $|\ell|\leq L_k$, then $\mathrm{Trem}_{g_t h_s \omega_1}(\ell)$ is very close to $\Ocal_1$ for all $s\in S_{t_k}(\varrho_k)$.

    \item Using the results of Section~\ref{sec:near teich}, we show that, for any finite collection of Teichm\"uller curves, there exists $\ell_k$ with $|\ell_k|\leq L_k$ so that for most $s\in B_{t_k}(\varrho_k)$, the point $\mathrm{Trem}_{g_t h_s \omega_1}(\ell_k)$ is a definite distance away from the chosen collection of Teichm\"uller curves.
\end{enumerate}

These results, and the fact that the proper $\GLp$-orbit closures of $\HHu$ consist of countably many Teichm\"uller curves, imply that we can choose a sequence $\nu_{t_i,a_i}$ whose weak-$*$ limit $\nu_0$
 is a non-trivial convex combination of a measure supported on $\mathcal{O}_1$ and a measure that gives zero weight to any proper $\SL$-orbit closure in $\HH_1$.
 
 Since $\nu_0$ is a limit of horocycle-invariant measures on periodic horocycles, it is horocycle-invariant as well.
 Moreover, it can be shown using by-now standard arguments that the part of the measure that lives on $\Ocal_1$ must be $\mu_{\Ocal}$ (note however that these periodic horocycles are not contained in $\Ocal$).
 
 However, we take a different approach that shows that both measures are $\SL$-invariant.
 To do so, in Section \ref{sec:soft}, we use results of Eskin-Mirzakhani-Mohammadi~\cite{EMM}, through a result of Forni~\cite{ForniDensity1}, to show that further pushing our weak-$\ast$ limit measure $\nu_0$ by the geodesic flow (along a subsequence) gives an $\SL$-invariant measure $\nu$ in the limit.
 This limiting measure must be a convex combination of $\mu_\Ocal$ and a measure which gives zero weight to any proper $\SL$-orbit closure, i.e.~$\mu_{\mrm{MV}}$.
Note this summary leaves out some issues in the proof and in particular, some of the statements we made require restricting to suitable compact sets.

\subsection{Accumulation on the octagon locus}\label{sec:hug octagon}

Recall the balanced tremor $\sigma\in\TT_{\w_1}^{\mrm{bal}}$ defined in Section~\ref{sec:octagon-tremors}.
We may assume it is normalized so that $\norm{\s}_{\w_1}=1$. Recalling that $\TT_{\w_1}^{\mrm{bal}}$ is contained in the horizontal space, we let
\begin{equation*}
    \sigma_t=\mrm{KZ}(g_tu_s,\w_1)\s = \mrm{KZ}(g_t,u_s\w_1)\s,
\end{equation*}
 where the action of the cocycle on the horizontal space is defined in~\eqref{eq:cocycle on horizontal}. The second equality follows from the cocycle property and~\eqref{eq:horocycle fixes tremor}.
 To simplify notation, we drop the dependence on $s$ in $\s_t$.
    Observe that 
    $\nu_{t,a}$, the $u_s$-invariant probability measure on the periodic $u_s$-orbit through $g_t\mathrm{Trem}_{\omega_1,\s}(a)$ is given by

  \begin{align}\label{eq:def. nu}
        \nu_{t,a}(\mathcal{A})&:=|\{s\in[0,1]:g_t\mathrm{Trem}_{u_s\omega_1,\sigma}(a)\in \mathcal{A}\}|
        \nonumber \\
        &= 
        |\{s\in [0,1]:\mathrm{Trem}_{g_t u_s\omega_1,\sigma_t}(e^{t}a)\in \mathcal{A}\}|
    \end{align}    
    for all measurable $\mathcal{A}\subset \HH_1$, where we used Lemma~\ref{lem:renorm} for the second equality.
    Here and throughout the remainder of the article, we continue to drop the trivialization maps from our notation. 
    
    We record the following basic fact about limits of $\nu_{t,a}.$
    \begin{lem}\label{lem:hor invar} 
    Let $t_i$ be a sequence going to infinity and let $a_i$ be an arbitrary sequence. Then, any  weak-$*$ limit of $\nu_{t_i,a_i}$ is horocycle-invariant.
    \end{lem}
    \begin{proof}
    This is because the set of invariant measures for a (continuous) flow is closed in the weak-$*$ topology and each of the $\nu_{t,a_i}$ is horocycle-invariant (in fact given by a periodic horocycle orbit). 
    \end{proof}

    Given $\varrho>0$, let $t_\varrho>0$ be the constant provided by Theorem~\ref{thm:sufficient}. Recall the sets $B_t(\varrho)$ and $S_t(\varrho) $ defined in~\eqref{eq:B and S}. 
    
       \begin{equation}\label{eq:nu'}
      \nu'_{t,a}(\varrho)(\mathcal{A})=\frac{1}{|B_t(\varrho)|} |\{s\in B_t(\varrho):g_t\mathrm{Trem}_{u_s\omega_1,\sigma}(a)\in \mathcal{A}\}|. 
    \end{equation}
    
    Similarly, let
        \begin{equation}\label{eq:nu''}
        \nu''_{t,a}(\varrho)(\mathcal{A})= \frac{1}{|S_t(\varrho)|} |\{s\in S_t(\varrho):g_t\mathrm{Trem}_{u_s\omega_1,\sigma}(a)\in \mathcal{A}\}|.
    \end{equation}
    
   The results of the previous subsection imply the following corollary via Fubini's theorem.
    \begin{cor}\label{cor:no teich} 
    
    Given $0<\varrho<1$, let $t_\varrho>0$ and $C_{t_\varrho}>0$ be the constants provided by Theorem~\ref{thm:sufficient}.
    Then, for every $0<\varrho<1$, there exists $a_\varrho$ with $|a_\varrho|\leq e^{-t_\varrho}/C_{t_\varrho}$ such that the following holds. 
    For every sequence $\varrho_k\to 0$ such that the measures $\nu'_{t_{\varrho_k},a_{\varrho_k}}$ converge to a measure $\nu'_\infty$ in the weak-$\ast$ topology, we have that $\nu'_\infty$ gives zero mass to all Teichm\"uller curves in $\HH_1$.
    
    \end{cor}

    The idea of the proof is as follows.
    We consider the $2$-parameter family of surfaces $g_{t_\varrho}\mrm{Trem}_{ u_s\w_1, \sigma}(\ell)$
    parametrized by $(s,\ell)\in B_{t_\varrho}(\varrho)\times \R$. Proposition~\ref{prop:near teich} says that for a fixed $s$, tremor orbits (corresponding to $\set{s}\times \R$) do not concentrate near any Teichm\"uller curve. The corollary will follow by an application of Fubini's theorem. We use that $s\in B_{t_\varrho}(\varrho)$ to ensure that the tangent vector to the tremor orbit has a definite size. Indeed, the statement does not hold for $s\in S_{t_\varrho}(\varrho)$, cf. Lemma~\ref{lem:mass on octagon}.
    
\begin{proof}[Proof of Corollary~\ref{cor:no teich}] By~\cite{McMRig}, any stratum of abelian differentials contains countably many Teichm\"uller curves.
    In genus two, this fact is due to~\cite{McM,Ca}.
    Let $\VV_1,\VV_2,\dots$ be an enumeration of the Teichm\"uller curves in $\HHu$.
    Let $\Kcal_1\subset \Kcal_2\subset\dots$ be an exhaustion of $\HHu$ by compact sets with non-empty interior.
    For every $n\geq 1$, let $P_n = \bigcup_{k=1}^n\VV_k$ and for each $\d>0$, denote by $P_n^\d$ the union of the $\d$-neighborhoods of $\VV_k$ for $1\leq k\leq n$.
    Denote by $\chi_n^\d$ the indicator function of $P_n^\d\cap\Kcal_n$.
    For $\varrho>0$ and $s\in\R$, let 
\begin{equation*}
    x(\varrho,s)=g_{t_\varrho}u_s\omega_1, \qquad \sigma(\varrho,s)=\mathrm{KZ}(g_{t_\varrho},u_s\w_1)\s, \qquad
     \bar{\sigma}(\varrho,s)=\frac{\sigma(\varrho,s)}{\norm{\sigma(\varrho,s)}}_{x(\varrho,s)}.
\end{equation*}
    
    We claim that Proposition~\ref{prop:near teich} implies that for each $n\geq 1$, we can find $\d_n>0$ so that for all $\varrho>0$,
    \begin{align}\label{eq:vertical lines}
       2C_{t_\varrho} \int_{|\ell|\leq C_{t_\varrho}^{-1}} \chi_n^{\d_n}\left( \mrm{Trem}_{x(\varrho,s), \sigma(\varrho,s)}(\ell) \right) \;d\ell \leq 1/n, \qquad \forall s\in B_{t_\varrho}(\varrho).
    \end{align}
    Indeed, note that for all $q\in\HHu$, $0\neq \b\in\TT^{\mrm{bal}}_q$ and $\ell\in\R$,
    \begin{equation*}
        \mrm{Trem}_{q,\b}(\ell) = \mrm{Trem}_{q,\frac{\b}{\norm{\b}_q}}(\ell\norm{\b}_q).
    \end{equation*}
    For simplicity, we will use $\norm{\cdot}$ to denote $\norm{\cdot}_{x(\varrho,s)}$.
    By a change of variable, we obtain
    
    \begin{align*}
        \int_{|\ell|\leq C_{t_\varrho}^{-1}} \chi_n^{\d_n}\big( \mrm{Trem}_{x(\varrho,s), \sigma(\varrho,s)}&(\ell) \big)\;d\ell \\
         &= \int_{|\ell|\leq C_{t_\varrho}^{-1}} \chi_n^{\d_n}\left( \mrm{Trem}_{x(\varrho,s), \bar{\sigma}(\varrho,s)}(\ell \norm{\sigma(\varrho,s)}) \right) \;d\ell \\
         &= \frac{1}{\norm{\sigma(\varrho,s)}} \int_{|\ell|\leq\frac{\norm{\sigma(\varrho,s)}}{C_{t_\varrho}}} \chi_n^{\d_n}\left( \mrm{Trem}_{x(\varrho,s), \bar{\sigma}(\varrho,s)}(\ell ) \right) \;d\ell.
    \end{align*} 
    Let $T = 
    \frac{\norm{\sigma(\varrho,s)}}{C_{t_\varrho}}$
     and note that since $s\in B_{t_{\varrho}}(\varrho)$, we have that $T\geq 1/\varrho >1$. 
    We may then apply Proposition~\ref{prop:near teich} for each $\VV_k$, $1\leq k\leq n$, with this choice of $T$ and with $\e = 1/n^2$, $q= x(\varrho,s)$, $\b = \bar{\sigma}(\varrho,s)$, and $\Kcal=\mrm{interior}(\Kcal_n)$ to obtain $\d_n>0$ which satisfies~\eqref{eq:vertical lines}.
    Note that we are using Lemma~\ref{lem:renorm} to ensure that $\bar{\sigma}(\varrho,s) $ belongs to the tremor space of $x(\varrho,s)$.
    
    In view of~\eqref{eq:vertical lines}, by Fubini's theorem,
\begin{equation*}
        2C_{t_\varrho}  \int_{|\ell|\leq C^{-1}_{t_\varrho}}
        \Big(\frac{1}{|B_{t_\varrho}(\varrho)|}
        \int_{B_{t_\varrho}(\varrho)} \chi_n^{\d_n}\left( \mrm{Trem}_{x(\varrho,s), \sigma(\varrho,s)}(\ell) \right)\; ds\Big) d\ell \leq 1/n.
    \end{equation*}

    Hence, for each $\varrho$ with $2^{-n} \leq \varrho < 2^{-n+1}$, we can find $\tilde{a}_\varrho$ with $|\tilde{a}_\varrho| \leq C^{-1}_{t_\varrho}$ such that the inner average (viewed as a function of $\ell$) is at most $1/n$.
    Let $a_{\varrho} = e^{- t_\varrho} \tilde{a_{\varrho}}$ and note that

  \begin{equation*}
        g_{t_\varrho} \mrm{Trem}_{u_s\w_1, \sigma}(a_{\varrho}) = 
        \mrm{Trem}_{x(\varrho,s), \sigma(\varrho,s)}(\tilde{a}_{\varrho}) 
    \end{equation*}
   by Lemma~\ref{lem:renorm}.
   One then checks that such a choice of $a_\varrho$ satisfies the corollary.

\end{proof}

    For $\tq\in\HHm$, denote by $E^{\mrm{u}}(\tq)\subset \mrm{T}_{\tq}$ the unstable subspace of the tangent space for the Teichm\"uller geodesic flow.
    More explicitly, we let $E^{\mrm{u}}(\tq)$ denote the horizontal space $\mrm{H}^1(S,\Sigma;\R_x)$, viewed as a subspace of the tangent space $\mrm{T}_{\tq}$ under the identification
    $\mrm{T}_{\tq}\cong \mrm{H}^1(S,\Sigma;\R_x)\oplus \mrm{H}^1(S,\Sigma;\R_y)$; cf. Section~\ref{sec:background} for definitions.

    Avila-Gou\"ezel defined an analogue of the exponential map, denoted by $\Psi_q$, from a neighborhood of $0$ in $E^{\mrm{u}}(\tq)$ to $\HHm$ as follows. Given a path $\k:[0,1]\to \HHm$ with $\k(0)=\tq$ such that $\k'(t)=v$ for all $t$, one defines $\Psi_{\tq}(v)=\k(1)$.
    It is shown in~\cite[Proposition 5.3]{AG} that for all $\tq\in\HHm$, $\Psi_{\tq}$ is well-defined on a ball of radius $1/2$ in $E^{\mrm{u}}(\tq)$ in the sup-norm.
    It also follows by~\cite[Proposition 5.3]{AG} that
    \begin{equation}\label{eq:bound distance by norm}
        \distm(\tq,\Psi_{\tq}(w)) \leq 2\norm{w}_{\tq},
    \end{equation}
    for all $w\in E^{\mrm{u}}(\tq)$ with $\norm{w}_{\tq}\leq 1/2$.
    By the description of tremor in holonomy coordinates in~\eqref{eq:tremor holonomy}, it follows that for a path given by $\k(t) = \mrm{Trem}_{\tq,\b}(t)$, we have $\k'(t)=\b$ for all $t$.
    In particular, for $\b\in \TT_{\tq}^{\mrm{bal}}$ with $\norm{\b}_{\tq}\leq 1/2$, we have
    \begin{equation}\label{eq:tremor is exponential map}
        \Psi_{\tq}(\b) = \mrm{Trem}_{\tq,\b}(1).
    \end{equation}

    \begin{lem}\label{lem:mass on octagon}
     For each $\varrho>0$, let $t_\varrho>0$ and $C_{t_\varrho}>0$ be the constants provided by Theorem~\ref{thm:sufficient}. Let $a_\varrho$ be the constant provided by Corollary~\ref{cor:no teich}.
    Let $\varrho_k\to 0$ be an arbitrary sequence such that the measures $\nu''_{t_{\varrho_k},a_{\varrho_k}}(\varrho_k)$ converge to a measure $\nu_\infty''$ in the weak-$\ast$ topology. Then, $\nu_\infty''$ 
    is $u_s$-invariant and satisfies $\nu_\infty''(\HH_1\backslash\Ocal_1)=0$.
    \end{lem}

    \begin{proof}[Proof of Lemma~\ref{lem:mass on octagon}]
    Let $\nu''_\infty$ be one such limit measure along a sequence $\varrho_n\to 0$. We first show that $\nu''_\infty(\HH_1 \backslash \Ocal_1) = 0$.
    Recall that $\distu$ and $\distm$ refer to the sup-norm metrics on $\HHu$ and $\HHm$ respectively.

    Let $\varrho\leq 1/2$, $s\in S_{t_\varrho}(\varrho)$ and denote $g_{t_\varrho}u_s\w_1$ by $x(\varrho,s)$.
    Let $\tilde{x}(\varrho,s)$ denote its unique lift to our fixed fundamental domain in $\HHm $.
    Let

  \begin{equation*}
        T(\varrho,s):= g_{t_\varrho} \mrm{Trem}_{u_s\w_1,\sigma}(a_\varrho) = \mathrm{Trem}_{x(\varrho,s),\sigma_{t_\varrho}}(e^{{t_\varrho}}a_\varrho) \in \HH_1,
    \end{equation*}
    where the second equality follows by  Lemma~\ref{lem:renorm}.
    Define $\tilde{T}(\varrho,s)\in \HH_{\mrm{m},1}$ similarly using $\tilde{x}(\varrho,s)$ in place of $x(\varrho,s)$.
    Note that $\tilde{T}(\varrho,s)$ is a lift to $\HH_{\mrm{m},1}$ of $T(\varrho,s)$.
    Hence, we have 
    \begin{align*}
    \mrm{dist}_u(x(\varrho,s),T(\varrho,s)) \leq \mrm{dist}_m(\tilde{x}(\varrho,s), \tilde{T}(\varrho,s)).
    \end{align*}
    As $s\in S_t(\varrho)$ and $e^{{t_\varrho}}a_\varrho\leq C_{t_\varrho}^{-1}$ (by Corollary~\ref{cor:no teich}), we have $\norm{e^{t_\varrho}a_\varrho\s_{t_\varrho}}_{\tilde{x}(\varrho,s)}\leq\varrho$.
    Since $\varrho \leq 1/2$, we may apply with~\eqref{eq:bound distance by norm} and~\eqref{eq:tremor is exponential map} to obtain

    \begin{align*}
        \mrm{dist}_u(x(\varrho,s),T(\varrho,s))
        \leq 2 C_{t_\varrho}^{-1} \norm{\sigma_{t_\varrho}}_{\tilde{x}(\varrho,s)} 
        < 2 \varrho.
    \end{align*}
    The above being true for all $s\in S_{t_\varrho}(\varrho)$ and since $x(\varrho,s)\in\Ocal_1$, it follows that $\nu''_\infty(\HH_1\backslash\Ocal_1)=0$.

    Since $\nu_\infty''$ lives on $\Ocal_1$, to show that it is $u_s$-invariant, it suffices to prove that $\nu_\infty''(u_s A) = \nu_\infty''(A)$ for all Borel sets $A\subseteq \Ocal_1$. 
    Since the space of Borel measures of mass $\leq 1$ on $\HH_1$ is compact, by passing to a subsequence we may assume that $\nu_{t_{\varrho_n},a_{\varrho_n}}$ converges to a measure $\nu_\infty$.
    By Lemma~\ref{lem:hor invar}, $\nu_\infty$ is $u_s$-invariant.
    Moreover, since $\varrho_n\to 0$, Theorem~\ref{thm:sufficient} shows that $\nu_\infty$ is a (non-trivial) convex combination of $ \nu_\infty' $ and $\nu_\infty''$, where $\nu_\infty'$ is some weak-$\ast$ limit of the measures $\nu'_{t_{\varrho_n},a_{\varrho_n}}$.
    
    Let $A\subseteq \Ocal_1$. 
    By Corollary~\ref{cor:no teich}, we have that $\nu_\infty'(\Ocal_1) =0$.
    Since $\Ocal_1$ is $u_s$-invariant, we have $u_sA\subseteq \Ocal_1$.
    Hence, we get
    \[
        \nu_\infty''(u_s A) = \nu_\infty(u_s A) = \nu_\infty(A) = \nu''_\infty(A).
        \qedhere
    \]
    \end{proof}

    Lemma~\ref{lem:mass on octagon} implies the following corollary.
    \begin{cor}\label{cor:nu' horocycle invariant}
    Let $\nu_\infty'$ be a limit measure as in Corollary~\ref{cor:no teich}. Then, $\nu_\infty'$ is $u_s$-invariant.
    \end{cor}
    \begin{proof}
    The proof is completely analogous to the argument at the end of the proof of Lemma~\ref{lem:mass on octagon} and relies on the fact that $\nu_\infty'$ lives on the complement of all Teichm\"uller curves which is a $u_s$-invariant set as follows by Corollary~\ref{cor:no teich}.
    \end{proof}

\subsection{Non-escape of mass}\label{sec:tight}

    We show that the collection of measures $\nu_{t,a}$ constructed above is tight.
    \begin{prop}
    \label{prop:compact} 
        For all $\varrho>0$, there exists a compact set $\Kcal \subset \HHu$ (depending on $\varrho$) so that for all $t\geq 0,a\in \mathbb{R}$, we have 
        \begin{equation*}
            \nu_{t,a}(\Kcal) \geq 1-\varrho.
        \end{equation*}
    \end{prop}

    We deduce Proposition~\ref{prop:compact} from the following lemma, which is due to~\cite{EskinMasur}; cf.~\cite[Lemma 3.5]{AAEKMU} for the version below.
    
  \begin{lem}
  \label{lem:integrability of the height function}
    There exists a proper function $\alpha:\HHu \to [0,\infty)$, $t_0>0$ and $b$ so that $\int_{0}^1\alpha(g_tu_sx)ds\leq \alpha(x)+b$ for all $x \in \HHu$ and $t\geq t_0$.  
  \end{lem}

\begin{proof}[Proof of Proposition~\ref{prop:compact}]
    Recall that $\mathcal{C}:=\{u_s\omega_1\}_{s\in [0,1)}$ and let $S$ be the set of tremors of $\mathcal{C}$. 
     Then $S$ is the image of a two-torus 
    under translation equivalence. Indeed, we can consider the measure $\tau'$ giving full weight to the cylinder of area $a$ on $\omega_1$. Because a twist of one of the cylinders will eventually return it to its initial position, there exists a minimal $r>0$ so that $\mathrm{Trem}_{\omega_1,\tau'}(r)=\omega_1$. By definition $u_1\omega_1=\omega_1$. Observe that the tremor of any point in $\mathcal{C}$ can be written as $u_s \mathrm{Trem}_{\omega_1,\tau'}(\ell)$ for some $s,\ell$. By Lemma~\ref{lem:commute}, tremors and the horocycle flow commute.
    Hence, we have that $S$ is the  two-torus (with fundamental domain $[0,1] \times [0,r]$) or its image under quotienting out by translation equivalence.  
     Thus, there is a fixed compact set, $K'$, so that $S\subset K'$.

    Let $\alpha$ be a function satisfying Lemma~\ref{lem:integrability of the height function}.
    Since $S$ is compact, we have $A=\max_{\omega \in S}\alpha(\omega)<\infty$. Given $\varrho$, let $\Kcal=\alpha^{-1}[0,\frac{2(A+b)}\varrho]$. So for every $x \in S$ and $t$, we have 
    \[
    \left|\{s\in [0,1]:g_t u_sx\notin \Kcal\}\right|\leq \frac\varrho2. \qedhere \] 
\end{proof}

\subsection{Proof of Theorem~\ref{thm:main} and Corollary~\ref{cor:nongen}}\label{sec:soft}

    As a first step to proving our main results, by a straightforward diagonal argument we have: 
\begin{lem}\label{lem:closed}
The set of measures that arise as weak-$*$ limits of $\nu_{t_i,a_i}$ for any choice of $t_i,a_i$ is closed in the weak-$*$ topology. 
\end{lem}
Since $g_\ell \nu_{t,a}=\nu_{\ell+t,e^{\ell}a}$, we have: 
\begin{lem}\label{lem:g push} The set of probability measures that can be obtained as weak-$*$ limits of $\nu_{t_i,a_i}$ is closed under the pushforwards $(g_\ell)_\ast$ for any $\ell$.
\end{lem}

    To upgrade the horocycle invariance of our measures to invariance under all of $\SL$, we need the following result which is deduced from the work of Eskin-Mirzakhani~\cite{EM} and Eskin-Mirzakhani-Mohammadi~\cite{EMM} via a result of Forni~\cite{ForniDensity1}.

    \begin{thm}
    \label{thm:g limit} 
    Let $\mu$ be a horocycle-invariant measure and $\mathcal{M}$ an $\mrm{SL}_2(\mathbb{R})$-orbit closure so that $\mu(\mathcal{M})=1$ and 
    for any $\mrm{SL}_2(\mathbb{R})$-orbit closure $\mathcal{M}'\subsetneq \mathcal{M}$, $\mu(\mathcal{M}')=0$. Then, there exists an unbounded set $S \subset \mathbb{R}_+$ so that 
    \begin{equation*}
        \lim_{t\in S, t\to\infty}  (g_t)_\ast\mu = \mu_{\mathcal{M}},
    \end{equation*}
     where $\mu_{\mathcal{M}}$ is the unique $\mrm{SL}_2(\mathbb{R})$-invariant Lebesgue measure whose support is $\mathcal{M}$.
    \end{thm}

    \begin{proof}
    Note that since $\mu_{\mathcal{M}}$ is $\mrm{SL}_2(\mathbb{R})$-ergodic, it is ergodic for the horocycle flow by Mautner's phenomenon.
    Hence, by \cite[Theorem 1.1]{ForniDensity1}
    it suffices to show that $\frac 1 T \int_0^T(g_t)_\ast\mu \;dt$ converges to $\mu_{\mathcal{M}}$. This follows by Eskin-Mirzakhani and Eskin-Mirzakhani-Mohammadi. To see this, we claim that: \begin{equation} \label{eq:meas push}
    \mu_{\mathcal{M}}= \lim_{T\to\infty} \frac 1 T \int_0^T(g_t)_\ast\mu \;dt.
    \end{equation} First, \cite[Theorem 2.7]{EMM} shows 
    $$\mu_{\mathcal{M}}=\lim \frac 1 T \int_0^T\int_0^1 \d_{g_tu_sx} ds dt$$
    for any $x\in \mathcal{M}$ that is not contained in a closed $\mrm{SL}_2(\mathbb{R})$-invariant sublocus of $\mathcal{M}$. 
    This establishes the result for $U$-invariant measures, whose support is contained in $\mathcal{M}$ and that give zero mass to the union of all closed $\mrm{SL}_2(\mathbb{R})$-invariant sets contained in $\mathcal{M}$. 
    By \cite[Proposition 2.16]{EMM}, there are only countably many of these, and so this is implied by the in principle weaker assumption that $\mu$ gives zero mass to each closed $\mrm{SL}_2(\mathbb{R})$-invariant sublocus of $\mathcal{M}$. We have established \eqref{eq:meas push}. 
    \qedhere

    \end{proof}  

    We are now ready for the proof of the main theorem.
\begin{proof}[Proof of Theorem \ref{thm:main} assuming Theorem~\ref{thm:sufficient}] 
    Recall the measures defined in equations~\eqref{eq:def. nu},~\eqref{eq:nu'}, and~\eqref{eq:nu''}.
    Given $\varrho \in (0,1)$, we denote by $t_\varrho>0$ and $C_{t_\varrho}>0$ the constants provided by Theorem~\ref{thm:sufficient}. Let $a_\varrho $ be a constant satisfying Corollary~\ref{cor:no teich} so that $|a_\varrho|\leq  e^{-t_\varrho}/C_{t_\varrho}$.
    Since the space of Borel measures of mass at most $ 1$ on $\HH$ is compact, we can find a sequence $\varrho_n\to 0$ such that the sequence $\nu_{t_{\varrho_n},a_{\varrho_n}}$ converges to a measure $\nu_\infty$.
    By Lemma~\ref{lem:hor invar}, $\nu_\infty$ is $u_s$-invariant. 
    By Proposition~\ref{prop:compact}, $\nu_\infty$ is a probability measure.
    
    By passing to a subsequence, we may assume that $\nu'_{t_{\varrho_n},a_{\varrho_n}}$ and $\nu''_{t_{\varrho_n},a_{\varrho_n}}$ converge to measures $\nu'_\infty$ and $\nu''_\infty$ respectively.
    By Theorem~\ref{thm:sufficient}, 
    \begin{equation*}
        \nu_\infty = p \nu_\infty' + (1-p) \nu_\infty'',
    \end{equation*}
    for some $1/100\leq p \leq 99/100$.
    
    Since $\nu_\infty$ is a probability measure, both $\nu_\infty'$ and $\nu_\infty''$ are probability measures.
    By Lemma~\ref{lem:mass on octagon}, $\nu_\infty''$ is $u_s$-invariant and $\nu_\infty''(\Ocal)=1$.
    It is then a well-known result (cf.~\cite{KleinbockMargulis-quasiunipotent}) that
    \begin{equation*}
        (g_{t})_\ast \nu_\infty'' \xrightarrow{t\to\infty} \mu_{\Ocal},
    \end{equation*}
   where $\mu_{\Ocal}$ is the Haar probability measure on $\Ocal$.
   In view of Corollary~\ref{cor:no teich}, $\nu_\infty'$ gives $0$ mass to all Teichm\"uller curves inside $\HHu$.
   By the work of Calta~\cite{Ca} and McMullen~\cite{McM}, the only proper $\SL$-invariant subloci of $\HHu$ are Teichm\"uller curves.
   Hence, Theorem~\ref{thm:g limit} shows that there is a sequence $t_k\to\infty$ such that
   \begin{equation*}
       (g_{t_k})_\ast \nu_\infty' \xrightarrow{k\to\infty} \mu_{\mrm{MV}},
   \end{equation*}
   where $\mu_{\mrm{MV}}$ is the Masur-Veech measure on $\HH_1$.
   It follows that 
   \begin{equation*}
       (g_{t_k})_\ast \nu_\infty \xrightarrow{k\to\infty} p\mu_{\mrm{MV}} + (1-p) \mu_{\Ocal}.
   \end{equation*}
   To conclude the proof, we note that Lemmas~\ref{lem:closed} and~\ref{lem:g push} show that the above convex combination can be realized as a limit of of measures of the form $\nu_{t_i,a_i}$ (for a now  possibly unbounded sequence of $a_i$'s).

\end{proof}


\section{Oscillations of the KZ cocycle}\label{sec:oscillation}
    
    The goal of this section is to reduce the proof of Theorem~\ref{thm:sufficient} to Proposition~\ref{prop:partitions spread} below.
    We also introduce several preliminary results which we need for the proof of Proposition~\ref{prop:partitions spread}.

    \subsection*{Notational convention}
    Throughout the remainder of the article, we make an identification
    \begin{equation*}
        \bal{\Ocal} \cong \R^2.
    \end{equation*} 
    by choosing a basis for $\bal{\Ocal}$.
    We may then view the restriction of the KZ-cocycle to $\bal{\Ocal}$ as taking values in $\SL$.
    In the remainder of the article, we will use the same notation $\mrm{KZ}(\cdot,\cdot)$ for this restriction of the cocycle to $\bal{\Ocal}$.
    Additionally, for convenience, we fix a norm on $\R^2$ and denote it $\norm{\cdot}$.
    A convenient explicit choice of such basis will be made in Section~\ref{sec:pseudoanosov}.

\begin{prop}\label{prop:partitions spread}
    For all $0<\k<1$, there exists $t_0>0$ such that for all $t\geq t_0$ and all $v\in \R^2$ with $\norm{v}=1$, 
    \begin{equation*}
        \sup_{r\geq 0} \left|\set{s\in [0,1]: \|\KZ{g_t}{u_s \omega_1)v}\|\in \left[ \k r,\k^{-1} r \right] }\right|<\frac{49}{50} +\k.
    \end{equation*}

\end{prop}

    The proof of Proposition~\ref{prop:partitions spread} is given in Section~\ref{sec:proof of oscillations}. The main intermediate results needed for the proof are stated in this section and proved in Sections~\ref{sec:pseudoanosov} and~\ref{sec:flags}.
    
\subsection{Proof of Theorem~\ref{thm:sufficient} from Proposition~\ref{prop:partitions spread}}

Fix some $\varrho\in (0,1)$.
 Note that it suffices to establish Theorem~\ref{thm:sufficient} with our fixed norm $\norm{\cdot}$ in place of the sup-norm.
    Indeed, using Proposition~\ref{prop:horocycle recurrence} below, we can choose a compact set $F\subset\Ocal$, depending on $\varrho$, so that for all $s\in[0,1]$ outside of a set of measure at most $\varrho/2$, we have $g_tu_s\w_1\in F$.
    On $F$, the sup-norm is uniformly equivalent to our fixed norm $\norm{\cdot}$.

Moreover, it suffices to show that for every $\k>0$, there exists $t_\k >0$ so that for all $t\geq t_\k$, and $v\in \mathbb{R}^2$ with $\|v\|=1$,
    \begin{equation}\label{eq:non-iterated big}
    \left| \left\{ s\in [0,1]: \|\mrm{KZ}(g_t,u_s\omega_1)v\|
    \geq  C'_t \k^{-1} \right\} \right|
    \geq 1/100-\k, 
    \end{equation}    
        and 
    \begin{equation}\label{eq:non-iterated small}
        |\{s\in [0,1]: \|\mrm{KZ}(g_t,u_s\omega_1)v\|< C'_t \k \}|\geq 1/100-\k,
    \end{equation}
    for some $C'_t >0$.
    To see that this implies the assertion of Theorem~\ref{thm:sufficient}, let $N\in\N$ be a large integer to be chosen below depending on $\varrho$, and apply~\eqref{eq:non-iterated big} and~\eqref{eq:non-iterated small} with $\k=\varrho^{2N}$.
    Since the interval $[\k C'_t, \k^{-1}C'_t)$ is a disjoint union of $2N-1$ intervals of the form $[\varrho^{2k+2} C'_t, \varrho^{2k}C'_t) $ for $-N\leq k\leq N-2$, the pigeonhole principle, along with~\eqref{eq:non-iterated big} and~\eqref{eq:non-iterated small}, implies that
    \begin{align*}
        |\{s\in [0,1]: \|\mrm{KZ}(g_t,u_s\omega_1)v\|\in [ \varrho^{2k+2} C'_t, \varrho^{2k}C'_t)\}|\leq \frac{49/50 +2\k}{ 2N-1},
    \end{align*}
    for some $-N\leq k\leq N-2$.
    Hence, if $N$ is large enough, depending on $\varrho$, the above bound is at most $\varrho/2$.
    In light of~\eqref{eq:non-iterated big} and~\eqref{eq:non-iterated small}, Theorem~\ref{thm:sufficient} now follows by taking $C_t = \varrho^{2k+1}C'_t$, where $k$ is chosen so the above estimate holds.

    To show the estimates~\eqref{eq:non-iterated big} and~\eqref{eq:non-iterated small}, let $t_0>0$ be the constant provided by Proposition~\ref{prop:partitions spread}.
    Let $\l$ be the Borel probability measure on the real line given by 
     \begin{equation*}
       \l(A) = \left|\set{s\in [0,1]: \|\KZ{g_t}{u_s \omega_1)v}\|\in A}\right|,
     \end{equation*}
     for every measurable set $A$.
     Note that $\l$ is supported on the half line $[0,\infty)$.
     Let $F_\l$ be the function defined by $F_\l(r) = \l( [r,\infty))$, for every $r\in \R$, and
     let $r_0=\sup\{r\geq 0:F_\l(r) > 1/100-\k\}$.
     Proposition~\ref{prop:partitions spread} implies that $\l(\set{0})<49/50 +\k$, and, hence, $r_0$ is strictly positive.
     Moreover, from outer regularity of Borel measures, we deduce that the half closed interval $[r_0,\infty)$ has measure at least $ 1/100-\k$.
     This implies~\eqref{eq:non-iterated big} with $C'_t = \k r_0$.
     On the other hand, inner regularity implies that the open half interval $(r_0,\infty)$ has measure at most $ 1/100-\k$.
     Hence, by Proposition \ref{prop:partitions spread}, we have that $F_\l(r_0\k^2) = \l([r_0\k^2,r_0]) + \l((r_0,\infty)) \leq  1/100-\k+49/50+\k$. Thus,~\eqref{eq:non-iterated small} also holds for the same choice of $C'_t$.   
    
    \subsection{Flow boxes and local holonomy}
    \label{sec:flow charts}
    
    For $\e>0$ and a closed connected subgroup $H\subseteq \SL$, we denote by $H_\e$ the $\e$-neighborhood of identity in $H$.
    Denote by $A$, $U$, and $U^-$ the subgroup of diagonal, upper triangular, and lower triangular matrices of $\SL$ respectively.
    The product map gives a local diffeomorphism $U^-\times A\times U\r \SL$ with a Zariski-dense open image.
    Hence, there exists a contant $c_0\geq 1$ such that the image of $U^-_\e\times A_\e \times U_\e$, which we denote $B_\e$, is contained in the $c_0\e$-neighborhood of identity in $\SL$.

    Given $x\in \Ocal$ and $0<\e\leq \inj(x)/c_0$, then the map $g\mapsto gx$ embeds $B_\e$ isometrically inside $\Ocal$.
    We write $B_\e x$ for the image of this embedding.
    Hence, for every $y= \hat{u}_s g_t u_rx\in B_\e$, we can write
    \begin{equation*}
        \hat{u}(y) := s, \qquad a(y):= t, \qquad u(y) :=r.
    \end{equation*}
    Note that the dependence on $\e<\inj(x)/c_0$ is surpressed.
    In particular, $\hat{u},a, $ and $u$ give coordinates on $B_\e x$ so that we may use $(s,t,r) \in B_\e x$ to denote the point $\hat{u}_s g_tu_r x$.
    We refer to these coordinates as \emph{flow adapted coordinates}.
    
    Given $y\in B_\e x$, we denote by $W^{\mrm{u}}_{\mrm{loc}}(y)$ the local unstable leaf of $y$ inside $B_\e x$.
    More precisely, 
    \begin{equation*}
        W^{\mrm{u}}_{\mrm{loc}}(y) = \set{z\in B_\e x: \hat{u}(z)=\hat{u}(y), a(z)=a(y) }.
    \end{equation*}
    The weak stable leaf through $y$, denoted $\Wcs(y)$, consists of those points $z\in B_\e y$ with $u(z)=u(y)$.
    
    Given $y,z\in B_\e x$, the weak stable holonomy, denoted $\Psi^{\mrm{cs}}_{y,z}:\Wu(y)\r \Wu(z)$ is defined as follows: for all $y'\in \Wu(y)$, $\Psi_{y,z}^{\mrm{cs}}(y') $ is defined to be the unique point in $\Wu(z)\cap \Wcs(y')$.
    Whenever $\e>0$ is small enough 
    (independently of $x$ so long as $\Psi_{\bullet}^{\mrm{cs}}$ is defined), the maps $\Psi_{\bullet}^{\mrm{cs}}$ are absolutely continuous with respect to the Lebesgue measure on $\Wu$.
    Moreover, the Jacobians of $\Psi_{\bullet}^{\mrm{cs}}$ tend to $1$ uniformly as $\e \r 0$.
    More precisely, for every $\d>0$, there is $\e_0>0$ so that for all $0<\e<\e_0$ and all $x\in \Ocal$ with $\inj(x)/c_0>\e$, the Jacobians of $\Psi_\bullet^{\mrm{cs}}$ in the flow box $B_\e x$ are within $\d$ from $1$.
    These facts follow readily for instance from the following computation:
    \begin{equation*}
        u_s\hat{u}_r=\begin{pmatrix}1&0\\ \frac{r}{1+s r}&1\end{pmatrix}\begin{pmatrix} 1+s r &0\\0& \frac{1} {1+s r}\end{pmatrix}
    \begin{pmatrix}1&\frac{s}{1+s r}\\0&1\end{pmatrix}.
    \end{equation*}
    Indeed, if $p^-=\hat{u}_rg_t \in  U^-_\e A_\e$ satisfies $y=p^- z$, and $u_s\in U_\e$ is such that $u_s y\in W_{\mrm{loc}}^{\mrm{u}}(y)$, then the above computation shows that $\Psi_{y,z}^{\mrm{cs}}(u_s y)= u_{e^{-2t}s/(1+sr)}z$. In particular, in flow adapted coordinates, the Jacobian of $\Psi_{y,z}^{\mrm{cs}}$ is the Jacobian of the map $s\mapsto e^{-2t}s/(1+sr)$.
   
    \subsection{Equidistribution and recurrence of horocycles}
    
    We recall classical results on the asymptotic behavior of translates of horocycles on finite volume quotients of $\mrm{SL}_2(\R)$ by its lattices.
    
   \begin{prop}[Proposition 2.2.1,~\cite{KleinbockMargulis-quasiunipotent}]
    \label{prop:horocycle equidist}
    For every compactly supported $\vp\in L^2([0,1])$ with integral $1$, $f\in C_c(\Ocal_1)$, and all $x\in \Ocal_1$,
    \begin{equation*}
        \lim_{t\r\infty} \int_0^1 f(g_tu_sx)\vp(s)\;ds = \int f\;d\mu_{\Ocal}.
    \end{equation*}
    Moreover, for a given $\vp$, the convergence is uniform as $x$ varies over compact subsets of $\Ocal_1$.
    \end{prop}

    \begin{prop}\label{prop:horocycle recurrence}
    For every $\e>0$ and compact set $\Kcal \subset \Ocal_1$, there exists a compact set $\Omega\subset\Ocal_1$ such that
    \begin{equation*}
        \left| \set{s\in[0,1]: g_tu_s x\notin \Omega }\right| \leq \e,
    \end{equation*}
    for all $t\geq 0$ and $x\in \Kcal$.
    \end{prop}
    \begin{proof}
        We provide a proof of this known result which holds more generally for quotients of $\mrm{SL}_2(\R)$ by a lattice for the reader's convenience.
       We indicate a proof using the integrability of the height function constructed by Eskin and Masur given by Lemma~\ref{lem:integrability of the height function}.
       Let $\a: \Ocal_1 \r \R_+$ be a function as in Lemma~\ref{lem:integrability of the height function}.
       Let $A$ denote the supremum over $x\in\Kcal$ of $\a(x)$ and let $B= A+b$, where $b>0$ is the constant provided by the lemma.
       Then, $A$ is finite by the semicontinuity of $\a$.
       Let $C = B/\e$ and $\Omega = \a^{-1}([0,C])$.
       Then, $\Omega$ is compact since $\a$ is proper.
       Moreover, for any $s\in [0,1]$ and $x\in\Kcal$, if $g_tu_s x\notin \Omega$, then $\a(g_t u_s x) > C$.
       Hence, the result follows by an application of Chebyshev's inequality.
        \qedhere        
    \end{proof}

 \subsection{Periodic orbits with distinct Lyapunov exponents}
 
 The fact underlying the presence of oscillations in Proposition~\ref{prop:partitions spread} is the existence of two periodic orbits for the geodesic flow over which the cocycle has sufficiently different growth rates. The following proposition describes the properties these periodic orbits need to satisfy. It is proved in Section~\ref{sec:pseudoanosov}.

 \begin{prop}\label{prop:pseudo anosovs}
     There exists $\e_0>0$, depending only on $\Ocal_1$, such that the following holds.
    For every $0<\e<\e_0$, there are $\w_a,\w_b\in \Ocal_1$ with periodic geodesic flow orbits with (not necessarily primitive) periods $\ell_a$ and $\ell_b$ respectively such that the following hold:
    \begin{enumerate}
        \item \label{item:close} $d(\w_a,\w_b)\leq \e$.
        
        \item \label{item:norm discrepancy} $\norm{\KZ{g_{\ell_b}}{\w_b}} \leq \e \norm{\KZ{g_{\ell_a}}{\w_a}}$.
        
        \item \label{item:close periods} $|\ell_a -\ell_b| <1$.
        
        \item \label{item:nesting} The lifts of $\w_a$ and $\w_b$ to our fundamental domain $\Dcal_0$ (cf. Section~\ref{sec:octagon}) are each at a distance at least $\e_0$ from $\partial \Dcal_0$.
        
        \item \label{item:injectivity} The injectivity radii at $\w_a$ and $\w_b$ is at least $\e_0$.

    \end{enumerate}
    \end{prop}
    
    We remark that Proposition~\ref{prop:pseudo anosovs} is essentially the only place in our arguments where we use the fact that we are working over the octagon locus.

 \subsection{Non-atomic boundary measures}
 \label{sec:outline of oscillations}
 Recall the notational convention from the beginning of the section.
 To prove Proposition \ref{prop:partitions spread}, we need to prove that certain asymptotic flags vary along stable and unstable horocycles. Briefly, we will ``match" points $\omega',\omega''\in \mathcal{C}$, so that $\KZ{g_t}{\omega'}=A\cdot \KZ{g_{\ell_a}}{\omega_a}\cdot B$ and $\KZ{g_t}{\omega''}=A\cdot \KZ{g_{\ell_b}}{\omega_b}\cdot B$. 
 Here, $A$ and $B$ will be elements of $ \SL$ occuring as common values of the cocycle along segments of the orbits of $\w'$ and $\w''$.
 Naively, for a given vector $v$ one would suspect that a discrepancy in operator norms of the form
$$ \|\KZ{g_{\ell_a}}{\omega_a}\|\ll\|\KZ{g_{\ell_b}}{\omega_b}\|$$ would imply that
$$\norm{\KZ{g_t}{\omega'}v}\ll \norm{\KZ{g_t}{\omega''}v}.$$ 
To make this work we need to show there are not coincidences between any of 
\begin{itemize}
\item $v$ and the most contracted input direction of $B$,
\item the most expanded output direction of $B$ and the most contracted input direction of $\KZ{g_{\ell_b}}{\omega_b}$
\item the most expanded output direction of $\KZ{g_{\ell_b}}{\omega_b}$ and the most contracted input direction of $A$.
\end{itemize}

    The result that we use to rule out such coincidences is the following proposition, proved in Section~\ref{sec:flags}.
    Let $K=\mrm{SO}_2(\R)$ and $A^+$ be the diagonal subsemigroup of $\SL$ with the larger eigenvalue in the top left corner. 
    For $A = \ell a k\in\SL$, with $k,\ell\in K$ and $a\in A^+$, 
    define maps $\xi_{\mrm{in}}, \xi_{\mrm{out}}:\SL\to \RP 
    $ by setting  \begin{equation}\label{eq:contracting input def}
    \xi_{\mrm{in}}(A) = k^{-1}\cdot e_2, \qquad
    \xi_{\mrm{out}}(A)=\ell\cdot k^{-1}\cdot e_1,
    \end{equation}
    where $e_1,e_2$ are the standard basis vectors of $\R^2$.
    In particular, $\xi_{\mrm{in}}(A)$ is the most contracted singular input direction and $\xi_{\mrm{out}}(A)$ be the most expanded output singular direction.

Recall the fundamental domain $\Dcal_0$ of $\Ocal$ given in Section~\ref{sec:octagon}.

\begin{prop}\label{prop:input nonatomic} 
For all $\e>0$, $v\in \RP$, and compact sets $\Kcal\subset \tilde{\Ocal}$, such that $\Kcal$ is contained in the interior of the fundamental domain $\Dcal_0$ of $\Ocal$, there exists $\d>0$ so that for any interval $I\subseteq [0,1]$, we can find $t_0>0$ such that for all $\omega\in \Kcal$ and $t\geq t_0$,
\[
|\{s \in I: \measuredangle\big(\xi_{\mrm{in}}(\mrm{KZ}(g_t,u_s\omega)),v\big)<\d\}|<\e |I|.
\]
\end{prop}
In this result, we restrict to compact sets contained in the interior of the fundamental domain to avoid technical issues arising from discontinuities of the cocycle at the boundary.


\section{Choice of pseudo-Anosovs}
\label{sec:pseudoanosov}

   The goal of this section is to prove Proposition~\ref{prop:pseudo anosovs}. 
   This section, together with Proposition~\ref{prop:strong irreducibility}, are the only places where we use specific properties of the octagon locus.
   We give a more algebraic description of the splitting of the tangent space over the octagon locus into the tautological and balanced subspaces.
   In particular, we show the splitting is defined (at the regular octagon) over a quadratic field and the two subspaces are Galois conjugates of one another.
   Using the Galois conjugate of the canonical basis of the tautological space (cf.~Section~\ref{sec:background}), we find explicit hyperbolic matrices in the Veech group of the octagon giving rise to the two periodic orbits satisfying Proposition~\ref{prop:pseudo anosovs}.

 \subsection{Monodromy over the octagon locus}   
In this subsection we analyze the octagon locus $\mathcal{O}$ and its monodromy.
Recall the notation introduced in Section~\ref{sec:octagon}. In particular,  recall that we chose a translation surface $M_0$ corresponding to the octagon. In this subsection, we interpret the Veech group $\Gamma$ as the group of affine automorphisms of $M_0$. If $\alpha\in\Gamma$ then we denote its derivative by $D\alpha$. The homomorphism that takes $\alpha$ to $D\alpha$ is the Veech homomorphism. We denote by $\a_\ast$ the induced action on the homology of $M_0$ and by $\a^\ast$ the action on cohomology groups $H^1(M_0,\Sigma;\R^2)$ or $H^1(M_0,\Sigma;\R)$. We use the term {\em monodromy} to refer to the right action of the Veech group on cohomology.

We begin by recalling the following standard result; cf.~\cite[ch.~3.1]{Hatcher}. 
\begin{lem}[Universal coefficient theorem]
\label{lem: universal coefficients}The natural homomorphism sending cohomology classes in $H^1(M_0,\Sigma;\R)$ to $\Hom_\Q(H_1(M_0,\Sigma;\Q),\R)$
is an isomorphism. Let $\alpha$ be an automorphism of $M_0$. The action of $\alpha^*$ on $H^1(M_0,\Sigma;\R)$ corresponds to the action of $\alpha_\ast$ by precomposition on $\Hom_\Q(H_1(M_0,\Sigma;\Q),\R)$.
\end{lem}

\begin{figure}
\includegraphics[scale=0.4]{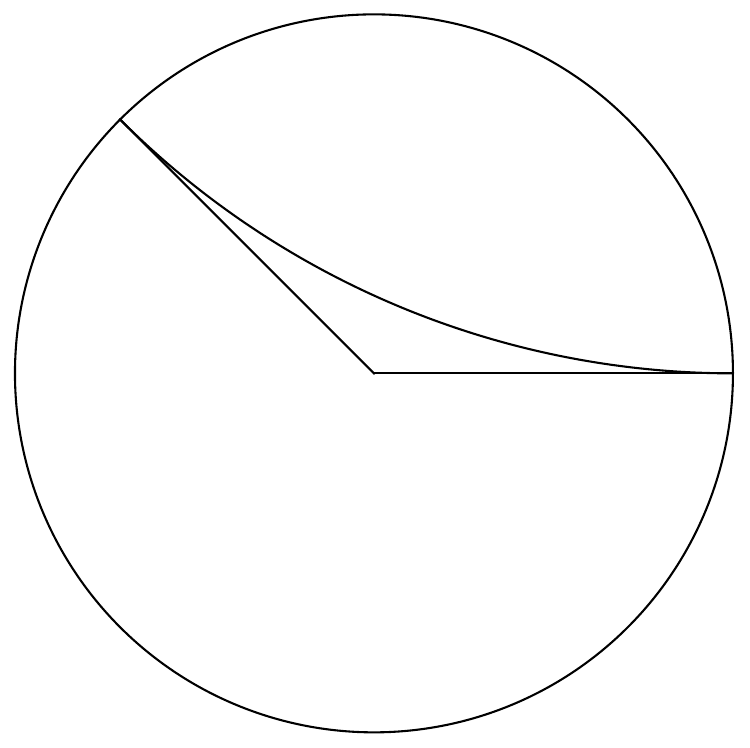}
\caption{The conjugate action of the extended Veech group is generated by reflections in the side of this triangle.}
\label{fig: conj}
\end{figure}

\begin{lem} \label{lem: Veech homomorphism} Let $\hol_0: H_1(M_0,\Sigma;\Q)\to \R^2$ be the holonomy map of $M_0$. This map is equivariant with respect to the action of the Veech group in that $\hol_0 \circ \alpha_\ast=D\alpha\circ \hol_0$ where $\alpha\in\Gamma$.
\end{lem}

\begin{proof} Let $\alpha$ be an element of $\Gamma$ which we think of as represented by an affine automorphism of $M_0$ with derivative $D\alpha$. Let $\sigma$ be an oriented saddle connection. Let $[\sigma]$ be the homology class of $\sigma$ in $H_1(M_0,\Sigma;\Z)$. The holonomy of the image of $\sigma$ under $\a$ is $D\alpha(\hol_0([\sigma]))$ thus we have $\hol_0\circ\alpha_\ast[\sigma]=D\alpha\circ \hol_0([\sigma])$. Since $H_1(M_0,\Sigma;\Z)$ is generated by saddle connections we have $\hol_0\circ \alpha_\ast=D\alpha\circ \hol_0$.
\end{proof}

The lemma tells us that the following square is commutative. 
\[\begin{tikzcd}
    H_1(M,\Sigma;\Q)\arrow{d}{\hol} \arrow{r}{\alpha_\ast}  & H_1(M,\Sigma;\Q) \arrow{d}{\hol} \\
     \R^2  \arrow{r}{D\alpha} & \R^2
  \end{tikzcd}
\]

We now recall standard facts regarding Galois conjugate representations of Veech groups; cf.~\cite{GutkinJudge}.

\begin{lem}  \label{lem: k^2 image} Let $k=\Q(\sqrt{2})\subset\R$. The holonomy map $\hol_0: H_1(M_0,\Sigma;\Q)\to\R^2$ is injective and its image is $k^2\subset\R^2$. \end{lem}

\begin{proof} Since $M_0$ has genus two $H_1(M_0,\Sigma;\Q)$ is a 4-dimensional vector space over $\Q$. Let us write 
\[H_1=H_1(S,\Sigma;\Q).\] 
By considering horizontal and vertical saddle connections contained in the regular octagon we see that $M_0$ has saddle connections with holonomy $(1,0)$, $(1+\sqrt{2},0)$, $(0,1)$ and  $(0,1+\sqrt{2})$. These generate a 4 dimensional $\Q$-subspace of $k^2$ hence all of $k^2$. Since $H_1$ has $\Q$ dimension 4 the image of $H_1$ is equal to $k^2$. 
Since the holonomy map has a 4-dimensional image its kernel is $0$ so the map is injective.
\end{proof}

According to the lemma, we can view the holonomy as taking values in either $\R^2$ or in $k^2$. We write $\mrm{hol}_k$ for the map $\mrm{hol}_k:H_1\to k^2$. Similarly we write $D\alpha_k$ for the map $D\alpha_k:k^2\to k^2$.

Lemmas \ref{lem: Veech homomorphism} and \ref{lem: k^2 image} give us the commutativity of the following diagram.
\[\begin{tikzcd}
    H_1\arrow{d}{\hol_k} \arrow{r}{\alpha_\ast}  & H_1 \arrow{d}{\hol_k} \\
     k^2  \arrow{r}{D\alpha_k} & k^2
  \end{tikzcd}
\]

Let us write $\phi_1$ for the inclusion of $k$ into $\R$. This map is a field embedding. There is a second real valued field embedding which we denote by $\phi_2:k\to\R$. If $\mrm{gal}$ denotes the Galois automorphism of $k$ which takes $\sqrt{2}$ to $-\sqrt{2}$ then $\phi_2=\phi_1\circ \mrm{gal}$.

 Lemma \ref{lem: k^2 image} implies that the image of the holonomy map, which is a priori a $\Q$ vector space, in fact has a $k$ vector space structure. If $a\in k$ and $v\in H_1$, we define $av\in H_1$ to  be $\hol_0^{-1}(a\cdot \hol_0(v))$. We use the injectivity of the holonomy to invert $\hol_0$.
  Let us write $\Hom_j(H_1,\R)$ for the $\R$-vector space of $\phi_j$-linear maps from $H_1$ to $\R$. These are $\Q$-linear maps $f$ with the additional property that $f(\lambda v)=\phi_j(\lambda) v$ for $\lambda\in k$.
  
  Let $\pi_1$ and $\pi_2$ be the coordinate projections on $\R^2$ and denote $dx_0=\pi_1\circ \hol_0$ and $dy_0=\pi_2\circ \hol_0$.
    It follows from Lemma \ref{lem: k^2 image} that $dx_0$ and $dy_0$ viewed as $\Q$-linear maps from $H_1$ to $\R$ take values in $k$. When we want to emphasize that we are dealing with $k$ valued functions we write $dx_k$ and $dy_k$ for the corresponding map from $H_1$ to $k$.
    In terms of this notation we have $dx_k=\pi_1\circ \hol_k$ and $dy_k=\pi_2\circ \hol_k$. So $dx_0=\phi_1\circ dx_k$ and $dy_0=\phi_1\circ dy_k$. So the following diagram commutes as does the analogous diagram where we replace $\pi_1$ by $\pi_2$.
    \[\begin{tikzcd}
    H_1 \arrow{r}{\hol_k}\arrow{d}{id}  & k^2 \arrow{r}{\pi_1}\arrow{d}{\phi_1\times\phi_1} &k \arrow{d}{\phi_1}\\
    H_1  \arrow{r}{\hol_0} & \R^2 \arrow{r}{\pi_1}& \R
  \end{tikzcd}
\]

It will be useful to introduce notation for the real valued ``Galois conjugate" versions of $dx_0$ and $dy_0$; namely $\overline{dx_0}=\phi_2\circ dx_k$ and $\overline{dy_0}=\phi_2\circ dx_k$.
The proof of the following result follows the same lines as~\cite[Section 2.4]{AvilaDelecroix}.

\begin{prop}\label{prop:Galois}  We have the following:
\begin{enumerate}
\item The monodromy preserves the splitting of $H^1(M_0,\Sigma;\R^2)$ into horizontal and vertical subspaces.
\item The 4-dimensional space $H^1(M_0,\Sigma;\R)$ splits as a direct sum of two 2-dimensional subspaces which we identify with $\Hom_1(H_1,\R)$ and $\Hom_2(H_1,\R)$.  
\item This splitting is  invariant under the $\Gamma$ monodromy.
\item\label{item:standard matrices} The action of  $\Gamma$ on $\Hom_1(H_1,\R)$ with respect to the basis given by $dx_0$ and $dy_0$ is given by the Veech homomorphism matrices. These matrices have coefficients in $k$.
\item\label{item:Galois on balanced} The action of $\Gamma$ with respect to the basis given by $\overline{dx_0}$ and $\overline{dy_0}$ on $\Hom_2(H_1,\R)$ is given by Galois conjugates of the matrices in Item~\eqref{item:standard matrices}.
\item $\Hom_1(H_1,\R)$ corresponds to the tautological subspace of $H^1(M_0,\Sigma;\R)$.
\end{enumerate}
\end{prop}

  In assertion (6) of the previous proposition we idenitfy $\Hom_1(H_1,\R)$ with the tautological subspace. The point of the following lemma is to identify the other summand, $\Hom_2(H_1,\R)$, with the balanced subspace.
 
\begin{lem}\label{lem:Galois=balanced}
The real vector spaces $\Hom_1(H_1,\R)$ and $\Hom_2(H_1,\R)$ are symplectically perpendicular (when viewed as subspaces of real cohomology). Thus we can identify $\Hom_2(H_1,\R)$ with the balanced subspace which is defined to be the symplectic complement of the $\GLd$-space.
\end{lem}
 
 \begin{proof}
 In order to prove this we give an alternate construction of the $k$ vector space structure on $H_1$.
Let $\alpha$ be the affine automorphism of $M_0$ corresponding to rotating the octagon by an angle of $\pi/4$ counterclockwise. Let $A$ be the corresponding rotation of $\R^2$. 
The trace of $A$ is $\sqrt{2}$ and, in particular, belongs to $k$. 
We define a linear map $L$ from $H_1(M_0,\Sigma;\Q)$ to itself  by $L(v)=(\alpha_\ast+\alpha^{-1}_\ast)(v)$ for $v\in H_1(M_0,\Sigma;\Q)$.

Let $v\in H_1(S,\Sigma;\Q)$. We calculate: 
\[\hol(L(v))= \hol \circ (\alpha_\ast+\alpha^{-1}_\ast)(v)=(A+A^{-1})(v)=\mrm{tr}(A)\cdot v=\sqrt{2}\cdot v
\]
 This uses the observation that $A$ satisfies its characteristic polynomial so $A^2-\mrm{tr}(A)\cdot A+I=0$ hence $A-\mrm{tr}(A)\cdot I+A^{-1}=0$  and $A+A^{-1}=\mrm{tr}(A)\cdot I.$

Now consider $\theta\in \Hom_1(H_1,k)$ and $\tau\in \Hom_2(H_1,k)$. For $\gamma\in H_1$ we have $\theta((\alpha_\ast+\alpha_\ast^{-1})\gamma)=\theta(\mrm{tr}(A)\cdot \gamma)=\mrm{tr}(A)\cdot\theta(\gamma)$ by the $k$-linearity of $\theta$.  We also have $\tau((\alpha_\ast+\alpha_\ast^{-1})\gamma)=\tau(\mrm{tr}(A)\cdot \gamma)= \overline{\mrm{tr}(A)}\cdot\tau(\gamma)$ by the $k$-antilinearity of $\tau$.

 We denote the  symplectic pairing coming from the cup product by $\langle\cdot,\cdot\rangle$ and we use the fact that it is invariant under the action of $\alpha$ and $(\alpha^*)^{-1}$. 

 \begin{align*} 
\mrm{tr}(A)\cdot \langle \theta,\tau\rangle&= \langle \mrm{tr}(A)\cdot\theta,\tau\rangle
=\langle (\alpha^*+(\alpha^*)^{-1}) \theta,\tau\rangle\\
&= \langle \alpha^*\theta,\tau \rangle+\langle (\alpha^*)^{-1}\theta,\tau \rangle
= \langle \theta,(\alpha^*)^{-1}\tau \rangle+\langle \theta,\alpha^*\tau \rangle\\
&=\langle \theta,\alpha^*\tau+(\alpha^*)^{-1}\tau\rangle
=\langle \theta,\overline{\mrm{tr}(A)}\cdot\tau\rangle=\overline{\mrm{tr}(A)}\cdot \langle \theta,\tau\rangle
 \end{align*}
Hence, since $\mrm{tr}(A)=\sqrt{2}\ne-\sqrt{2}=\overline{\mrm{tr}(A)}$, we have $\langle \theta,\tau\rangle=0$. 
  \end{proof}

    \subsection{Proof of Proposition~\ref{prop:pseudo anosovs}}

    In the standard continued fraction algorithm there is a correspondence between periodic sequences and affine automorphisms of the torus. The reference \cite{SU-Coding} describes a continued fraction algorithm for the regular octagon.
    This associates to a direction in the octagon a sequence of natural numbers between $1$ and $7$.  In this case there is also a correspondence between periodic sequences and affine automorphisms of the regular octagon. 
    Let

    \[
        \gamma:=
        \begin{pmatrix}
        -1&2+2\sqrt{2}\\    0&1
        \end{pmatrix},
        \qquad
        \nu_3:=
        \begin{pmatrix}
        0&1\\
        1&0
        \end{pmatrix}, \qquad
        \nu_4:=
        \begin{pmatrix}
        0&1\\
        -1&0
        \end{pmatrix}. 
    \]
These matrices occur as derivatives of (possibly orientation reversing) affine automorphisms of the regular octagon~\cite{SU-Geodesic,SU-Coding}.
The  matrices $\nu_j$ correspond to the image of the dihedral group that acts as isometries of the regular octagon while $\gamma$ is an affine involution which corresponds to a hidden symmetry in the language of Veech~\cite{Veech-hidden}.
In particular, the elements\footnote{Squaring makes $\t_2$ orientation-preserving, which simplifies some arguments.} $\t_1 = \g\nu_3$ and $\t_2 = (\g\nu_4)^2$ belong to the image of the Veech homomorphism of $\G$.

   Note that the trace of $\t_1$ and $\t_2$ is $2+2\sqrt{2}$ and $8+4\sqrt{2}$ respectively.
   In particular, the trace of each of $\t_1$ and $\t_2$ is strictly greater than $2$ while their determinant is $1$.
   Hence, they are both hyperbolic matrices corresponding to closed geodesics in $\Ocal$.
   
   We claim that there exists $N\geq 1$ such that for all $n,m\geq 1$, the matrices $\t_2^{mN}\t_1^{nN}$ are hyperbolic.
   This claim is a special instance of the general well-known fact concerning the existence of Zariski-dense Schottky subgroups inside discrete Zariski-dense subgroups of $\mathrm{SL}_2(\R)$; cf.~\cite[Prop.~4.3]{Benoist-Schottky}.
   
   Indeed, let $\t_i^{\pm}$ be the attracting and repelling fixed points of $\t_i$ on the boundary of $\mathbb{H}^2$.
   The only fixed points of $\t_2^{-1}\t_1$ are $0$ and $\infty$, neither of which is fixed by $\t_1$ and $\t_2$.
   Hence, the sets $\set{\t_1^{\pm}}$ and $\set{\t_2^{\pm}}$ are disjoint.
   Thus, given $4$ closed disjoint complex disks, $B_1^{\pm}$ and $B_2^{\pm}$, centered at $\t_1^{\pm}$ and $\t_2^{\pm}$ respectively, we can find a large enough $N$ so that $\t_i^k$ maps $\mathbb{H}^2 \setminus B_i^{-}$ into $B_i^+$ for all $k\geq N$ and for $i=1,2$.
   An application of the ping-pong Lemma then implies that $\t_1^N$ and $\t_2^N$ generate a convex cocompact (i.e. contains no parabolic elements) Schottky subgroup $\G_0$ of $\G$.
   Indeed, this can be seen by noting that any non-trivial, cyclically reduced, word $w$ in $(\t_i^{ N})^{\pm}$ has exactly two fixed points on the boundary, contained in two disjoint closed arcs which are contracted by the first letter and the inverse of the last letter in $w$ respectively. As every word is conjugate to a cyclically reduced one, this shows that all elements of $\G_0$ are hyperbolic.
   In particular, it follows that $\t_2^{mN}\t_1^{nN}$ is hyperbolic for all $n,m\geq 1$ as claimed.
   
   By Proposition~\ref{prop:Galois}, our choice of basis implies that $\G\subset \mrm{SL}_2(\Q(\sqrt{2}))$.
   Denote by $\s:\G \r \SL$ the entrywise Galois conjugation.
   We claim that
   \begin{equation}\label{eq:ratio of norms to 0}
       \frac{\norm{\s(\t_1^n)}}{\norm{\s(\t_2^m\t_1^n)}}
       \xrightarrow{m\r\infty} 0,
   \end{equation}
   uniformly in $n$.
   Indeed, observe that $\s(\t_1) $ has trace $2-\sqrt{2}$ which is strictly less than $2$ and that $\t_1$ has an infinite order. 
   Hence, $\s(\t_1)$ is an elliptic matrix of infinite order as well. 
   Moreover, the trace of $\s(\t_2)$ is strictly greater than $2$ implying it is a hyperbolic matrix.
   Hence, we can find a unit norm eigenvector $ v\in\R^2$ of $\s(\t_2)$ with eigenvalue $\l$ so that $|\l|>1$.
   Since $\s(\t_1)$ is elliptic, we have that $\norm{\s(\t_1)^{n}v}\asymp 1$ uniformly over $n\in \Z$.
    It follows that $\norm{\s(\t_2^m \t_1^n)} \gg |\l|^m $, while $\norm{\s(\t_1^n)}$ is uniformly bounded over $n\in\Z$.
    This implies~\eqref{eq:ratio of norms to 0}.

   For $m,n\in\N$, let $q_{m,n}\in\Ocal$ be a point with periodic geodesic flow orbit of primitive period $\ell_{m,n}>0$, corresponding to $\t_2^{mN}\t_1^{nN}$.
   We will find $n_0, m_0\in\N$ and $t_1,t_2\geq 0$ so that $\w_b = g_{t_1} q_{0,1}$ and $\w_a = g_{t_2} q_{m_0,n_0}$ satisfy the proposition.
   In view of Proposition~\ref{prop:Galois}\eqref{item:Galois on balanced} and Lemma~\ref{lem:Galois=balanced}, the restriction of the KZ-cocycle to $\bal{\Ocal}$ takes values in $\s(\G)$ in the basis chosen in Proposition~\ref{prop:Galois}.
   Hence, recalling that this restriction of the cocycle is also denoted $\mrm{KZ}$, we have $\mrm{KZ}(g_{\ell_{m,n}},q_{m,n}) =\s(\t_2^{mN}\t_1^{nN})$. 
   In particular,
   Item~\eqref{item:norm discrepancy} follows by~\eqref{eq:ratio of norms to 0} if $m_0$ is chosen large enough depending on $\e$. Fix one such $m_0$. 
   
   Convex cocompactness of $\G_0$ implies that there is a compact set $\Kcal\subset \Ocal$ containing $g_t q_{m,n}$ for all $n,m\in\N$ and $t\in\R$.
   Indeed, $\Kcal$ can be chosen to be the image of the compact non-wandering set for the geodesic flow on $\mrm{T}^1\mathbb{H}^2/\G_0$ under the projection $\mrm{T}^1\mathbb{H}^2/\G_0\to \Ocal$.
   Taking $\e_0$ to be smaller than the injectivity radius at all points in $\Kcal$ yields Item~\eqref{item:injectivity}.

   We can find a compact set $\tilde{\Kcal}\subset\mrm{T}^1\HHH^2$ and lifts $\tilde{q}_{m,n}\in \tilde{\Kcal}$ of $q_{m,n}$ for all $n,m$. Let $\pi:\mrm{T}^1\HHH^2 \r \HHH^2$ be the natural projection and let $x_{m,n}=\pi(\tilde{q}_{m,n})$. 
   In our notation, $\G$ acts on the right by isometries on $\HHH^2$.
   As $n\r\infty$, $y\cdot\t_1^{nN}$ converges to $\t_1^+$ on the boundary of $\HHH^2$, uniformly as $y$ varies in the compact set $\tilde{\Kcal}\cdot \t_2^{m_0 N}$.
   Hence, along a subsequence, the geodesic segment joining $x_{m_0,n}$ to $(x_{m_0,n}\cdot \t_2^{m_0 N})\cdot\t_1^{nN}$ converges (in the Hausdorff topology) to the geodesic ray $r$ joining some point $y\in\tilde{\Kcal}$ to $\t_1^+$ as $n\to \infty$.

   Since $g_t \tilde{q}_{0,1}$ converges to $\t_1^+$ as $t\r\infty$, we can find a vector $v$ tangent to the ray $r$ and at a distance at most $\e/2$ from $g_{t_1} \tilde{q}_{0,1}$ for some $t_1>0$.
   Hence, we can find $t_2$ and $n_0$ large enough so that the distance between $g_{t_2}\tilde{q}_{m_0,n_0}$ and $g_{t_1} \tilde{q}_{0,1}$ is at most $\e/2$.
   In particular, Item~\eqref{item:close} holds for $\w_b=g_{t_1} \tilde{q}_{0,1}$ and $\w_a=g_{t_2}\tilde{q}_{m_0,n_0}$.
   Recall that $\ell_{m,n}$ denoted the primitive periods of the periodic geodesics $q_{m,n}$.
   
   By Dirichlet's theorem, we can find positive integers $p$ and $q$ such that $q>2\ell_{m_0,n_0}$ and $|p -q\ell_{0,1}/\ell_{m_0,n_0}|\leq 1/q$.
   Taking $\ell_a = p\ell_{m_0,n_0}$ and $\ell_b=q\ell_{0,1}$, we obtain Item~\eqref{item:close periods}.
   
   For Item~\eqref{item:nesting}, first note that the periodic geodesic through $q_{0,1}$ is independent of $\e$.
   Moreover, the geodesics containing the boundary of our chosen fundamental domain $\Dcal_0$ of $\G$ on $\HHH^2$ connect parabolic points on the boundary at infinity; cf. Section~\ref{sec:octagon}.
   In particular, no side of $\Dcal_0$ is contained in a lift of this geodesic to $\HHH^2$. It follows that the distance of the lift of $\w_b$ to $\Dcal_0$ to the boundary $\partial\Dcal$ is uniformly bounded below by a constant $\d>0$ depending only on $\t_1$.
   Hence, Item~\eqref{item:nesting} holds whenever $\e_0<\d/2$.

    
\section{Non-atomic boundary measures}
    \label{sec:flags}

    The goal of this section is to prove Proposition~\ref{prop:input nonatomic}.
     We use the fact that the Oseledets subspaces vary along horocycles to show that the ``bad coincidences" discussed above Proposition~\ref{prop:input nonatomic} occur on a set of negligible measure.
    This issue appears frequently and the random walk approach developed in~\cite[Appendix C]{EM} is convenient for our purposes. We recall their set up and then use their results to obtain our desired conclusion in our closely related setting of the geodesic flow.
    In Proposition~\ref{prop:strong irreducibility}, we verify the strong irreducibility hypothesis of the subbundles we study.

\subsection{Boundary measures and strong irreducibility}
Let $\nu$ be a compactly supported
and $\mrm{SO}_2(\mathbb{R})$-bi-invariant probability measure on $\SL$. 
Let
\begin{align*}
    X=\SL^\Z \times \Ocal_1 ,\qquad
    \bar{\nu} = \nu^\Z.
\end{align*}
Define $T:X\to X$ by $T({\bf g},\omega)=(S{\bf g},g_0\omega)$ where $S$ is the left shift, $\mbf{g}=(\dots, g_{-1},g_0,g_1,\dots) \in \SL^\Z, $ and $\w\in\Ocal_1$.

\begin{lem}\label{lem:ergodic skew product}
    The measure $\bar{\nu}\times \mu_\Ocal$ is $T$-ergodic.
\end{lem}
\begin{proof}
By~\cite[Prop.~2.14]{BenoistQuintBook}, the lemma follows once ergodicity of $\mu_\Ocal$ as a $\nu$-stationary measure is established.
    To this end, denote by $\mrm{P}_\nu$ the averaging operator on $\mrm{L}^2(\mu_\Ocal)$ associated to $\nu$.
    That is $\mrm{P}_\nu(f)(x):=\int f(gx)\;d\nu(g)$. 
    Let $f$ be a $\mrm{P}_\nu$-invariant function. 
    We show that $f$ is constant by showing that
    \begin{align}\label{eq:orth. to non-constant}
    \langle f,\vp\rangle =0,    
    \end{align}
    for every $\vp\in \mrm{L}^2(\mu_\Ocal)$ with $\int \vp\;d\mu_\Ocal =0$.
    Fix one such $\vp$.
    Since the support of $\nu$ generates $\SL$, it follows from Oseledets' theorem and Furstenberg's positivity of the top Lyapunov exponent of the random walk on $\R^2$ generated by $\nu$ (cf.~\cite[Cor.~4.32]{BenoistQuintBook}) that $\norm{g_n\cdots g_0}$ tends to infinity for $\bar{\nu}$-almost every $\mbf{g}$.
    Hence, the Howe-Moore theorem implies that the matrix coefficient $\langle f\circ g_n\cdots g_0,\vp\rangle$ tends to $0$ almost surely.
    The dominated convergence theorem then gives that $\langle \mrm{P}_\nu^n (f),\vp\rangle \to 0$ as $n\to \infty$.
    This implies~\eqref{eq:orth. to non-constant} since $\mrm{P}_\nu(f)=f$.
\end{proof}

Recall the definition of the map $\xi_{\mrm{in}}$ in~\eqref{eq:contracting input def}.
\begin{lem}\label{lem:flag existence}
    For $\bar{\nu}\times\mu_\Ocal$-almost every $(\mbf{g},\w)\in X$, the limit as $n\to\infty$ of $\xi_{\mrm{in}}(\mrm{KZ}(g_n....g_0,\omega))$ exists in $\RP$. 
\end{lem}
\begin{proof}
    This is a consequence of Oseledets' theorem and positivity of the top Lyapunov exponent as we now explain.

    The restriction of the $\mrm{KZ}$-cocycle to the ($2$-dimensional) balanced space induces a cocycle denoted $\a$ over the dynamical system $(X,T, \bar{\nu}\times\mu_\Ocal)$ defined by $\a(1,(\mbf{g},\w)) = \KZ{g_0}{\w}$.
    It follows by~\cite[Lemma 2.1']{Forni} that $\norm{\KZ{g}{\w}}\ll \norm{g}$ uniformly for all $(g,\w)\in X$.
    In particular, $\log$-integrability of the cocycle $\a$ follows since $\nu$ is compactly supported.
    Combined with Lemma~\ref{lem:ergodic skew product}, this means that Oseledets' theorem provides a well-defined top Lyapunov exponent for $\a$.

    Moreover, it follows by~\cite[Theorem 15.1]{Bainbridge-LyapunovGenus2} that this top exponent is positive.
    Indeed, the quoted result concerns the exponent with respect to $\mu_\Ocal$ of the cocycle over the geodesic flow on $\Ocal_1$.
    That this implies positivity of the exponent of $\a$ follows from the fact that random products of the form $g_n\cdots g_0$ almost surely shadow a geodesic, up to sublinear error; cf.~\cite[Lemma 4.1]{ChaikaEskin} and~\cite[Remark 7.4]{AAEKMU} for more details on the relationship between these two cocycles.
    
    Since the balanced space is two dimensional and the cocycle is $\SL$-valued, Oseledets' theorem then implies that the second (i.e.~bottom) Lyapunov space is a well-defined, one-dimensional space occuring as the limit of $\xi_{\mrm{in}}(\mrm{KZ}(g_n....g_0,\omega))$ almost surely.
    This concludes the proof.
    \end{proof}

By slight abuse of notation, we set
\begin{align*}
    \xi_{\mrm{in}}({\bf g},\omega)=\lim_{n\to\infty} \xi_{\mrm{in}}(\mrm{KZ}(g_n....g_0,\omega)),
\end{align*}
on the full measure set where the limit exists. Note that $\xi_{\mrm{in}}(\mbf{g},\w)$ depends only on the ``future" of $\mbf{g}$, i.e.~it depends only on the non-negative coordinates of $\mbf{g}$.

Let $\hat{\nu}$ be the measure on $\mathcal{O}_1\times \RP$ defined by 
$$\hat{\nu}(\mathcal{S}\times U)=\int_{\mathcal{O}_1}\chi_{\mc{S}}(\omega)\bar{\nu}(\{{\bf g} :\xi_{\mrm{in}}({\bf g},\omega)\in U\})d\mu_{\Ocal}(\w).$$ 

Let $\SL$ act on $\mathcal{O}_1 \times \RP$ by $g(\omega,v)=(g\omega, \mrm{KZ}(g,\omega)v)$, where the $\mrm{KZ}$ cocycle acts projectively on the second coordinate.
The following lemma follows from the equivariance identity: $\xi_{\mrm{in}}({\bf g},\omega)=\mrm{KZ}(g_0,\omega)^{-1}\xi_{\mrm{in}}(\mrm{KZ}(S{\bf g},g_0\omega))$.

\begin{lem} $\hat{\nu}$ is a $\nu$-stationary measure. That is 
\[\hat{\nu}(\mathcal{S} \times U)=\int_{\SL}\hat{\nu}(g(\mathcal{S}\times U))d\nu(g).\]  
\end{lem}

A measurable almost invariant splitting is a finite collection of measurable maps $V_i:\Ocal \to \RP$ such that
\begin{enumerate}
    \item $V_i(x)\neq V_j(x)$ for all $i\neq j$ and $\mu_{\Ocal}$-almost every $x$.
    \item For $\nu\times\mu_\Ocal$-almost every $(g,x)\in \SL\times \Ocal_1$ and every $i$, there is $j$ such that $\KZ{g}{x}V_i(x)=V_j(gx)$.
\end{enumerate}
We say the cocycle is \textbf{strongly irreducible} (relative to $\nu$ and $\mu_{\Ocal}$) if it does not admit any measurable almost invariant splitting. Strong irreducibility implies the following result.

\begin{lem}[Lemma C. 10,~\cite{EM}]
\label{lem:EMM C10} If $\mrm{KZ}$ is strongly irreducible on $\mathbb{R}^2$, then for almost every $\omega \in \mathcal{O}_1$ and every $v\in \RP$ we have 
$\bar{\nu}(\{{\bf g}:\xi_{\mrm{in}}({\bf g},\omega)=v\})=0$.

\end{lem} 

We remark that the above definition of strong irreducibility differs slightly from the one given in~\cite{EM}, however the above definition is the one used in their proof.

In order to apply the above results, we verify the strong irreducibility of the $\mrm{KZ}$ cocycle in our setting using the results of Filip~\cite{Filip-semisimplicity} and Smillie-Ulcigrai~\cite{SU-Geodesic,SU-Coding}.

    \begin{prop}\label{prop:strong irreducibility}
        The restriction of the $\mathrm{KZ}$ cocycle over the octagon locus $\mc{O}$ to the subbundle with fibers the balanced subspace is strongly irreducible (relative to $\nu$ and $\mu_{\mc{O}}$).
    \end{prop}

    \begin{proof}
    
    Suppose the restriction of the cocycle is not strongly irreducible.
    Then, we can find a measurable almost invariant splitting $L_1,\dots,L_k:X\to\RP$.
    Let $Y = (\RP)^k/S_k$, where $S_k$ is the symmetric group acting by permutations of the coordinates.
    Let $H=\SL$.
    Note that the diagonal action of $H$ on $Y$ is smooth in the sense of~\cite[Definition 2.1.9]{Zimmer1984}, i.e. $Y/H$ admits a countable collection of Borel sets which separate points.
    Indeed, by a result of Borel-Serre (cf.~\cite[Theorem 3.1.3]{Zimmer1984}), since $H$ acts algebraically on the variety $(\RP)^k$, $H$-orbits are locally closed. This remains true on $Y$ since $S_k$ is a finite group.
    
    Recall that $X=\SL^\Z\times \Ocal_1$.
    As in Lemma~\ref{lem:flag existence}, we define a cocycle $\mrm{\a}: \Z\times X \to H$ over the skew-product on $X$ given by $\a(1,({\bf g},x)) = \KZ{ g_0}{ x}$. 
    Let $\vp:X\to Y$ be given by $\vp({\bf g},x)=[L_1( x),\dots,L_k(x)]$, for $\bar{\nu}\times\mu_\Ocal$-almost every $({\bf g},x)$.
    Since the $L_i$'s are an almost invariant splitting, $\vp$ is invariant by the cocycle $\a$, i.e. 
    \[\vp(T({\bf g},x)) = \vp(S{\bf g}, g_0 x) = [L_1( g_0x)\dots, L_k( g_0x)] = \a(1,({\bf g},x)) \vp({\bf g},x). \]
    In particular, $\vp$ is invariant in the sense of~\cite[Definition 4.2.17]{Zimmer1984} (where in our notation the acting group is $\Z$ by the skew-product on $X$).

    Recall that $\bar{\nu}\times \mu_\Ocal$ is $T$-ergodic by Lemma~\ref{lem:ergodic skew product}.
    Hence, we can apply Zimmer's cocycle reduction lemma~\cite[Lemma 5.2.11]{Zimmer1984} (since the cocycle is $H$-valued), to get that $\a$ is equivalent to a cocycle taking values in the stabilizer in $H$ of a point $y\in Y$, denoted $H_y$.
    This means that, up to a measurable change of basis, we may assume that $\KZ{g}{x}\in H_y$ for $\nu$-a.e.~$g$ and for $\mu_\Ocal$-a.e.~$x$.

    We wish to apply the results of Filip~\cite[Section 8.1]{Filip-semisimplicity}.
    We recall his terminology.
    The above discussion implies that the measurable algebraic hull $F$ of the cocycle is, up to conjugacy, contained in the group $H_y$.
    In other words, there is a measurable $\SL$-equivariant section $\s$ of the principal $H$-bundle over $\Ocal_1$ whose fibers are isomorphic to $H/F$.
    By~\cite[Theorem 8.1]{Filip-semisimplicity}, $\s$ agrees almost everywhere with a real analytic section of this principal bundle and in particular we may assume that $\s$ is everywhere defined.

    This implies that there is a finite collection of continuous sections $\ell_1,\dots,\ell_j$ of the bundle over $\Ocal_1$ with fibers the balanced space and which are permuted by the cocycle.
    By the work of Smillie-Ulcigrai~\cite{SU-Geodesic,SU-Coding}, the Veech group contains a pseudo-Anosov element $\t_1$ whose monodromy action on the balanced space is given by an elliptic element of infinite order, cf.~the proof of Proposition~\ref{prop:pseudo anosovs}.
    Let $x\in\Ocal_1$ be a periodic point for the geodesic flow, with period $t_0$, corresponding to $\t_1$.
    Since $g_{t_0}x=x$, we get that $\KZ{g_{t_0}}{x}=\t_1$ must permute the finite collection of lines $\ell_i(x)$. This contradicts the fact that an elliptic element of infinite order acts strongly irreducibly, i.e. does not permute any finite collection of lines.
    \end{proof}

\subsection{Non-atomicity of the boundary measures for the geodesic flow}

    We now translate the above results into results about the geodesic flow.
    First, note that Oseledets' theorem implies that
    \begin{equation*}
        \xi_{\mrm{in}}(x) := \lim_{t\r\infty} \xi_{\mrm{in}}\left(\mrm{KZ}(g_t,x) \right)
    \end{equation*}
    exists in $\RP$ for $\mu_\Ocal$-almost every $x$.
    Lemma~\ref{lem:EMM C10} and Proposition~\ref{prop:strong irreducibility} yield the following corollary.
    
\begin{cor}\label{cor:nonatomic flag} For all $\e>0$, there exists $\d>0$ so that for any $v \in \RP$ we have $\mu_{\mathcal{O}}(\{\omega:\measuredangle(v,\xi_{\mrm{in}}(\omega))<\d\})<\e$. 
\end{cor}

\begin{proof}
Given $\w\in \Ocal_1$, let $\nu_\omega$ be the measure on $\RP$ given by $\nu_\omega(U)=\bar{\nu}(\{{\bf g}:\xi_{\mrm{in}}({\bf g},\omega)\in U\})$.
By the $\mrm{SO}_2(\mathbb{R})$-invariance of $\nu$, $\nu_\omega=\nu_{r_\theta\omega}$.
By sublinear tracking (cf.~\cite[Lemma 4.1]{ChaikaEskin}), there exists a measure $\sigma$ on $\mrm{SO}_2(\R)$ so that $$\sigma(\{r_\theta\in\mrm{SO}_2(\R) : \underset{t \to \infty}{\lim} \xi_{\mrm{in}}(\mrm{KZ}(g_t,r_\theta \omega))\in U\})=\nu_\w(U).$$ 
 Since $\nu_\omega=\nu_{r_\theta\omega}$, we have that $\sigma$ is  $\mrm{SO}_2(\mathbb{R})$-invariant and, hence, is the Haar measure on $\mrm{SO}_2(\R)$.
  Recall the Iwasawa decomposition $\SL=U^-AK$, where $K=\mrm{SO}_2(\R)$, $A$ is the diagonal group and $U^-$ is the lower triangular unipotent group.
  In particular, the measure $\mu_\Ocal$ is locally equivalent to the product of the Haar measures on these $3$ subgroups.
  Moreover, for any $\w\in\Ocal_1$ and for $a$ and $b$ small enough depending on the injectivity radius at $\w$, $\xi_{\mrm{in}}(g_t,\hat{u}_ag_br_\th\omega)$ depends only on $r_\th$.
 Thus, the estimate in the corollary follows by Lemma \ref{lem:EMM C10}, since the cocycle is strongly irreducible by Proposition~\ref{prop:strong irreducibility}. 
\end{proof}

    We need the following lemma before starting the proof of Proposition~\ref{prop:input nonatomic}. It relates the singular vectors of a product of two matrices to those of the matrices in the product. 
    
    \begin{lem}\label{lem:singular product}
    Let $A,B \in \SL$. Let $v= \frac{B^{-1}\xi_{\mrm{in}}(A)}{\norm{B^{-1}\xi_{\mrm{in}}(A)}}$ and let $w=\xi_{\mrm{in}}(AB)^\perp$.
    Then, 
    \begin{equation*}
        \left| \langle v,w\rangle \right| \leq \left( \frac{\norm{B}}{\norm{A}} \right)^2, \qquad
        \left|\langle \xi_{\mrm{in}}(B), w\rangle \right| \leq 
        \left( \frac{\norm{A}}{\norm{B}}\right)^2.
    \end{equation*}
    In particular, for all $A,B,C \in \SL$,
    \begin{equation*}
        \left|\sin \measuredangle (\xi_{\mrm{in}}(ABC), C^{-1}\xi_{\mrm{in}}(B)) \right|^2 \leq 
        \frac{\norm{A}^2+\norm{C}^2}{\norm{B}^2}.
    \end{equation*}
    \end{lem}
    \begin{proof}
    Recall that $\norm{M^{-1}} = \norm{M}$ for all $M\in\SL$.
    Since $\set{w,w^\perp}$ provides an orthonormal basis of $\R^2$ and $\set{AB \cdot w, AB\cdot w^\perp}$ is an orthogonal basis, it follows that
    \begin{equation}\label{eq:orthonormal estimate}
        \norm{AB}^2\left| \langle v,w\rangle \right|^2 \leq \norm{AB\cdot v}^2 = 
        \left( \frac{1}{\norm{A} \norm{B^{-1}\xi_{\mrm{in}}(A)}} \right)^2.
    \end{equation}
    To estimate the denominator, note that $\norm{B^{-1}\xi_{\mrm{in}}(A)} \geq 1/\norm{B}$. We also have that
    \begin{equation}\label{eq:product norm lower}
        \norm{AB \cdot \frac{B^{-1}\xi_{\mrm{in}}(A)^\perp}{\norm{B^{-1}}}} = \frac{\norm{A}}{\norm{B}} \leq \norm{AB}.
    \end{equation}
    The estimates in~\eqref{eq:orthonormal estimate} and~\eqref{eq:product norm lower} imply the first assertion.
    Similarly, we observe that
    \begin{equation*}
        \norm{AB}^2  \left|\langle \xi_{\mrm{in}}(B), w\rangle \right|^2
        \leq \norm{AB\cdot \xi_{\mrm{in}}(B)}^2 \leq 
        \left( \frac{\norm{A}}{\norm{B}}\right)^2.
    \end{equation*}
    Hence, for the second assertion, it suffices to note that $\norm{AB} \geq \norm{A^{-1}}^{-1} \norm{B} = \norm{B}/\norm{A}$.
    
    For the final assertion, note that for any $2$ unit vectors, $v,w$, we have $|\sin \measuredangle (v,w)| = \left| \langle v,w^\perp\rangle\right|$. Moreover, $|\sin \measuredangle (v,w)|$ gives a distance on $\RP$. In particular, it satisfies the triangle inequality. The estimate then follows from the previous two inequalities.
    \end{proof}

\subsection{Proof of Proposition \ref{prop:input nonatomic}}
    Let $\e>0$ and $\Kcal\subset \Ocal_1$ be given. 
    For concreteness, we will prove the proposition for the interval $I=[0,1]$.
    The general case follows by minor modifications.
    We note that it suffices to show that for every $\w\in \Kcal$, we can find $\d_\w, t_\w>0$ and a neighborhood $\Ucal_\w$ of $\w$ in $\Ocal_1$ so that for all $x\in \Ucal_\w$ and all $v\in \RP$,
    \begin{equation}\label{eq:prove for neighborhood}
        \bigg|\bigg\lbrace s\in [0,1]: \xi_{\mrm{in}} (\mrm{KZ}(g_t,u_s x))\in B_{\RP}\left(v,\d_\w\right) \bigg\rbrace\bigg| 
        \leq \e , \qquad
        \forall t\geq t_\w.
    \end{equation}
    Indeed, by picking a finite subcover of $\Kcal$ of $\set{\Ucal_\w:\w\in \Kcal}$, say $\set{\Ucal_{\w_1},\dots,\Ucal_{\w_n}}$, and taking $t=\max\set{t_{\w_i}:1\leq i\leq n}$ and $\d=\min\set{\d_{\w_i}:1\leq i \leq n}$, we obtain the result.

    Let $p=\hat{u}_r g_z$.
    Define $r_p(s) = r/(1+s r)$ and $t_p(s) = 2\log(1+sr)+z$. We also define $f_p(s) := e^{-2z}s/(1+sr)$.
    We then have that 
    \[u_s p = \hat{u}_{r_p(s)} g_{t_p(s)} u_{f_p(s)}. \]
    The cocycle property then implies that
    \begin{align}\label{eq:cocycle property}
        &\mrm{KZ}(g_t,u_sp\w) =\\ &\underbrace{\mrm{KZ}(\hat{u}_{e^{-t}r_p(s)}g_{t_p(s)},g_{t}u_{f_p(s)}\w)}_{A(\w,s,p)}
        \mrm{KZ}(g_{t},u_{f_p(s)}\w)
        \underbrace{\mrm{KZ}(\hat{u}_{r_p(s)}g_{t_p(s)},u_{f_p(s)}\w)^{-1}}_{C(\w,s,p)}.
        \nonumber
    \end{align}
    Using the above formula and Lemma~\ref{lem:singular product}, we will transfer the measure estimate at a point $\w\in \Kcal$ to points of the form $p\w$, for $p$ in the following set:
    \[
        \mc{W}_k = \set{p=\hat{u}_r g_z: |r|,|z|<1/k}.
    \]
    
    We first require several preliminary estimates.
    Since the cocycle is bounded on $B\times\Kcal$, for any bounded set $B\subset \SL$, it follows that the set
    \begin{equation*}
        \mc{Z} := \set{C(\w,s,p): \w\in\Kcal,s\in [0,1], p\in\mc{W}_2}
    \end{equation*}
    is finite and depends only on $\Kcal$. It follows that the constant
    \begin{equation*}
        C:= \sup \set{\norm{\g}^2:\g\in\mc{Z}}
    \end{equation*}
    is finite. In particular, elements of $\mc{Z}$ are uniformly Lipschitz on $\RP$. Hence, we can find a constant $D \geq 1$ such that for all $v\in\RP$, $\g\in \mc{Z}$, and $\d>0$,
    \begin{equation}\label{eq:proj stretch}
         B_{\RP}(\g v, D^{-1}\d) \subseteq
          \g \cdot  B_{\RP}(v,\d) \subseteq B_{\RP}(\g v, D\d).
    \end{equation}

    Endow $\mc{Z}$ with the Borel $\s$-algebra induced from its discrete topology. 
    For every $\w\in \Kcal$ and $p\in\mc{W}_2$, denote by $\Pcal_p(\w)$ the pull-back $\s$-algebra under the map $s\mapsto C(\w,s,p)$.
    Since the cocycle is locally constant, the set of atoms $P_p(\w)$ of $\Pcal_p(\w)$ consists of sub-intervals of $[0,1]$ satisfying
    \begin{equation*}
        \eta(\w) := \inf \set{|I| : I\in P_p(\w),  p\in\mc{W}_2} >0.
    \end{equation*}
    We note that $\eta(\w)$ depends on the closeness of $\w$ to the boundary of the fundamental domain defining the cocycle.
    In particular, since $\Kcal$ is contained in the interior of the fundamental domain, then  $\eta(\w)$ is uniformly bounded below over all $\w\in\Kcal$ by some constant $\eta>0$.
     
   By Proposition~\ref{prop:horocycle recurrence}, there is a compact set $\Omega \subset \Ocal_1$, depending only on $\Kcal$ and $\e$, so that for all $t\geq 0$ and $\w\in \Kcal$,
    \begin{equation}\label{eq:recurrence}
        \left|\set{s\in[0,1]: g_tu_s\w\notin\Omega} \right| < \frac{\e \eta}{100}.
    \end{equation}
    Since $\Omega$ is compact, it follows that the constant
    \begin{equation*}
        A := \sup \set{ \norm{\mrm{KZ}(\hat{u}_{e^{-t}r_p(s)}g_{t_p(s)},x)}^2: 
        t\geq 0, s\in [0,1], p\in\mc{W}_2, x\in \Omega}
    \end{equation*}
    is finite.
    
    Fix some $\w\in \Kcal$. It follows by the positivity of the top exponent of the cocycle and the Oseledets' genericity result in~\cite[Theorem 1.2]{ChaikaEskin}\footnote{Note that~\cite[Theorem 1.2]{ChaikaEskin} is phrased for $r_\theta \omega$ for almost every $\theta \in S^1$, but it is straightforward that this implies our claim.}
    that for almost every $s$, $\norm{\mrm{KZ}(g_t,u_s\w)} \r\infty$ as $t\r\infty$.
    Hence, by Egoroff's theorem, given $\d>0$, we can find $t(\w,\d)>0$ such that 
    \begin{equation}\label{eq:large norm}
        \left|\set{s\in[0,1]: \norm{\mrm{KZ}(g_t,u_s\w)} \leq \frac{100 D(A+C)}{\d}} \right| < \frac{\e\eta}{100}, \quad \forall t\geq t(\w,\d).
    \end{equation}

    Fix $\g\in \mc{Z}$ and $p\in\mc{W}_2$ and let $I\in P_p(\w)$ be such that $C(\w,\cdot,p) \equiv \g$ on the interval $I$.
    We claim that, given $\d>0$ and $t(\w,\d)$ as in~\eqref{eq:large norm}, for all $t \geq t(\w,\d)$, we have
    \begin{align}\label{eq:transfer measure}
       \bigg\lbrace s\in I: \xi_{\mrm{in}} &(\mrm{KZ}(g_t,u_s p \w))\in B_{\RP}\left(v,\frac{\d}{2D}\right), g_t u_s \w \in \Omega \bigg\rbrace \\
        &\subseteq 
         \bigg\{s\in [0,1]:  
        \xi_{\mrm{in}}(\mrm{KZ}(g_{t}, u_{f_p(s)}\w)) \in B_{\RP}(\g v, \d) \bigg\}. \nonumber
    \end{align}
    Indeed, Lemma~\ref{lem:singular product},~\eqref{eq:cocycle property}, and the choice of the constants $A$ and $C$ imply that
    \begin{equation*}
        d_{\RP}(\xi_{\mrm{in}}(\mrm{KZ}(g_t,u_sp\w)), \g^{-1}\cdot \mrm{KZ}(g_t,u_{f_p(s)}\w)) \leq \frac{A+C}{\norm{\mrm{KZ}(g_t,u_{f_p(s)}\w)}^2},
    \end{equation*}
    whenever $g_t u_s \w \in \Omega$ and $t \geq t(\w,\d)$.
    Hence, the estimate in~\eqref{eq:large norm} and the triangle inequality imply
    \begin{equation*}
        d_{\RP}(\g v, \xi_{\mrm{in}}\left(\mrm{KZ}(g_t,u_{f_p(s)}\w)\right)) \leq \d,
    \end{equation*}
    whenever $\xi_{\mrm{in}} (\mrm{KZ}(g_t,u_s p \w))\in B_{\RP}\left(v,\frac{\d}{2D}\right)$.
    This proves~\eqref{eq:transfer measure}.
    
    Let $B_\w = \set{u_\ell g\cdot \w: 0\leq\ell\leq 1, g\in \mc{W}_{100}}$. Then, $B_\w$ has non-empty interior and, in particular, $m_w:=\mu_{\Ocal_1}(B_\w) >0$.
    Applying Corollary~\ref{cor:nonatomic flag} with $\e \eta m_\w/100$ in place of $\e$, we obtain $\d_\w >0$ so that for all $v\in\RP$,

\begin{equation*}
        \mu_{\Ocal}(x\in B_\w: \xi_{\mrm{in}}(x) \in B_{\RP}(v,\d_w)) < \e\eta m_\w/200.
    \end{equation*}
    Moreover, by Oseledets' theorem and the dominated convergence theorem, we have
    \begin{align*}
        \mu_{\Ocal}(x\in B_\w: &\xi_{\mrm{in}}(x) \in B_{\RP}(v,\d_w))\\
        &= \lim_{t\r\infty}
        \mu_{\Ocal}(x\in B_\w: \xi_{\mrm{in}}(\mrm{KZ}(g_t,x)) \in B_{\RP}(v,\d_w)).
    \end{align*}
    Hence, taking $t_\w$ to be large enough, we obtain that for all $t\geq t_\w$ and $v\in\RP$,
    \begin{equation*}
        \mu_{\Ocal}\Big(\{x\in B_\w: \xi_{\mrm{in}}(\mrm{KZ}(g_t,x)) \in B_{\RP}(v,\d_w)\}\Big)
        \leq \e\eta m_\w/100.
    \end{equation*}
    It then follows by Fubini's theorem, that we can find $p_0\in\mc{W}_{100}$, so that
    \begin{equation}\label{eq:estimate at omega}
        \left| \set{s\in[0,1]: d_{\RP}(v,\xi_{\mrm{in}}(\mrm{KZ}(g_t,u_sp_0\w))) <\frac{\d_w}{2 D} } \right| < \frac{\e\eta}{100},
    \end{equation}
    for all $t\geq t_\w$.
    
    Using the above estimates, we show that
    \begin{equation}\label{eq:transfer+recurrence}
        \bigg|\bigg\lbrace s\in I: \xi_{\mrm{in}} (\mrm{KZ}(g_t,u_s p \w))\in B_{\RP}\left(v,\frac{\d_\w}{2D}\right) \bigg\rbrace\bigg| 
        \leq \e |I|,
    \end{equation}
    for all $p\in \mc{W}_{100}$ and $I\in P_p(\w)$. Let $p\in\mc{W}_{100}$ be fixed.
    We apply~\eqref{eq:transfer measure} with $\d=\d_\w$, with $p_0\w$ in place of $\w$, and $pp_0^{-1}$ in place of $p$, noting that $pp_0^{-1}\in\mc{W}_2$. 
    By enlarging $t_\w$ if necessary, we may assume it is larger than the constant $t(\w,\d_\w)$ provided in~\eqref{eq:large norm}.
    Then, the set on the right side of~\eqref{eq:transfer measure} is the preimage under $f_p$ of the set of $s$ so that $ \xi_{\mrm{in}}(\mrm{KZ}(g_{t}, u_sp_0\w)) \in B_{\RP}(\g v, \d_\w)$.
    The latter set has measure at most $\e \eta/100$ by~\eqref{eq:estimate at omega}.
    Moreover, the derivative of $f_p$ lies in the interval $[1/50, 50]$. It follows that the set on the left side of~\eqref{eq:transfer measure} has measure at most $\e\eta /2 \leq \e |I|/2$.
    Combined with the recurrence estimate in~\eqref{eq:recurrence}, we obtain~\eqref{eq:transfer+recurrence}.
    
    Since~\eqref{eq:transfer+recurrence} holds for all intervals $I$ in the partition $P_p(\w)$, then it also holds for any union of partition elements and, in particular, for the entire interval $I=[0,1]$.
    This implies that for all $x=g \w$, where $g= u_{\t}p$ with $|\t|<\e$ and $p\in\mc{W}_{100}$, we have
    \begin{equation*}
        \bigg|\bigg\lbrace s\in [0,1]: \xi_{\mrm{in}} (\mrm{KZ}(g_t,u_s x))\in B_{\RP}\left(v,\frac{\d_\w}{2D}\right) \bigg\rbrace\bigg| 
        \leq 2\e, \qquad\forall t\geq t_\w.
    \end{equation*}
    Running the above argument with $\e/2$ in place of $\e$ and 
    taking $\Ucal_\w = \set{u_\t p:|\t| <\e, p\in \mc{W}_{100}}$, we obtain~\eqref{eq:prove for neighborhood}.

    
 \section{Proof of Proposition~\ref{prop:partitions spread}}
    \label{sec:proof of oscillations}

    \subsection{Strategy of the proof}
    Fix some $v_0 \in \R^2$ with $\norm{v_0}=1$.
    Denote by $C(r,t)$ the set of points along the horocycle for which the norm of the cocycle at time $t$ is concentrated in the interval $[r,r/\k]$. More precisely, we let
    \begin{equation*}
        C(r,t) := \set{s\in [0,1]: \norm{\KZ{g_t}{u_s\w_1}v_0} \in \left[\k r, \frac{r}{\k}\right] }.
    \end{equation*}
    Suppose $s_1,s_2\in C(r,t)$. Then, 
    \begin{equation}\label{eq:tight}
       \k^2\leq \frac{\norm{\KZ{g_t}{u_{s_1}\w_1}v_0} }{\norm{\KZ{g_t}{u_{s_2}\w_1}v_0}}
       \leq \k^{-2}.
    \end{equation}
    Based on this observation, we use the following matching procedure
    to bound the measure of $C(r,t)$ whenever $t$ is sufficiently large, depending only on $\k$.
    We find subsets $\Pcal_a(t,v_0)$ and $\Pcal_b(t,v_0)$ of $[0,1]$ and a matching map $\mathbb{M}:\Pcal_a(t,v_0) \r \Pcal_b(t,v_0)$ such that
    \begin{enumerate}
        \item \label{item:strategy1}
        $\Pcal_a(t,v_0) \cap \Pcal_b(t,v_0) =\emptyset$ and the measures of each of $\Pcal_a(t,v_0)$ and $\Pcal_b(t,v_0)$ are at least $ 1/50$.
        
        \item \label{item:strategy2} 
        The Jacobian of $\mathbb{M}$ with respect to the Lebesgue measure is within $\kappa$ from $1$. More precisely, given a Borel set $A\subseteq \Pcal_a(t,v_0)$, we have
        \begin{equation*}
            (1+\k)^{-1} |A| \leq  |\mathbb{M}(A)| \leq (1+\k) |A|.
        \end{equation*}
        
        \item \label{item:strategy3}
        For all $s\in \Pcal_a(t,v_0)$, 
        \begin{equation*}
            \norm{\KZ{g_t}{u_s\w_1}v_0} < \k^2 \norm{\KZ{g_t}{u_{\mathbb{M}(s)}\w_1}v_0 }.
        \end{equation*}
        In particular, $\mathbb{M}(C(r,t)\cap \Pcal_a(t,v_0)))$ is disjoint from $C(r,t)\cap \Pcal_b(t,v_0)$.
        
    \end{enumerate}

        The existence of sets $\Pcal_a(t,v_0)$ and $\Pcal_b(t,v_0)$ with the above properties immediately implies that the measure of $|C(r,t)|$ is at most $49/50+\k$. Indeed, 
        \begin{align*}
            |C(r,t)| &\leq |[0,1]\setminus \left(\Pcal_a(t,v_0)\sqcup \Pcal_b(t,v_0) \right)| \\
            &\qquad \qquad  \qquad
            + |C(r,t)\cap \Pcal_a(t,v_0)| + |C(r,t)\cap \Pcal_b(t,v_0)| 
            \\ &\leq 1- |\Pcal_a(t,v_0)|-|\Pcal_b(t,v_0)| \\
            & \qquad \qquad  \qquad
            + \left(1+\k\right) |\mathbb{M}(C(r,t)\cap \Pcal_a(t,v_0))| + |C(r,t)\cap \Pcal_b(t,v_0)| \\ 
            &\leq  1- |\Pcal_a(t,v_0)|+  \k \leq 49/50+\k ,
        \end{align*}
        where we used the fact that the set $\mathbb{M}(C(r,t)\cap \Pcal_a(t,v_0)))$ is disjoint from $C(r,t)\cap \Pcal_b(t,v_0)$.

    Briefly, the essential property of the sets $\Pcal_\bullet(t,v_0)$, $\bullet\in\{a,b\}$, which allows us to verify Item~\eqref{item:strategy3} above is as follows.
    We produce points $\w_a$ and $\w_b$ with periodic geodesic flow orbits of periods $\ell_a$ and $\ell_b$ respectively and which satisfy Proposition~\ref{prop:pseudo anosovs}.
    The sets $\Pcal_a$ and $\Pcal_b$ are constructed so that an orbit segment $(g_r u_s\w_1)_{0\leq r\leq t}$ for some $s\in \Pcal_a(t,v_0)$ fellow travels with the orbit $(g_r u_{\mathbb{M}(s)}\w_1)_{0\leq r\leq t}$ for all times in $[0,t]$, except for a window of time of the form $[\t,\t+\ell_a]\subset [0,t]$, where the first orbit tracks the periodic geodesic of $\w_a$ while the second tracks $\w_b$. Fellow travelling ensures the value of the cocycle along the two orbit segments of $u_s\w_1$ and $u_{\mathbb{M}(s)}\w_1$ essentially only differs during the shadowing window $[\t,\t+\ell_a]$. The discrepancy in the cocycle norm along those two orbits will then follow from Proposition~\ref{prop:pseudo anosovs}~\eqref{item:norm discrepancy}.
    
    One of the difficulties with implementing this strategy comes from the fact that the norm of the cocycle for the window $[0,t]$ is in general not proportional to the product of the norms of the cocycle before, during, and after fellow travelling.
    However, this proportionality holds when the singular vectors of the relevant matrices are in general position with respect to one another.
    Using results of Section~\ref{sec:flags}, we show that this general position property holds outside an exceptional set of small measure along the horocycle.
    We now carry out the details.
    
    \subsection{Setup}
    We let $\e,\e_2>0$ be two parameters which we will make small over the course of the proof to satisfy certain requirements. 
    Towards the end of the proof, we will make $\e$ small depending on $\e_2$; cf.~\eqref{eq:matching discrepancy} and the discussion following it. 
    Hence, for clarity, we will highlight the quantities depending on these two constants.
    
    Let $\w_a,\w_b \in \Ocal$ be two points with periodic geodesic flow orbits satisfying Proposition~\ref{prop:pseudo anosovs} with $\e/500$ in place of $\e$.
    For $\bullet=a,b$, define the following Bowen balls:
    \begin{equation*}
        B(\bullet, \ell_\bullet, 10^{-3}\e_2) = 
        \set{x\in\Ocal: d(g_r \w_\bullet, g_rx) \leq 10^{-3}\e_2, \text{ } \forall 0\leq r\leq \ell_\bullet }.
    \end{equation*}

  We may assume that $\e_2$ is smaller than the constant $\e_0$ in Proposition~\ref{prop:pseudo anosovs} so that the injectivity radii at $\w_a$ and $\w_b$ are at least $\e_2$.
    In particular, we can find a flow box of the form $B_{\e_2}x$ containing $\w_a$ and $\w_b$ for some $x\in \Ocal$.
    Using the flow adapted coordinates on $B_\e x$, it follows from Proposition~\ref{prop:pseudo anosovs} \eqref{item:close} that
    \begin{equation}\label{eq:a and b close on unstable}
        |u(\w_a)-u(\w_b)| \leq 10^{-3}\e_2.
    \end{equation}

    Denote by $M_\bullet = \KZ{g_{\ell_\bullet}}{\w_\bullet}\in \SL$.
    Note that
    Proposition~\ref{prop:pseudo anosovs}(\ref{item:nesting}) implies that $\KZ{g}{\w_\bullet}$ is the identity matrix whenever $g\in B_{\e_2}$, $\bullet=a,b$.
    Moreover, we have $|u(\w_\bullet)- u(y)|$ is at most 
    $\e_2 e^{-2\ell_\bullet}/10^3$, 
    for all $y\in B(\bullet, \ell_\bullet,\e_2/10^3)$, $\bullet=a,b$.
    In particular, using the cocycle property, we see that for all $y \in B(\bullet, \ell_\bullet,\e_2/10^3)$,
    \begin{equation}\label{eq:cocycle on Bowen balls}
        \KZ{g_{\ell_\bullet}}{y} = M_\bullet.
    \end{equation}
    
    Recall that we are fixing $\w_1\in \Ocal$ to be a point with periodic horocycle orbit.
    For every $\t\geq 0$, $\bullet=a,b$, a subset $J\subseteq [0,1]$ of positive Lebesgue measure, we define $P_\bullet(\t,J)$ to be the set of points in $J$ whose orbit starts tracking $\w_\bullet$ at time $\t$.
    More precisely, we let
    \begin{equation}
        P_\bullet(\t,J) := \set{s\in J: g_\t u_s \w_1 \in B(\bullet,\ell_\bullet,\e_2/10^3) }.
    \end{equation}
    Given $s \in P_\bullet(\t,J)$ for some $\t>0$, 
    denote by $\Wcal(s)$ the local unstable leaf through $g_{\t+\ell_a} u_{s}\w_1$ inside the flow box $B_{\e_2/500} x$.
    We define $W(s)$ to be the pullback of $\Wcal(s)$ to our original horocycle through $\w_1$, i.e.
    \begin{equation}\label{eq:component definition}
        W(s) = \set{ s'\in J: g_{\t+\ell_a}u_{s'}\w_1 \in \Wcal(s)}.
    \end{equation}
    Our matching map $\mathbb{M}$ will send a component of the form $W(s)$ for $s\in P_a(\t,J)$ to a component $W(s')$ for a suitable $s'\in P_b(\t,J)$.

    Note that since $g_t$ contracts $u_s$ by $e^{-2|t|}$ for $t\leq 0$, we see that
    \begin{equation}\label{eq:length of a component}
        |W(s)| \leq \e_2 e^{-2(\t+\ell_a)}/500.
    \end{equation}

    In establishing the various properties of the matching map, we make use of several consequences of the equidistribution of translates of horocycle arcs of definite length by the geodesic flow.
    To this end, we introduce a thickening of the sets $P_\bullet(\t,J)$ in the $u$-direction, denoted by $P_\bullet^+(\t,J)$, and defined as follows:
    \begin{equation}
        P^+_\bullet(\t,J) = \bigcup_{s\in P_\bullet(\t,J)} W(s).
    \end{equation}
    The sets $P^+_b(\t,J)$ are defined analogously to $P_a^+(\t,J)$.
    Note that the local product structure on $B_{\e_2/500}x$ implies that for all $s_1,s_2\in P_a(\t,J)$, the sets $W(s_1)$ and $W(s_2)$ are either disjoint or equal.
    We shall refer to a subset of the form $W(s)$ as a \emph{component} of $P^+_\bullet(\t,J)$.
    We have the following estimate on the measures of $P_a^+(\t,J)$.

   \begin{lem}\label{lem:bound P_a^+ above}
    If $\e_2>0$ is sufficiently small, then the following holds.
    Suppose $J\subseteq [0,1]$ is a positive measure subset.
    Then, there exists $\t_0>0$, such that
    \begin{equation*}
        |P_a^+(\t,J)| \leq 40|P_a(\t,J)|,
    \end{equation*}
    for all $\t\geq \t_0$.
    \end{lem} 
    
    We will deduce Lemma~\ref{lem:bound P_a^+ above} from the following doubling property of Bowen balls.
    \begin{lem}\label{lem:doubling Bowen}
    For all sufficiently small $\rho>0$,
    \begin{equation*}
        \mu_\Ocal\left(B(a,\ell_a,\rho) \right) \leq 30\mu_\Ocal \left( B(a,\ell_a,\rho/2)\right).
    \end{equation*}
    \end{lem}

    \begin{proof}
        Recall that for $\d>0$ and $H$ a subgroup of $\SL$, we use $H_\d$ to denote the $\d$-neighborhood of identity.
        Then, for all $\d>0$, 
        \begin{equation}\label{eq:Bowen sandwich}
          U^-_{\d/3} A_{\d/3} U^+_{\d e^{-2\ell_a}/3} \cdot \w_a
          \subseteq B(a,\ell_a,\d) \subseteq U^-_{\d} A_{\d} U^+_{\d e^{-2\ell_a}}\cdot \w_a.
        \end{equation}
        The lemma then follows from the fact that the Haar measure is absolutely continuous with respect to the product of the Lebesgue measures on the groups $U^-, A,$ and $U^+$ and that the Jacobian of the product map $U^-\times A\times U^+ \r \SL$ is $1$ at identity. The factor $30$ arises from the fact that the radii of the factors in the product sets approximating $B(a,\ell_a,\d)$ in~\eqref{eq:Bowen sandwich} differ by a factor $1/3$ each, which scales the measure by $27$. Indeed, we can choose the multiplicative factor in our conclusion to be any number $\geq 27$.
        
    \end{proof}
    
    \begin{proof}[Proof of Lemma~\ref{lem:bound P_a^+ above}]
        Note that
        \begin{equation*}
            P_a^+(\t,J)\subseteq \set{s\in J: g_\t u_s\w_1 \in 
            B(a,\ell_a,\e_2/500)}.
        \end{equation*}
        Hence, Proposition~\ref{prop:horocycle equidist} implies that, if $\t$ is large enough, then 
        \[
        |P_a^+(\t,J)|\leq 1.1 \mu_\Ocal(B(a,\ell_a,\e_2/500))|J|.
        \]
        Similarly, for all sufficiently large $\t>0$, the measure of $P_a(\t,J)$ is at least $.9 \mu_\Ocal(B(a,\ell_a,10^{-3}\e_2))|J|$.
        The lemma then follows by Lemma~\ref{lem:doubling Bowen}.
    \end{proof}

    The following lemma allows us to make the matching procedure well-defined.    
    \begin{lem}\label{lem:disjoint bowen balls}
    Let $\e_0$ be the constant provided by Proposition~\ref{prop:pseudo anosovs}.
    If $\e$ is sufficiently small, depending only on $\e_0$, then
    the Bowen balls $B(a,\ell_a,10^{-3}\e_2)$ and $B(b,\ell_b,10^{-3}\e_2)$ are disjoint.
    \end{lem}
    \begin{proof}
    By Proposition~\ref{prop:pseudo anosovs}\eqref{item:injectivity}, the injectivity radii at $\w_a$ and $\w_b$ are at least $\e_0$.
    Hence, if these Bowen balls are not disjoint, using the cocycle property, one can show that there is a matrix $A$ in a uniformly bounded neighborhood of identity of $\SL$, depending only on $\e_0$, such that 
        \begin{equation*}
            \KZ{g_{\ell_a}}{\w_a} = \KZ{A}{\w_b} \KZ{g_{\ell_b}}{\w_b}.
        \end{equation*}
        In particular, there is a constant $C_0\geq 1$, depending only on $\e_0$, such that $\norm{\KZ{A}{\w_b}} \leq C_0$.
        Since $\norm{\KZ{g_{\ell_a}}{\w_a}} \leq \norm{\KZ{A}{\w_b}}\norm{ \KZ{g_{\ell_b}}{\w_b}} $,
        Proposition~\ref{prop:pseudo anosovs}~\eqref{item:norm discrepancy} implies
        that $\norm{\KZ{A}{\w_b}} \geq  \e^{-1}$.
        This is a contradiction as soon as $\e$ is sufficiently small.
    \end{proof}

    We record a useful estimate for later parts of the argument.
    In what follows, denote by $m_0$ the following constant:
    \begin{equation}\label{eq:m_0}
        m_0 = \min\set{\mu_\Ocal(B(\bullet, \ell_\bullet,\e_2/10^3)): \bullet =a,b}.
    \end{equation}
    By Proposition~\ref{prop:horocycle equidist} and a standard approximation argument, we can find $t_0>0$ so that for all $\t\geq t_0$, $\bullet=a,b$, we have
    \begin{equation}\label{eq:apply equidist}
        |P_\bullet(\t,J)| \geq .9\mu_\Ocal(B(\bullet,
        \ell_\bullet,\e_2/10^3))|J| \geq .9 m_0|J|.
    \end{equation}

    \subsection{Alignment of flags}
    Ultimately, we would like to decompose the cocycle $\KZ{g_t}{u_s\w_1}$, for suitably large $t$ into a product of large matrices corresponding to the value of the cocycle up to a time $\t$ at which the orbit of $u_s\w_1$ enters $B(a,\ell_a,\e_2/10^3)$, the value of the cocycle along the orbit segment shadowing $\w_a$, and finally the value of the cocycle from time $\t+\ell_a$ up to $t$.
    In order for the norm of this product to be comparable to the product of the norms of the indiviual matrices, we have to ensure the singular vectors of the relevant matrices are properly aligned.
    We do so by using the results of Section~\ref{sec:flags} to discard subsets of small measure of $P_a^+(\t,J)$ with poorly aligned flags from the matching procedure.

  To this end, recall the maps $\xi_{\mrm{in}},\xi_{\mrm{out}}$ defined in~\eqref{eq:contracting input def}. Let $\d>0$ be a small parameter to be chosen after Claim~\ref{claim:alignment before} below.
  Given $t>\t+\ell_a$, we let $P_a(t;\t,J)\subseteq P_a(\t,J)$ denote the subset of points $s$ in $P_a(\t,J)$ satisfying:
  \begin{align}\label{eq:avoid output flag}
       \xi_{\mrm{in}}(
        \KZ{g_{t-(\tau+\ell_a)}}
        {g_{\t+\ell_a}u_s\omega_1} )   
               \notin B_{\RP}(\xi_{\mrm{out}}(\KZ{g_{\t+\ell_a}}{u_s\w_1}),\d),
  \end{align}
  and
  \begin{align}\label{eq:avoid input flag}
    \xi_{\mrm{out}}\left(\KZ{g_\t}{u_s \w_1 }\right)
        \notin  B_{\RP}\left(\xi_{\mrm{in}}(M_a),  \d\right).      
  \end{align}
  Then, we note that if $s\in P_a(t;\t,J)$, then the first condition implies that
    \begin{equation}
     \label{eq:factor cocycle after shadowing 1}
        \norm{\KZ{g_t}{u_s\w_1}} \geq C_\d^{-1} \norm{\KZ{g_{t-\t-\ell_a}}{g_{\t+\ell_a}u_s\w_1}} \norm{\KZ{g_{\t+\ell_a}}{u_s\w_1}},
    \end{equation}
    for a constant $C_\d\geq 1$, depending only on $\d$.
  Moreover, the second condition implies that
    \begin{equation}
     \label{eq:factor cocycle before shadowing}
        \norm{\KZ{g_{\t+\ell_a}}{u_s\w_1}} \geq C_\d^{-1}
        \norm{M_a}\norm{\KZ{g_\t}{u_s\w_1}}.
    \end{equation}
    The goal for the rest of this subsection is to show that a definite proportion of the measure of $P_a^+(\t,J)$ satisfies the above conditions if $\d$ is small enough and $t$ is large enough (depending on $\d$).
    More precisely, we will show
    \begin{equation}\label{eq:measure without bad flags}
        |P_a(t;\t,J)| \geq \frac{1}{41} |P_a^+(\t,J)|
    \end{equation}

    \subsubsection{Post-shadowing flags}
    First, we show that Condition~\eqref{eq:avoid output flag} is satisfied on a set of definite measure.
    We do so by applying Proposition~\ref{prop:input nonatomic} to control the points with poorly aligned singular vectors after exiting our fixed Bowen balls.
    Recall the parameter $t_0$ chosen above~\eqref{eq:apply equidist}.
    By enlarging $t_0>0$ and using Proposition~\ref{prop:input nonatomic}, we can find $\d>0$ so that for all local unstable leaves $\Wcal \subset B_{\e_2/500}x$, $v\in\RP$, and all $\t\geq t_0$,
    \begin{equation}\label{eq:alignment after}
        \left|\set{w\in \Wcal: \xi_{\mrm{in}}\left(\KZ{g_\t}{w}\right) \in B_{\RP}(v,\d) } \right| \leq \e_2  |\Wcal|.
    \end{equation}
    Here we use the Lebesgue measure on $\Wcal$ viewed as an orbit segment of $u_s$.
    
    The next observation is that the cocycle $\KZ{g_{\t+\ell_a}}{u_s\w_1}$ is constant as $s$ varies in a given connected component $W$ of $P_a^+(\t,J)$.
    This follows from the cocycle property along with the fact that $\w_a$ is at distance at least $\e_2$ from the boundary of our chosen fundamental domain, by Proposition~\ref{prop:pseudo anosovs}(\ref{item:nesting}). 
    
    Fix a connected component $W$ of $P_a^+(\t,J)$ and let $\Wcal$ be the corresponding local unstable leaf.
    Let $t> \t+\ell_a+t_0$.
    Define $W^\diamond(t)\subseteq W$ as follows:
    \begin{align}
        W^\diamond(t) = 
        \left\{s\in W: \text{~\eqref{eq:avoid output flag} holds}
        \right \}.
    \end{align}
     
    As explained above, if $s\in W^\diamond(t)$ and $A\in \SL$ is the common value of the cocycle $\KZ{g_{\t+\ell_a}}{u_s\w_1}$ for $s\in W$, then we have the following refinement of~\eqref{eq:factor cocycle after shadowing 1}
    \begin{equation}\label{eq:factor cocycle after shadowing}
        \norm{\KZ{g_t}{u_s\w_1}} \geq C_\d^{-1} \norm{\KZ{g_{t-\t-\ell_a}}{g_{\t+\ell_a}u_s\w_1}} \norm{A}.
    \end{equation}
    Let $P_a^\diamond(t;\t,J)$ denote the union of the sets $W^\diamond(t)$ as $W$ varies over the connected components of $P^+_a(\t,J)$.
    The estimate in~\eqref{eq:alignment after} implies that
    \begin{equation}\label{eq:diamond measure}
        |P_a^\diamond (t;\t,J)| \geq (1- \e_2) |P^+_a(\t,J)|.
    \end{equation}

    \subsubsection{Pre-shadowing flags}
    We now turn our attention to Condition~\eqref{eq:avoid input flag}.
    To ensure proper alignment of the relevant flags at the time the orbits start shadowing $\w_a$, we have the following:
    \begin{claim}\label{claim:alignment before}
    There exists a uniform constant $C\geq 1$ so that the following holds. There exist $\d>0$ so that for all $J\subset [0,1]$, we can find $t_0>0$ such that for all $\t>t_0$ and $v\in \RP$, 
    \begin{equation}\label{eq:alignment before}
    \left|\set{s\in P_a^+(\t,J):  \xi_{\mrm{out}}\left(\KZ{g_\t}{u_s \w_1 }\right)        \in  B_{\RP}\left(v,  \d\right)}\right|
    \leq C \e_2  |P^+_a(\t,J)|.
    \end{equation}
    \end{claim}
    
    The above claim is an analogue of Proposition~\ref{prop:input nonatomic} for output flags. The deduction of this claim from Proposition~\ref{prop:input nonatomic} is somewhat involved and is postponed until the end of the section.

     Let $\d>0$ be sufficiently small and $t_0$ be sufficiently large so that~\eqref{eq:alignment after} and Claim~\ref{claim:alignment before} hold.
    For every $t >\t + \ell_a$, define
    \begin{equation*}
        P_a^+(t;\t,J) := \set{s\in P_a^\diamond(t;\t,J):  \text{~\eqref{eq:avoid input flag} holds} }.
    \end{equation*}
    Then,~\eqref{eq:diamond measure} and Claim~\ref{claim:alignment before} imply that for all $\t>t_0$ and $t>\t+\ell+t_0$,
    \begin{equation}
        |P_a^+(t;\t,J)| \geq (1- (C+1)\e_2) |P_a^+(\t,J)|.
    \end{equation}
    Let $P_a(t;\t,J)$ denote the intersection of $P_a^+(t;\t,J)$ with $P_a(\t,J)$.
    It follows by Lemma~\ref{lem:bound P_a^+ above} that
    \begin{equation}\label{eq:measure without bad flags 1}
        |P_a(t;\t,J)| \geq \left( \frac{1}{40}- (C+1)\e_2\right) |P_a^+(\t,J)| \geq \frac{1}{41} |P_a^+(\t,J)|,
    \end{equation}
    whenever $\e_2$ is small enough.
    This proves~\eqref{eq:measure without bad flags}.   
    
    \subsection{Measure estimates and discarded sets}
    Recall the constant $m_0$ defined in~\eqref{eq:m_0}.
    In what follows, we fix a parameter $N\in \N$ satisfying
    \begin{equation*}
        (1-.9 m_0)^{N+1} \leq \e_2.
    \end{equation*}
    We choose a sequence of times $0<\t_1<\cdots<\t_N$ and subsets $J_i\subset [0,1]$, $i=0,\dots,N$ inductively as follows.
    Set $J_0=[0,1]$ and $\t_0 =0$. If the subset $J_k$ has been chosen for $k\geq 0$, we let $\t_{k+1}>\t_k+1$ be sufficiently large so that~\eqref{eq:alignment after} and the conclusion of Claim~\ref{claim:alignment before} hold for $J=J_k$ and $\t= \t_{k+1}$.
    If $J_k$ and $\t_{k+1}$ have been defined, we let 
    \begin{equation*}
        J_{k+1} = J_k \setminus  P_a(\t_{k+1},J_k) .
    \end{equation*}
    If for some $1\leq k_0 \leq N-1$, the set $J_{k_0+1}$ has measure $\leq \e_2$, we stop the inductive procedure and set $N=k_0$. 
    Note that the above choices of $\t_k$ depend only on $\e_2$.
    For every $t > \ell_a +\t_N$, define\footnote{The sets $\Pcal_b(t)$ will be defined differently, cf.~\eqref{eq:Pcal_b}. In particular, we do not restrict the sets $P^+_b(\t,J)$ to subsets with properly aligned singular vectors. Roughly, the reason is we only need to bound \emph{from above} the norm of the product of cocycle factors for points in $\Pcal_b(t)$ by the product of their respective norms, which holds in general.}
   \begin{equation}
       \Pcal_a(t) = \bigcup_{k=1}^N P_a(t;\t_k,J_{k-1}), \qquad
       \Qcal^+_a(t) = \bigcup_{k=1}^N P^+_a(\t_k,J_{k-1}).
   \end{equation}
   Then, by construction, the members of the union defining $\Pcal_a(t)$ are disjoint\footnote{Recall that $P_a(t;\t_k,J_{k-1})$ is a subset of $P_a(\t_k,J_{k-1})$, which is in turn a subset of $J_{k-1}$.}.
   Hence, it follows by~\eqref{eq:measure without bad flags} that
   \begin{equation*}
       |\Pcal_a(t)| \geq \frac{1}{41} |\Qcal^+_a(t)|.
   \end{equation*}
   Moreover, by our choice of $N$, we have that $|\Qcal^+_a(t)|\geq 1-\e_2$.
   Hence, if $\e_2$ is small enough, we obtain
   \begin{equation}\label{eq:exhaust interval}
      |\Pcal_a(t)| \geq 1/45.
   \end{equation}

   Recall that we are fixing a unit vector $v_0\in \R^2$. By Proposition~\ref{prop:input nonatomic}, if $t$ is large enough and $\d$ is small enough, depending on $\e_2$, we have
   \begin{equation}\label{eq:separate inputs from v_0}
       \left| \set{s\in [0,1]: \xi_{\mrm{in}}\left(
       \KZ{g_t}{u_s\w_1}
       \right) 
       \in B_{\RP}(v_0,\d)
       } \right| \leq \e_2.
   \end{equation}
   Moreover, by Proposition~\ref{prop:horocycle recurrence}, applied with $\Kcal=\set{\w_1}$, we can find a compact set $\Omega\subset \Ocal$, such that
   \begin{equation}\label{eq:restrict to recurrent points}
       \left| \set{s\in[0,1]:g_\t u_s\w_1 \notin \Omega }\right| \leq\e_2, \qquad\forall  \t\geq 0.
   \end{equation}

   Define $\Pcal_a(t,v_0)$ to be the subset of $\Pcal_a(t)$ consisting of points who orbit lands in a compact set at time $t$ and for whom the input expanding singular direction of the cocycle is transverse to $v_0$. More precisely, we define $\Pcal_a(t,v_0)$ as follows:
   \begin{equation}
   \Pcal_a(t,v_0) := \set{ s\in \Pcal_a(t): g_t u_s\w_1 \in \Omega, 
       \xi_{\mrm{in}}\left( \KZ{g_t}{u_s\w_1}\right) 
       \notin B_{\RP}(v_0,\d)}.
   \end{equation}
   By~\eqref{eq:exhaust interval},~\eqref{eq:separate inputs from v_0}, and~\eqref{eq:restrict to recurrent points}, if $\e_2$ is small enough, we have the lower bound
   \begin{equation}\label{eq:exhaust interval with v_0}
       |\Pcal_a(t,v_0)| \geq 1/49,
   \end{equation}
   The advantage of introducing the set $\Pcal_a(t,v_0)$ is that for all $s\in \Pcal_a(t,v_0)$, we have that
   \begin{equation}\label{eq:vector grows like norm}
       \norm{\KZ{g_t}{u_s\w_1}v_0} \geq C_\d^{-1} \norm{\KZ{g_t}{u_s\w_1}}.
   \end{equation}

    \subsection{Matching and proof of Proposition~\ref{prop:partitions spread}}
    Define sets $\Qcal^+_b(t)$ analogously to $\Qcal^+_a(t)$ as follows:
    \begin{equation*}
        \Qcal_b^+(t) = \bigcup_{k=1}^N P_b^+(\t_k,J_{k-1}).
    \end{equation*}
    Recall the definition of the components of $P_\bullet^+(\t,J)$ of the form $W(s)$ introduced in~\eqref{eq:component definition}. 
    We wish to define a matching map $\mathbb{M}:\Pcal_a(t,v_0) \r \Qcal^+_b(t)$.
    We require $\mathbb{M}$ to be the restriction of a map sending each component $W(s)$ of $P^+_a(\t,J)$ to a chosen component of $P^+_b(\t,J)$ via the weak stable holonomy inside $B_{\e_2/500}x$.
    Hence, our first step is to pair the components of each set.
    
    Observe that it follows from~\eqref{eq:a and b close on unstable} and the discussion following it that for every $s_a\in P_a(\t,J)$, there is $s_b \in P_b(\t,J)$ such that 
    \begin{equation}\label{eq:paired components are close}
        |s_a-s_b| \leq\e_2 e^{-2\t}/10^3.
    \end{equation}
    The restriction of the matching map $\mathbb{M}$ to $W(s_a)$ is given by the weak stable holonomy sending $W(s_a)$ to $W(s_b)$ as we detail below. In fact, the arguments below can be used to show that there is a unique such choice of components $W(s_b)$.
    
    Next, we claim that the above pairing of components is one-to-one, i.e. if $W(s_{a_1})$ and $W(s_{a_2})$ are distinct components of $P^+(\t,J)$ then the corresponding components $W(s_{b_1})$ and $W(s_{b_2})$ of $P^+_b(\t,J)$ chosen as above are also distinct.
    This follows as soon as we show that the distance between $W(s_{a_1})$ and $W(s_{a_2})$ is at least $\e_0 e^{-2\t}/10^3\gg \e_2e^{-2\t}/10^3$, whenever $s_{a_1}, s_{a_2}\in P_a(\t,J)$ and $W(s_{a_1})\neq W(s_{a_2}) $.
    This latter separation follows from the fact that the injectivity radius at $\w_a$ is at least $\e_0$ so that any two connected components of the set
    \begin{equation*}
        \set{s\in[0,1]: g_{t'}u_s\w_1 \in B_{\e_2/10^3}\w_a }
    \end{equation*}
    are separated by at least $\e_0 e^{-2t'}$ for any $t'\geq 0$.

    Hence, we can define a map $\Psi_\t^J: P_a^+(\t,J)\r P_b^+(\t,J)$ as follows.
    Let $W$ be a given connected component of $P_a^+(\t,J)$
    and let $W'$ be the connected component of $P_b^+(\t,J)$ that we paired to $W$. Let $y = g_\t u_s\w_1$ and $z= g_\t u_{s'}\w_1$ for some $s\in W$ and $s'\in W'$.
    Recall the definition of the weak stable holonmy map $\Psi_{y,z}^{\mrm{cs}}$ (inside $B_{\e/500}x$) in Section~\ref{sec:flow charts}.
    Define $\tilde{\Psi}_\t^J := g_{-\t-\ell_a} \circ \Psi_{y,z}^{\mrm{cs}} \circ g_{\t+\ell_a} $.
    Then, $\tilde{\Psi}_\t^J$ is a map between the two segments of the closed horocycle of $\w_1$ corresponding to $W$ and $W'$.
    We define $\Psi_\t^J$ by composing $\tilde{\Psi}_\t^J$ with the identification of the closed horocycle with the interval $[0,1)$ via $s\mapsto u_s\w_1$. Note that the Jacobians of $\Psi_\t^J$ and $\Psi_{y,z}^{\mrm{cs}}$ agree.

    We define $\mathbb{M}:\Pcal_a(t)\r\Qcal_b^+(t)$ by:
    \begin{equation*}
        \mathbb{M}(s) = \Psi_{\t_k}^{J_{k-1}}(s), \qquad \text{ if } s\in P_a(t;\t_k,J_{k-1}).
    \end{equation*}
    The map $\mathbb{M}$ is well-defined since $\Pcal_a(t)$ is a disjoint union of sets of the form $P_a(t;\t_k,J_{k-1})$.
    We let 
    \begin{equation}\label{eq:Pcal_b}
        \Pcal_b(t) := \mathbb{M}(\Pcal_a(t)), \qquad \Pcal_b(t,v_0) := \mathbb{M}(\Pcal_a(t,v_0)).
    \end{equation}

    If $\e_2$ is small enough, the Jacobian of the individual maps $\Psi_\t^J$ become uniformly close to identity (cf.~Section~\ref{sec:flow charts}).
    It follows that the same holds for $\mathbb{M}$.
    In particular, if $\e_2$ is sufficiently small, then~\eqref{eq:exhaust interval} implies that
    \begin{equation*}
        |\Pcal_b(t,v_0)| \geq 1/50.
    \end{equation*}
    Moreover, by Lemma~\ref{lem:disjoint bowen balls}, the sets $\Pcal_a(t)$ and $\Pcal_b(t)$ are disjoint.
    Hence, the same holds for the sets $\Pcal_a(t,v_0)$ and $\Pcal_b(t,v_0)$.
    This completes the verification of Items~\eqref{item:strategy1} and~\eqref{item:strategy2} of the strategy outlined at the begining of the section.
    
    To verify Item~\eqref{item:strategy3}, let $ \Pcal_a(t,v_0)$. Then, $s\in P_a(t;\t_k,J_{k-1})$ for some $k$.
    We claim that
    \begin{equation}\label{eq:matching discrepancy}
        \norm{\KZ{g_t}{u_{\mathbb{M}(s)}\w_1}v_0}\leq \tilde{C}\e  \norm{\KZ{g_t}{u_s\w_1}v_0},
    \end{equation}
    for some constant $\tilde{C}\geq 1$, depending only on $\e_2$. Note that this claim finishes the proof by choosing $\e$ to be small enough so that $\tilde{C}\e < \k^2$.

    To show~\eqref{eq:matching discrepancy}, let $C_\d \geq 1$ be chosen large enough, depending only on $\d$ to satisfy~\eqref{eq:factor cocycle after shadowing},~\eqref{eq:factor cocycle before shadowing}, and~\eqref{eq:vector grows like norm}.
    Recall that $M_\bullet = \KZ{g_{\ell_\bullet}}{\w_\bullet}$.
    By~\eqref{eq:factor cocycle after shadowing},~\eqref{eq:factor cocycle before shadowing} and~\eqref{eq:vector grows like norm}, since $\norm{v_0}=1$, we have
    \begin{align}\label{eq:factoring}
    &\frac{\norm{\KZ{g_t}{u_{\mathbb{M}(s)}\w_1}v_0}} {\norm{\KZ{g_t}{u_s\w_1}v_0}}
     \leq C_\d
        \frac{\norm{\KZ{g_t}{u_{\mathbb{M}(s)}\w_1}}}{\norm{\KZ{g_t}{u_s\w_1}}} 
        \nonumber\\
        &\leq C_\d^3 
        \frac{
        \norm{\KZ{g_{t-(\t_k+\ell_b)}}{g_{\t_k+\ell_b}u_{\mathbb{M}(s)}\w_1}
        }}
        {
        \norm{\KZ{g_{t-(\t_k+\ell_a)}}{g_{\t_k+\ell_a}u_s\w_1}}
        }
        \cdot \frac{\norm{M_b}}{\norm{M_a}} 
       \cdot \frac{
        \norm{\KZ{g_{\t_k}}{u_{\mathbb{M}(s)}\w_1}} 
        }
        {\norm{\KZ{g_{\t_k}}{u_s\w_1}} }
        \nonumber\\
        &\leq C_\d^3  \cdot
        \frac{
        \norm{\KZ{g_{t-(\t_k+\ell_b)}}{g_{\t_k+\ell_b}u_{\mathbb{M}(s)}\w_1}
        }}
        {
        \norm{\KZ{g_{t-(\t_k+\ell_a)}}{g_{\t_k+\ell_a}u_s\w_1}}
        }
        \cdot \e\cdot \frac{
        \norm{\KZ{g_{\t_k}}{u_{\mathbb{M}(s)}\w_1}} 
        }
        {\norm{\KZ{g_{\t_k}}{u_s\w_1}} }.
    \end{align} 
    In the last inequality, we used the fact that $\w_a$ and $\w_b$ are chosen to satisfy Proposition~\ref{prop:pseudo anosovs}\eqref{item:norm discrepancy}.
    
    Let $r= \mathbb{M}(s) - s$. Then, the cocycle property implies that
    \begin{equation*}
        \KZ{g_{\t_k}}{u_{\mathbb{M}(s)}\w_1}
        = \KZ{u_{e^{2\t_k}r}}{g_{\t_k}u_s\w_1} \KZ{g_{\t_k}}{u_s\w_1}
        \KZ{u_r}{u_s\w_1}.
    \end{equation*} 
    Moreover, the estimates in~\eqref{eq:length of a component} and~\eqref{eq:paired components are close} imply that $|r|\leq \e_2 e^{-2\t_k}$.
    Note further that $g_{\t_k}u_s\w_1$ belongs to $B_{\e_2/500}x$ by our choice of $\t_k$.
    In particular, the injectivity radii at $u_s\w_1$ and $g_{\t_k}u_s\w_1$ are uniformly bounded from below, depending only on the constant $\e_0$ given by Proposition~\ref{prop:pseudo anosovs} (since $\w_a,\w_b\in B_{\e/500}x$ and their injectivity radii are chosen to be larger than $\e_0\gg \e_2$ by Proposition~\ref{prop:pseudo anosovs}~\eqref{item:injectivity}).
    It follows that there is a uniform constant $C_0\geq 1$ such that
    \begin{equation}\label{eq:bounded difference before shadowing}
        \norm{\KZ{g_{\t_k}}{u_{\mathbb{M}(s)}\w_1}} \leq
        C_0 \norm{ \KZ{g_{\t_k}}{u_s\w_1}}.
    \end{equation}
    
    To control the ratio of the remaining terms, note that 
    \begin{equation*}
        g_{\t_k+\ell_b}u_{\mathbb{M}(s)}\w_1 = 
        \Psi_{y,z}^{\mrm{cs}}( g_{\t_k+\ell_a}u_{s}\w_1) = p^-   g_{\t_k+\ell_a}u_{s}\w_1,
    \end{equation*}
    for some $p^-\in A_{\e_2/500}U^-_{\e_2/500}\subset \SL$, and for a suitable choice of $y$ and $z$. Indeed, this follows from definitions of the local stable holonomy map and the map $\mathbb{M}$. By the cocycle property, we get
    \begin{align*}
        \KZ{g_{t-(\t_k+\ell_b)}}{g_{\t_k+\ell_b}u_{\mathbb{M}(s)}\w_1}
        = \KZ{g_{\ell_a-\ell_b}p'}{g_tu_s\w_1} 
        \KZ{g_{t-(\t_k+\ell_a )}}{ g_{\t_k+\ell_a}u_{s}\w_1}
        ,
    \end{align*}
    where $p' = g_{t-(\t_k+\ell_a)}p^- g_{-t + \t_k+\ell_a}$ and we used the fact that $\KZ{p^-}{ g_{\t_k+\ell_a}u_{s}\w_1}$ is the identity matrix since $B_{\e_2/500}x$ is in the interior of our chosen fundamental domain of $\Ocal$ for defining the cocycle.
    This is ensured by Proposition~\ref{prop:pseudo anosovs}~\eqref{item:nesting}.
    Next, we claim that 
    \begin{equation}\label{eq:bounded difference after shadowing}
        \norm{\KZ{g_{\ell_a-\ell_b}p'}{g_tu_s\w_1}} \leq C_1,
    \end{equation}
    for a constant $C_1\geq 1$, depending only $\e_2$.
    Indeed, by Proposition~\ref{prop:pseudo anosovs}~\eqref{item:close periods}, $|\ell_b-\ell_a| <1$. Hence, since $g_t$ contracts $U^-$, 
    \[g_{\ell_a-\ell_b}p' \in A_{2}U^-_{\e_2 e^{-2(t-(\t_k+\ell_a))}/500}.\]
    Moreover, by definition of $\Pcal_a(t,v_0)$, $g_tu_s\w_1 \in \Omega$, where $\Omega$ is the compact set chosen in~\eqref{eq:restrict to recurrent points}.
    This proves~\eqref{eq:bounded difference after shadowing}.
    
    Hence, the claimed estimate in~\eqref{eq:matching discrepancy} follows by~\eqref{eq:factoring},~\eqref{eq:bounded difference before shadowing}, and~\eqref{eq:bounded difference after shadowing} by taking $\tilde{C} = C_\d^3C_0 C_1$.
    Note that the constant $C_\d$ depends on $\e_2$ and not $\e$.
    To complete the proof of the proposition, it remains to verify Claim~\ref{claim:alignment before} which we do in the next subsection.

    \subsection{Proof of Claim~\ref{claim:alignment before}}
    The cocycle property implies that for all $y\in \Ocal$, $\mrm{KZ}(g_\t,y) = \mrm{KZ}(g_{-\t},g_{\t}y)^{-1}$.
    It follows that 
    \begin{equation}\label{eq:in =out}
        \xi_{\mrm{out}}\left(\KZ{g_\t}{u_s \w_1 }\right)=\xi_{\mrm{in}}\left(\KZ{g_{-\t}}{g_{\t}u_s \w_1 }\right).
    \end{equation}
    Moreover, note that output flags of $\KZ{g_\t}{u_s \w_1}$ do not change along local strong unstable leaves inside $B_{\e_2/500}x$.
    More precisely, suppose $s\in [0,1]$ satisfies $y_1:=g_\t u_s\w_1 \in B_{\e_2/500}x $ and let $r\in [0,1]$ denote any other point such that $g_\t u_r\w_1$ lies on the same local unstable leaf through $y_1$ in $B_{\e_2/500}x$.
    Then, $|r-s|\leq e^{-2\t}\e_2/500$.
    Hence, the cocycle property and  the fact that $B_{\e_2/500}x$ is an isometrically embedded neighborhood of identity in $\SL$ imply that $\xi_{\mrm{out}}\left(\KZ{g_\t}{u_s\w_1}\right) = \xi_{\mrm{out}}\left(\KZ{g_\t}{u_r\w_1}\right)$.
    
    In particular, to prove the claim, we count the number of points of intersection $w$ of $g_\t u_s\w_1, s\in [0,1]$ with a weak stable leaf in $B_{\e_2/500}x$ for which $\xi_{\mrm{in}}\left(\KZ{g_{-\t}}{w}\right)$ belongs to the ball $B_{\RP}(v,\d) $, for a suitably chosen $\d$.
    
    By Proposition~\ref{prop:input nonatomic}, applied to the backward geodesic flow, we find $\d>0$ so that if $t_0$ is large enough, then for all $\t\geq t_0$, all $v\in\RP$ and all local stable leaves $\Wcal^-\subset B_{\e_2/500}x$,
    \begin{equation*}
        \left|\set{w\in \Wcal^-: \xi_{\mrm{in}}\left(\KZ{g_{-\t}}{w}\right) \in B_{\RP}(v,\d) } \right| \leq  \e_2 m_0 |J|  |\Wcal^-|,
    \end{equation*}
    where $m_0$ is defined in~\eqref{eq:m_0}.
    Note that since $m_0$ depends only on $\e_2$, we have that $t_0$ depends only on $J$ and $\e_2$.

    Given a weak stable leaf $\Wcal^{0-}\subset B_{\e_2/500}x$, we endow it with the product Lebesgue measure via the local coordinates introduced in Section~\ref{sec:flow charts}.
    We also use $|\cdot|$ to denote such a measure.
    By Fubini's theorem and the product structure of the measure on $\Wcal^{0-}$, it follows that for any local weak stable leaf $\Wcal^{0-}\subset B_{\e_2/500}x$, we have
    \begin{equation}\label{eq:alignment weak stable leaf}
        \left|\set{w\in \Wcal^{0-}: \xi_{\mrm{in}}\left(\KZ{g_{-\t}}{w}\right) \in B_{\RP}(v,\d) } \right| \leq  \e_2 m_0 |J|  |\Wcal^{0-}|.
    \end{equation}

    Fix some weak stable leaf $\Wcal^{0-}\subset B_{\e_2/500}x$ which meets every positive measure connected component of the intersection of $\{s\in J:g_\t u_s\w_1\}$ with $B_{\e_2/500}x$. Define for each $\t>0$ the following finite sets:
    \begin{equation*}
        \Wcal^{0-}(\t) := \set{w\in\Wcal^{0-}: w = g_\t u_s \w_1 \text{ for some } s\in J  }.
    \end{equation*} 
    To apply the estimate in~\eqref{eq:alignment weak stable leaf}, we first note that 
    \begin{align}\label{eq:estimate via counting}
      |\{s\in &P_a^+(\t,J): \xi_{\mrm{out}}\left(\KZ{g_{\t}}{u_s \w_1}\right) \in B_{\RP}(v,\d/2) \}|
      \nonumber\\
        &\leq  (\e_2/500)\cdot  e^{-2\t}
        \cdot \# \set{h\in \Wcal^{0-}(\t): \xi_{\mrm{in}}\left(\KZ{g_{-\t}}{h}\right) \in B_{\RP}(v,\d/2)   }.
    \end{align} 
    Indeed, this follows from the relation between input and output flags in~\eqref{eq:in =out} and the fact that the flag $\xi_{\mrm{out}}\left(\KZ{g_{\t}}{u_s \w_1}\right)$ is constant for points $s\in J$ for which $g_\t u_s \w_1$ belong to the same local unstable leaf in $B_{\e_2/500}x$.
    The measure of the set of such points is at most $(\e_2/500)\cdot  e^{-2\t}$ since $g_\t$ expands $u_s$ by $e^{2\t}$.

    It remains to estimate the right side of~\eqref{eq:estimate via counting} using~\eqref{eq:alignment weak stable leaf}. 
    The key observation is that if $h\in \Wcal^{0-}$ satisfies $\xi_{\mrm{in}}\left(\KZ{g_{-\t}}{h}\right) \in B_{\RP}(v,\d/2) $ and $\t$ is large enough, then Lemma~\ref{lem:singular product} implies that the same holds for a whole neighborhood of $h$ in $\Wcal^{0-}$ of radius $e^{-2\t}$ along the stable coordinate, with $\d$ in place of $\d/2$.
    It follows that
    \begin{align}\label{eq:bound count by measure}
        (\e_2/500)\cdot  e^{-2\t}
        \cdot \# 
        \{h\in &\Wcal^{0-}(\t):  \xi_{\mrm{in}}\left(\KZ{g_{-\t}}{h}\right) \in B_{\RP}(v,\d/2)   \} 
        \nonumber\\
        &\leq
        C_\#
        \left|\set{w\in \Wcal^{0-}: \xi_{\mrm{in}}\left(\KZ{g_{-\t}}{w}\right) \in B_{\RP}(v,\d) } \right|,
    \end{align}
    where $C_\#\geq 1$ is the maximal cardinality of the set of points in $\Wcal^{0-}(\t)$ whose stable coordinates all lie within a fixed interval of size $e^{-2\t}$.
    
     We claim that
    \begin{equation}\label{eq:bound C sharp}
        C_\# \text{ is uniformly bounded independently of all parameters}.
    \end{equation}
    In view of~\eqref{eq:apply equidist}, we have $m_0|J|\leq |P_a^+(a,J)|/0.9$. Hence, combining~\eqref{eq:estimate via counting} and~\eqref{eq:bound count by measure} with the estimates~\eqref{eq:alignment weak stable leaf} and~\eqref{eq:bound C sharp} completes the proof of Claim~\ref{claim:alignment before}. 
    
    To prove~\eqref{eq:bound C sharp}, we need to understand the spacing of the points in $\Wcal^{0-}(\t)$.
    Roughly, we show that the stable coordinates of the points in $\Wcal^{0-}(\t)$ must be separated by at least $e^{-2\t}$. 
    More precisely, fix some arbitrary interval $V$ of size $e^{-2\t}$ and let $\mf{S}(V)\subset J$ be the set consisting of points $s\in J$ with $g_\t u_s \w_1 \in \Wcal^{0-}(\t)$ such that the stable coordinates of 
    $\set{g_\t u_s\w_1:s\in\mf{S}(V)}$ all lie in $V$.
    We will show that $\mf{S}(V)$ has uniformly bounded cardinality. 
    
    For each $\rho>0$, let $\mc{U}_\rho$ denote the neighborhood of identity in $\SL$ of radius $\rho$. Write $\w_1 = g\G$ and suppose $g_\t u_{s_1} g\g = h g_\t u_{s_2} g$ for some $h\in AU^-$, $\g\in \G$, and $s_1\neq s_2 \in [0,1]$ with $g_\t u_{s_i}\w_1\in \Wcal^{0-}(\t)$ for $i=1,2$.
    In particular, $h= g_v \hat{u}_r $, with $|v|\leq \e_2/500 $ and $|r|\leq e^{-2\t}$.
    Then, $g_{-\t}hg_\t\in \mc{U}_2$. Hence, since $g_{-\t}hg_\t =  u_{s_1} g \g g^{-1} u_{s_1}^{-1}\cdot  u_{s_1-s_2}$, we get
    \begin{equation}\label{eq:contract Gamma}
         u_{s_1} g \g g^{-1} u_{s_1}^{-1} \in \mc{U}_2\cdot  u_{s_1-s_2}^{-1} \subseteq \mc{U}_{3}.
    \end{equation}
    By discreteness of the lattice $g\G g^{-1}$, it follows that the set
    \begin{equation*}
        \mf{T} :=  \set{\g \in \G: \exists s\in [0,1]; u_s g\g g^{-1} u_s^{-1} \in \mc{U}_3 }
    \end{equation*}
    is finite and has a uniformly bounded cardinality independently of all the parameters.
    We show that 
    \begin{equation}\label{eq:bound cardinality of fancy S}
        \#\mf{S}(V)\leq \#\mf{T}.
    \end{equation}
    Since $C_\#=\sup_V \mf{S}(V)$ and $V$ was arbitrary, this proves~\eqref{eq:bound C sharp} and completes the proof of the claim. 
    
    To prove~\eqref{eq:bound cardinality of fancy S}, fix some point $x_0:=g_{\t}u_{s_0}\w_1\in \mc{W}^{0-}(\t)$ so that the stable coordinates of 
    $\set{g_\t u_s\w_1:s\in\mf{S}(V)}$ all lie in an interval of size $e^{-2\t}$ around the stable coordinate of $x_0$.
    Define a map $\Phi:\mf{S}(V)\to \mf{T}$ which assigns to each $s$ an element $\g\in \G$ such that 
    \begin{equation}\label{eq:the counting map}
        hg_\t u_s g\g = g_{\t}u_{s_0}g
    \end{equation}
     for some $h=g_{v}\hat{u}_r\in AU^-$ with $|v|\leq \e_2/500$ and $|r|\leq e^{-2\t}$.
    By the argument preceding~\eqref{eq:contract Gamma}, we note that such $\g$ must belong to $\mf{T}$.
    Moreover, recalling that the injectivity radius in the box $B_{\e_2/500}x$ is $\gg \e_0$ by Proposition~\ref{prop:pseudo anosovs}(\ref{item:injectivity}), there is a unique such $\g\in \G$ satisfying~\eqref{eq:the counting map} so that $\Phi$ is well-defined.
    Indeed, this can always be arranged by ensuring $\e$ is small enough and $\t$ is large enough, depending only on $\e_0$.

    Note that~\eqref{eq:bound cardinality of fancy S} follows if $\Phi$ is injective.
    Suppose not. Then, we can find $s_1\neq  s_2 \in \mf{S}(V)$ so that 
    \begin{equation*}
        h_1 g_\t u_{s_1}  = h_2 g_\t u_{s_2} ,
    \end{equation*}
    for some $h_1,h_2 \in AU^-$.
    This implies that $g_{-\t} h_2^{-1}h_1g_\t = u_{s_2-s_1}$ and, hence, that the groups $AU^-$ and $U^+$ intersect non-trivially, since $s_1\neq s_2$. This is a contradiction. This shows that $\Phi$ is injective and completes the proof of~\eqref{eq:bound cardinality of fancy S}.

    
    \section{Concluding remarks}
    We conclude the paper with two remarks.
We briefly outline the mechanism for our non-genericity result (Corollary \ref{cor:nongen}), which perhaps can be applied more generally. There is an action of $\SL$ on our space and a closed, $\SL$-invariant sublocus $Y$. The derivative cocycle for the geodesic flow (built from the KZ cocycle, cf.~Section~\ref{eq:G derivative}) has oscillations in a direction transverse to $Y$ along a horocycle segment. Considering points moved transversly to $Y$ from $g_t$ pushes of our fixed horocycle,
we use the fluctuations (and some auxiliary results, especially Section \ref{sec:near teich}) to obtain longer and denser horocycle segments, where the uniform measure on the segment is a definite distance from a fixed fully supported ergodic measure.

If we start with a closed horocycle in $Y$, one can hope to obtain an ``exotic" limit of closed horocycles if we can deform the horocycles transverse to $Y$ while keeping them closed.

In the case of strata of translation surfaces, we hope that our approach can be used to exhibit similar behavior to the one established here for $\HH(2)$. Indeed, we only make specific use of the octagon locus for producing pseudo-Anosovs with distinct Lyapunov exponents (Proposition \ref{prop:pseudo anosovs}) and for proving strong irreducibility of the balanced subbundle (Proposition~\ref{prop:strong irreducibility}). In particular, we hope that the renormalization arguments 
in Lemma \ref{lem:renorm} and Section \ref{sec:general} and the linearization of certain affine geodesics near certain $\mrm{GL}_2(\mathbb{R})$-orbit closures in Section \ref{sec:near teich} can be used to show more behaviors for the horocycle flow.

    \bibliography{bibliography}{} 
\bibliographystyle{amsalpha}

\end{document}